\documentclass[12pt,twoside]{amsart}
\usepackage{amssymb}
\usepackage{amscd}
\usepackage[abbrev,alphabetic]{amsrefs}
\usepackage{hyperref}
\usepackage{comment}
\usepackage{tikz}
\usetikzlibrary{cd}
\usepackage{here}
\usepackage{color}

\usepackage[margin=1.25in]{geometry}

\title[Bertini theorems admitting base changes]
{Bertini theorems admitting base changes} 
\author{Hiromu Tanaka} 
\address{Graduate School of Mathematical Sciences, 
The University of Tokyo, 
3-8-1 Komaba, Meguro-ku, Tokyo 153-8914, JAPAN} 
\email{tanaka@ms.u-tokyo.ac.jp}
\newcommand{\pr}[0]{{\operatorname{pr}}}

\newcommand{\red}[0]{{\operatorname{red}}}

\newcommand{\Ker}[0]{{\operatorname{Ker}}}

\newcommand{\Proj}[0]{{\operatorname{Proj}}}
\newcommand{\Spec}[0]{{\operatorname{Spec}}}
\newcommand{\Spm}[0]{{\operatorname{Spm}}}

\newcommand{\Bs}[0]{{\operatorname{Bs}}}
\newcommand{\Supp}[0]{{\operatorname{Supp}}}
\newcommand{\Pic}[0]{{\operatorname{Pic}}}

\newcommand{\univ}[0]{{\operatorname{univ}}}
\newcommand{\gen}[0]{{\operatorname{gen}}}
\newcommand{\Frac}[0]{{\operatorname{Frac}}}
\newtheorem{thm}{Theorem}[section]
\newtheorem{lem}[thm]{Lemma}
\newtheorem{cor}[thm]{Corollary}
\newtheorem{prop}[thm]{Proposition}

\newtheorem{step}{Step}

\theoremstyle{definition}
\newtheorem{ex}[thm]{Example}

\newtheorem{dfn}[thm]{Definition}

\newtheorem{rem}[thm]{Remark}
      
\newtheorem{nota}[thm]{Notation}         

\makeatletter
  
  \@addtoreset{equation}{thm}
  \makeatother

\newcommand{\cred}{\color{red}}

\newcommand{\MO}{\mathcal{O}}
\newcommand{\C}{\mathbb{C}}
\newcommand{\R}{\mathbb{R}}
\newcommand{\Q}{\mathbb{Q}}
\newcommand{\Z}{\mathbb{Z}}
\newcommand{\F}{\mathbb{F}}
\renewcommand{\P}{\mathbb{P}}
\newcommand{\p}{\mathfrak{p}}
\newcommand{\mfq}{\mathfrak{q}}
\newcommand{\m}{\mathfrak{m}}
\newcommand{\n}{\mathfrak{n}}
\begin{document}

\maketitle

\begin{abstract}
Given a base point free linear system on an algebraic variety, many classes of singularities are stable under taking suitable members after enlarging the base field. We establish analogous results when the base ring is  an excellent ring. 
\end{abstract}

\tableofcontents

\section{Introduction}

Given a smooth projective variety $X$ over an infinite field, 
the classical Bertini theorem asserts that general hyperplane sections of $X$ are smooth \cite[Ch. II, Theorem 8.18]{Har77}. 
In other words, if $L$ is a very ample divisor on $X$, then general members of its complete linear system are smooth. 
In characteristic zero, the same conclusion holds for arbitrary base point free linear systems \cite[Ch. III, Corollary 10.9]{Har77}, 
whilst such a generalisation is known to fail in positive characteristic (cf. Example \ref{e-Frob} and Example \ref{e-fib-P1}). 
This pathological phenomenon is one of typical obstructions to 
extend results in characteristic zero to ones in positive characteristic. 

The purpose of this paper is to establish results of Bertini type for base point free linear systems in positive and mixed characteristic. 
For a base point free linear system,  all the members can have bad singularities. 
On the other hand, 
we shall prove that the linear system will have good members after taking a suitable base change. 
To this end, we first introduce purely transcendental extensions $R(t_1, ..., t_n)$ 
for arbitrary rings $R$: 
\[
R(t_1, ..., t_n) := U_{R[t_1, ..., t_n]/R}^{-1}(R [t_1, ..., t_n]), 
\]
where $R[t_1, ..., t_n]$ denotes the polynomial ring over $R$ and 
\[
U_{R[t_1, ..., t_n]/R} := \bigcap_{\p \in \Spec\,R} (R [t_1, ..., t_n] \setminus \p R[t_1, ..., t_n]).
\]
When $R$ is a field, $R(t_1, ..., t_n)$ coincides with the purely transcendental extension in the usual sense. 
After taking the base change $X \mapsto X \times_R R(t_1, ..., t_n)$, 
base point free linear systems have good members as follows.

\begin{thm}[Theorem \ref{t-tr-ext-dim}, Theorem \ref{t-Bertini1}]\label{intro-Bertini1}
Fix $n \in \Z_{>0}$. 
Let $R$ be a noetherian ring and set $R(n) := R(t_1, ..., t_n)$. 
Then there exists a hyperplane $(\mathbb P^n_R)^{\gen}_{{\rm id}}$ of $\mathbb P^n_{R(n)}$ such that the following holds.

Let $\varphi : X \to \mathbb P^n_R$ be an $R$-morphism from a scheme $X$ of finite type over $R$. 
Set $X_{R(n)} := X \times_R R(n)$. 
Let (P) be a property for noetherian schemes which satisfies the following properties (I)--(III) 
(e.g. (P) can be regular, reduced, klt, or strongly $F$-regular). 
\begin{enumerate}
\renewcommand{\labelenumi}{(\Roman{enumi})}
\item 
Let $Y$ be a noetherian scheme and let $Y = \bigcup_{i \in I} Y_i$ be an open cover. 
Then $Y$ is (P) if and only if $Y_i$ is (P) for 
{every} $i \in I$. 
\item 
For a noetherian ring $A$, 
if $\Spec\,A$ is (P), then also $\Spec\,A[t]$ is (P), 
where $A[t]$ denotes the polynomial ring over $A$ with one variable. 
\item 
For a noetherian ring $A$ and a multiplicatively closed subset $S$ of $A$, 
if $\Spec\,A$ is (P), then also $\Spec\,S^{-1}A$ is (P). 
\end{enumerate}
If $X$ is (P), then also $X^{\gen}_{\varphi}$ is (P), where 
$X^{\gen}_{\varphi}$ is the scheme-theoretic inverse image of $(\mathbb P^n_R)^{\gen}_{{\rm id}}$, 
i.e. the following diagram is cartesian: 
\[
\begin{CD}
X^{\gen}_{\varphi} @>>> X_{R(n)}\\
@VVV @VV\varphi \times_R R(n) V\\
(\mathbb P^n_R)^{\gen}_{{\rm id}} @>>> \mathbb P^n_{R(n)}. 
\end{CD}
\]
\end{thm}

\begin{rem}
We shall establish analogous results for pairs $(X, \Delta)$ (Theorem \ref{t-Bertini3.1}, Theorem \ref{t-Bertini3.2}). 
\end{rem}

\begin{rem}
Theorem \ref{intro-Bertini1} is folklore for the case when $R$ is a field. 
Indeed, we can apply a similar argument to the one of \cite[Theorem 1]{CGM86} for this case. 
\end{rem}

Theorem \ref{intro-Bertini1} assures the existence of good members $X^{\gen}_{\varphi}$. 
If $k$ is an infinite field, then the classical Bertini theorem states that 
general hyperplane sections are good. 
We shall establish the following analogous result. 

\begin{thm}[Theorem \ref{t-dense-good}, Remark \ref{r-dense-good}]\label{intro-Bertini2}
Let $k$ be a field. 
Fix $n \in \Z_{\geq 0}$ and 
let $k'$ be a purely transcendental field extension with $n+1 \leq {\rm tr.deg}_k k' <\infty$. 
Then there exists a dense subset $\Lambda \subset (\mathbb P^n_{k'})^*(k')$ 
with respect to the Zariski topology such that the following holds. 

Let (P) be a property for noetherian schemes which satisfies the properties (I)--(III) 
listed in Theorem \ref{intro-Bertini1}. 
Let $\varphi : X \to \mathbb P^n_k$ be a $k$-morphism from a $k$-scheme $X$ of finite type.  
If $X$ is (P) and $\lambda \in \Lambda$, then also $X^{\lambda}_{\varphi}$ is (P), where 
$X^{\lambda}_{\varphi}$ denotes the scheme-theoretic inverse image of the hyperplane 
$H^{\lambda} \subset \mathbb P^n_{k'}$ corresponding to $\lambda$, 
i.e. the following diagram is cartesian: 
\[
\begin{CD}
X^{\lambda}_{\varphi} @>>> X \times_k k'\\
@VVV @VV\varphi \times_k k' V\\
H^{\lambda} @>>> \mathbb P^n_{k'}. 
\end{CD}
\]
\end{thm}



\medskip

\textbf{Overview of contents and sketches of proofs:} 
In Section \ref{s-prelim}, we summarise notation and recall some known results. 
Furthermore, we introduce the notion of $\mathbb A^n$-localising morphisms in Subsection \ref{ss-An-localise}, 
which frequently appears in this paper. 

In Section \ref{s-p-tr-ex}, 
we introduce purely transcendental extensions for rings: $R \subset R(t_1, ..., t_n)$ (Definition \ref{d-tr-ext}). 
We then establish some foundational properties. For example,  
\begin{enumerate}
\renewcommand{\labelenumi}{(\roman{enumi})}
\item the pullback map $\Spm\,R(t_1, ..., t_n) \to \Spm\,R$ is bijective (Theorem \ref{t-tr-ext-maximal}), and 
\item $R(t_1, ..., t_n) \simeq R(t_1, ..., t_{n-1})(t_n)$  (Proposition \ref{p-tr-ext-transitive}). 
\end{enumerate}
Theorem \ref{intro-Bertini1} claims that 
results of Bertini type hold after replacing $X$ by its base change $X \times_R R(t_1, ..., t_n)$. 
Then it is natural to compare $X$ and $X \times_R R(t_1, ..., t_n)$. 
In this direction, we shall establish the following properties under the assumption that $R$ is an excellent ring. 
\begin{enumerate}
\renewcommand{\labelenumi}{(\roman{enumi})}
\setcounter{enumi}{2}
\item $R \to R(t_1, ..., t_n)$ is essentially smooth and faithfully flat (Proposition \ref{p-tr-ext-ff}). 
In particular, $X$ is regular (resp. Cohen-Macaulay, resp. reduced) if and only if 
so is $X \times_R R(t_1, ..., t_n)$. 
\item $X$ is integral (resp. connected) if and only if so is $X \times_R R(t_1, ..., t_n)$  (Theorem \ref{t-tr-bc-conne}). 
\item $\dim X = \dim X \times_R R(t_1, ..., t_n)$  (Theorem \ref{t-tr-ext-dim}).  
\end{enumerate}
We here  give a sketch of the proof only for (v). 
By (ii), we may assume that $n=1$. Set $t :=t_1$, 
The key part is to determine the fibres of $X \times_R R(t) \to X$ over the closed points of $X$. 
More specifically, we show that $\kappa \otimes_R R(t) \simeq \kappa(t)$ 
if $\kappa$ is a field which is a finitely generated $R$-algebra (Theorem \ref{t-tr-fibre}). 
We now describe the outline of how to prove this. 
By the formula $(R/I) \otimes_R R(t) \simeq (R/I)(t)$ (Proposition \ref{p-tr-ext-residue}), 
we may assume that $R$ is a subring of $\kappa$. 
Applying Hilbert's Nullstellensatz and the Artin--Tate lemma, 
the problem is reduced to the case when $\kappa$ is the field of fractions of $R$: 
$\kappa =\Frac\,R$. 
Since $\Frac\,R$ is a finitely generated $R$-algebra, 
it follows from Ratliff's existence theorem that 
$R$ is a semi-local ring with $\dim R \leq 1$. 
Taking the normalisation $R^N$ of $R$, the problem is finally reduced to the case when 
$R$ is a principal ideal domain (Lemma \ref{l-tr-ext-normaln}). 
In this case, we can directly check the required equation $\kappa \otimes_R R(t) \simeq \kappa(t)$ (Lemma \ref{l-pid-frac}). 
For more details, see Subsection \ref{ss-tr-ext-bc}.

In Section \ref{s-main}, we prove a main theorem: Theorem \ref{intro-Bertini1}. 
We now explain the main idea of the proof of Theorem \ref{intro-Bertini1}, 
which is similar to the one of \cite[Theorem 1]{CGM86}. 
For simplicity, we assume that $k:=R$ is a field and (P) is regular. 
We have the following diagram in which each square is cartesian: 
\[
\begin{CD}
X^{\gen}_{\varphi} @>>> X_{k(n)} @>>> X\\
@VVV @VVV @VVV\\
(\mathbb P^n_k)^{\gen}_{{\rm id}} @>>> \mathbb P^n_{k(n)} @>>> \mathbb P^n_k.\\
\end{CD}
\]
In particular, it suffices to find a hyperplane 
$(\mathbb P^n_k)^{\gen}_{{\rm id}} (\simeq \mathbb P^{n-1}_{k(n)})$ of $\mathbb P^n_{k(n)}$ such that 
the induced morphism $(\mathbb P^n_k)^{\gen}_{{\rm id}} \to \mathbb P^n_k$ is essentially smooth. 
For the coordinates $\mathbb P^n_k = \Proj\,k[x_0, ..., x_n]$ and 
$(\mathbb P^n_k)^* \simeq \mathbb P^n_k = \Proj\,k[s_0, ..., s_n]$, we set 
\[
(\mathbb P^n_k)^{\univ}_{{\rm id}} := \{ s_0x_0+ \cdots + s_nx_n =0 \} \subset \mathbb P^n_k \times_k (\mathbb P^n_k)^*. 
\]
Identifying $k(n) =k(t_1, ..., t_n)$ with the function field of $(\mathbb P^n_k)^*$, 
we define $(\mathbb P^n_k)^{\gen}_{{\rm id}}$ as the generic fibre, i.e. the following cartesian square 
\[
\begin{CD}
(\mathbb P^n_k)^{\gen}_{{\rm id}} @>>>
(\mathbb P^n_k)^{\univ}_{{\rm id}}\\
@VVV @VVV\\
\Spec\,k(n) @>>> (\mathbb P^n_k)^*. 
\end{CD}
\]
Then it is almost straightforward to prove that $(\mathbb P^n_k)^{\gen}_{{\rm id}} \to \mathbb P^n_k$ 
is essentially smooth. For more details, see Proposition \ref{p-Pn-compute}.

In Section \ref{s-density}, we  prove Theorem \ref{intro-Bertini2}. 
We now overview the proof of Theorem \ref{intro-Bertini2}. 
Assume that $X$ is (P). 
For the universal family $X^{\univ}_{\varphi} \subset X \times_k (\mathbb P^n_k)^*$, 
we consider the generic fibre $X^{\gen}_{\varphi}$ of $X^{\univ}_{\varphi} \to (\mathbb P^n_k)^*$. 
Set  
\[
\Lambda := \{ \lambda \in (\mathbb P^n_{k'})^*(k')\,|\, {\rm (I), (II)}\}.
\]
\begin{enumerate}
\item[(I)] The image of $\lambda$ to $(\mathbb P^n_k)^*$ is equal to the generic point of $(\mathbb P^n_k)^*$.
\item[(II)] The field extension $K((\mathbb P^n_k)^*) \subset \kappa(\lambda)$ induced by (I) is a purely transcendental extension. 
\end{enumerate}
Then it suffices to prove the following. 
\begin{enumerate}
\item $\Lambda$ is a dense subset of $(\mathbb P^n_{k'})^*(k')$. 
\item If $\lambda \in \Lambda$, then $X^{\lambda}_{\varphi}$ is (P). 
\end{enumerate}

(1) 
Fix a complement of a hyperplane: $\mathbb A^n_k = \Spec\,k[z_1, ..., z_n] \subset (\mathbb P^n_k)^*$. 
For a transcendental basis $\{t_i \}_{i \in I}$ of the purely transcendental extension $k \subset k'$, 
it follows from  ${\rm tr.deg}_k k' \geq n+1$ that $|I| \geq n+1$. 
Then there exist distinct elements $i_0, i_1, ..., i_n \in I$. 
We define $\Lambda'$ as the set of the maximal ideals $\lambda$ of $k'[z_1, ..., z_n]$ of the form 
\[
\lambda = (z_1 - t_{i_1} - t_{i_0}^{d_1}, ..., z_n - t_{i_n} - t_{i_0}^{d_n})
\]
for $(d_1, ..., d_n) \in (\Z_{>0})^n$. 
Then we can check that $\Lambda' \subset \Lambda$, which implies that $\Lambda$ is dense in $(\mathbb P^n_{k'})^*(k')$. 

(2) 
Fix $\lambda \in \Lambda$. 
By the definition of $\Lambda$, 
$X^{\lambda}_{\varphi}$ is obtained as the base change of $X^{\gen}_{\varphi}$ by the 
purely transcendental field extension $K(\mathbb P(V)) \subset \kappa(\lambda)$, i.e. $X^{\lambda}_{\varphi} = X^{\gen}_{\varphi} \times_{K(\mathbb P(V))} \kappa(\lambda)$. 
Since (P) is stable under taking base changes by purely transcendental field extensions of finite degree, 
it is enough to show that $X^{\gen}_{\varphi}$ is (P). 
This is guaranteed by  Theorem \ref{intro-Bertini1}. 



\medskip

\textbf{Acknowledgements:} 
The author would like to thank Karl Schwede and Shunsuke Takagi for many constructive suggestions and answering questions.  
{He also thanks the referee for many constructive suggestions and reading the manuscript carefully.} 
The author was funded by JSPS KAKENHI Grant numbers JP18K13386 and JP22H01112.

\section{Preliminaries}\label{s-prelim}

\subsection{Notation}\label{ss-notation}

In this subsection, we summarise notation and terminologies used in this paper. 

\begin{enumerate}
\item Throughout this paper, $k$ denotes a field. 
We will freely use the notation and terminology in \cite{Har77}. 
We say that an $\F_p$-scheme $X$ is {\em $F$-finite} 
if the absolute Frobenius morphism $F: X \to X$ is a finite morphism. 
Given a ring $A$, we denote by $\Spm\,A$ the set of the maximal ideals of $A$. 
\item We say that a scheme $X$ is {\em excellent} if $X$ is a noetherian scheme such that $\MO_{X, x}$ is an excellent ring for any $x \in X$. 
\item Given a scheme $X$, $X_{\red}$ be its reduced structure, i.e. 
$X_{\red}$ is the reduced closed subscueme of $X$ such that the induced closed immersion $X_{\red} \to X$ is surjective. 
For an integral scheme $X$, its normalisation is denoted by $X^N$.  
\item For an integral scheme $X$, 
we define the {\em function field} $K(X)$ 
of $X$ as $\MO_{X, \xi}$ for the generic point $\xi$ of $X$. 
Given an integral domain $A$, 
we set $\Frac\,A := K(\Spec\,A)$, which is the  field of fractions of $A$.  
\item For a ring homomorphism $R \to R'$ and an $R$-scheme $X$, 
we set $X \times_R R' := X \times_{\Spec\,R} \Spec\,R'$. 
\item 
Given a ring $A$, we say that an $A$-algebra $B$ is a {\em localisation} of $A$ if 
$B$ is isomorphic to $S^{-1}A$ as $A$-algebras for some multiplicatively closed subset of $A$. 
In this case, also the induced ring homomorphism $A \to S^{-1}A, a \mapsto a/1$ is called a {\em localisation}. 
\item 
A {\em variety} $X$ (over $k$)  
is an integral scheme which is separated and of finite type over $k$. 
\item\label{sym-prod}  
Given a finite-dimensional $k$-vector space $V$, we set $\mathbb P(V) := {\rm Proj}\,(S(V^*))$, 
where $S(W) := \bigoplus_{d=0}^{\infty} S^d(W)$ denotes the symmetric algebra for $W$. 
Note that the definition of $\mathbb P(V)$ depends on literature. 
\item 
Given an integral scheme $X$, we define $K(X):=\MO_{X, \xi}$ 
where $\xi$ is the generic point of $X$. 
\item 
Let $X$ be a scheme of finite type over $k$ and $L$ an invertible sheaf. 
Take $\varphi_0, ..., \varphi_n \in H^0(X, L)$ and set $\mathbb P^n_k := \Proj\,k[t_0, ..., t_n]$. 
We define $\{ t_0\varphi_0+ \cdots + t_n\varphi_n =0\}$ as the closed subscheme of $X \times_k \mathbb P^n_k$ as follows. 
Take an affine open cover $X = \bigcup_{i \in I} X_i$ which trivialises $L$, 
and hence we have an isomorphism $\theta_i: L|_{X_i} \xrightarrow{\simeq} \MO_{X_i}$ for each $i \in I$. 
Then  $\{ t_0\varphi_0+ \cdots + t_n\varphi_n =0\}$ is defined as the closed subscheme of $X \times_k \mathbb P^n_k$ 
that satisfies the following equation for each $i \in I$: 
\[
\{ t_0\varphi_0+ \cdots + t_n\varphi_n =0\}|_{X_i \times_k \mathbb P^n_k} := 
\Proj\,\left(\frac{\MO_X(X_i)[t_0, ..., t_n]}{(t_0\theta_i(\varphi_0|_{X_i})+ \cdots + t_n \theta_i(\varphi_n|_{X_i})}\right).
\]
We can check that this definition does not depend on the choice of the open cover $X= \bigcup_{i \in I} X_i$. 
\item 
Let $X$ be a scheme of finite type over $k$, 
$L$ an invertible sheaf, and 
$V$ a $k$-vector subspace of $H^0(X, L)$. 
For an element $\varphi$ of $V$, 
we have the closed subscheme of $X$ defined as the zeros of $\varphi$, which is called a {\em member} of $V$. 
A {\em general member} of $V$ is a member of some fixed Zariski open subset of $V$. 
For $\varphi \in H^0(X, L)$, we set $Z(\varphi)$ to be the closed subset of $X$ that is the zero locus of $\varphi$. 
Set 
\[
\Bs(V) := \bigcap_{\varphi \in V} Z(\varphi), 
\]
which is a closed subset of $X$. 
Note that we have $\Bs(\{0\}) = X$ by definition. 
We say that $V$ is {\em base point free} if $\Bs(V)=\emptyset$. 
\item\label{negative-An} 
Let $A$ be a ring. 
For $n \in \Z_{\geq 0}$, $A[x_1, ..., x_n]$ denotes the polynomial ring with $n$ variables $x_1, ..., x_n$. 
For $n \in \Z_{<0}$, we set $A[x_1, ..., x_n] := \{0\}$, i.e. the zero ring. 
Similarly, if $n \in \Z_{<0}$, then we set $\mathbb A^n_k := \emptyset$ and $\mathbb P^n_{k} :=\emptyset$. 
\item 
Given an $\R$-divisor $D$ on a normal variety $X$, we define $\MO_X(D)$ by 
\[
\Gamma(U, \MO_X(D)) = \{ \varphi \in K(X)\,|\, ({\rm div}(\varphi) +D)|_U \geq 0\}
\] 
for any open subset $U$ of $X$. 
In particular, we have $\MO_X(D) = \MO_X( \llcorner D \lrcorner)$. 
\item 
{Unless otherwise specified, $R[t_1, ..., t_n]$ denotes the polynomial ring over a ring $R$ with variables $t_1, ..., t_n$. Similarly, $R(t_1, ..., t_n)$ denotes the purely transcendental extension 
of $R$ (of degree $n$) in the sense of Definition \ref{d-tr-ext}.}  
We frequently use the abbreviated notation $\underline{t} =(t_1, ..., t_n)$. 
For example, $k[\underline{t} ] :=k[t_1, ..., t_n]$ denotes the polynomial ring over $k$ 
with variables $t_1, ..., t_n$. 
\item\label{ss-n-hyperplane} 
Let $R$ be a ring. 
We say that $H$ is a {\em hyperplane} of $\mathbb P^n_R = \Proj\,R[x_0, ..., x_n]$ (over $R$) 
if $H$ is flat over $R$ and there exist $a_0, ..., a_n \in R$ such that 
$(a_0, ..., a_n) \neq (0, ..., 0)$ 
and $H = \Proj\,R[x_0, ..., x_n]/(a_0x_0+ \cdots + a_nx_n)$. 
For example, if $a_0 \in R^{\times}$, then $\Proj\,R[x_0, ..., x_n]/(a_0x_0+ \cdots + a_nx_n)$ 
is a hyperplane of $\mathbb P^n_R$. 
\end{enumerate}



\subsection{Singularities}\label{ss-singularities}

\subsubsection{Singularities for schemes}\label{sss-sing-scheme}

Let $X$ be a noetherian scheme. Fix $n \in \Z_{\geq 0}$. 

\begin{enumerate}
\renewcommand{\labelenumi}{(\roman{enumi})}
\item $X$ is $R_n$ (resp. $G_n$) if $\MO_{X, x}$ is regular (resp. Gorenstein) 
for any point $x \in X$ such that $\dim \MO_{X, x} \leq n$. 
$X$ is $S_n$ if 
 if $\MO_{X, x}$ is Cohen--Macaulay for any point $x \in X$ such that $\dim \MO_{X, x} \leq n$. 
\item 
We say that $X$ is a {\em Mori scheme} if 
$X$ is a reduced noetherian scheme such that 
the normalisation $X^N \to X$ is a finite morphism. 
We say that $X$ is {\em semi-normal} (resp. weakly normal) 
if $X$ is a Mori scheme such that $\MO_{X, x}$ is semi-normal (resp. weakly normal) 
for any $x \in X$. 
For basic properties on semi-normal (resp. weakly normal) schemes, 
we refer to \cite{GT80} (resp. \cite{Man80}). 
\item 
We say that $X$ is {\em resolution-rational} 
if 
\begin{itemize}
\item $X$ is a normal excellent scheme, 
\item $X$ admits a dualising complex, and 
\item 
there exists a proper birational morphism $f:X' \to X$ from a regular scheme $X'$ 
such that $Rf_*\MO_{X'} \simeq \MO_X$ and $Rf_*\omega_{X'} \simeq \omega_X$. 
\end{itemize}
It is known that this condition does not depend on the choice of $f:X' \to X$ 
{when $X$ is of finite type over a perfect field \cite[Theorem 1]{CR11}}.
\item 
We say that $X$ is {\em $F$-rational} 
if $X$ is a noetherian $\F_p$-scheme such that 
$\Gamma(U, \MO_X)$ is $F$-rational in the sense of \cite[Definition 4.1]{HH94} 
for any affine open subset $U$ of $X$. 
If $X$ is an $F$-rational excellent scheme, 
then $\MO_{X, x}$ is $F$-rational \cite[Theorem 3.1]{Vel95}. 
\item We say that $X$ is {\em $F$-injective} if 
$X$ is a noetherian $\F_p$-scheme such that $\MO_{X, x}$ is $F$-injective in the 
sense of \cite[Definition 2.1]{DM} for any $x \in X$ (cf. \cite[Proposition 3.3]{DM}). 
\item\label{Schwede-duBois} 
Assume that $X$ is a reduced scheme of finite type over 
a field $k$ is of characteristic zero. 
As for Du Bois singularities, we only use the following characterisation given by \cite[Theorem 4.6]{Sch07}. 
Assume that $X$ is separated and a closed subscheme of a smooth scheme $Y$ over $k$. 
Let $\pi:\widetilde Y \to Y$ is a log resolution of $X$ which is isomorphic over $Y \setminus X$. 
Set $E$ to be the reduced closed subscheme of $\widetilde Y$ that satisfies 
the set-theoretic equation $E = \pi^{-1}(X)$. 
Then $X$ is Du Bois if and only if $\MO_X \simeq R \pi_*\MO_E$ holds. 
\end{enumerate}

\subsubsection{Singularities of minimal model program}

Let $X$ be  an integral normal excellent scheme and 
let $\Delta$ be an $\R$-divisor on $X$. 
Assume that $K_X+\Delta$ is $\R$-Cartier 
and $X$ admits a dualising complex. 
For the singularities in minimal model program (terminal, canonical, klt, plt, dlt, lc), 
we refer to \cite[Definition 2.23]{BMP}. 
If $X$ is of finite type over a field and $\Delta$ is a $\Q$-divisor 
such that $K_X+\Delta$ is $\Q$-Cartier, 
this definition coincides with the one in \cite[Definition 2.8]{Kol13}.  
We say that $X$ is terminal, klt,...etc if so is $(X, 0)$.

\subsubsection{$F$-Singularities}


\begin{dfn}\label{d-F-sing}
Let $X$ be an integral normal noetherian $F$-finite $\F_p$-scheme and 
let $\Delta$ be an effective $\R$-divisor on $X$.  

\begin{enumerate}
\renewcommand{\labelenumi}{(\roman{enumi})}
\item 
The pair $(X, \Delta)$ is {\em globally $F$-regular} 
if, for every effective $\Z$-divisor $D$ on $X$, 
there exists $e \in \Z_{>0}$ such that 
\[
\MO_{X} \xrightarrow{F^e} F_*^e\MO_X \hookrightarrow F_*^e\MO_X( \ulcorner (p^e-1) \Delta \urcorner +D )
\]
splits as an $\MO_{X}$-module homomorphism. 
The pair $(X, \Delta)$ is {\em globally purely $F$-regular} 
if, for every effective $\Z$-divisor $D$ on $X$ such that $\Supp\,D$ contains no irreducible component of $\llcorner \Delta\lrcorner$, 
there exists $e \in \Z_{>0}$ such that 
\[
\MO_{X}\xrightarrow{F^e} F_*^e\MO_X \hookrightarrow F_*^e\MO_X( \ulcorner (p^e-1) \Delta \urcorner +D )
\]
splits as an $\MO_{X}$-module homomorphism. 
The pair $(X, \Delta)$ is {\em globally sharply $F$-split} if 
there exists $e \in \Z_{>0}$ such that 
\[
\MO_{X} \xrightarrow{F^e} F_*^e\MO_X \hookrightarrow F_*^e\MO_X( \ulcorner (p^e-1) \Delta \urcorner )
\]
splits as an $\MO_{X}$-module homomorphism. 
The pair $(X, \Delta)$ is {\em globally $F$-split} if 
\[
\MO_{X} \xrightarrow{F^e} F_*^e\MO_X \hookrightarrow F_*^e\MO_X( \llcorner (p^e-1) \Delta \lrcorner )
\]
splits as an $\MO_{X}$-module homomorphism for any $e \in \Z_{>0}$. 
\item 
Let $f : X \to Y$ be a morphism to a scheme $Y$. 
The pair $(X, \Delta)$ is {\em globally $F$-regular over $Y$} 
(resp. {\em globally purely $F$-regular over $Y$}, 
{\em globally sharply $F$-split over $Y$}, resp. {\em globally $F$-split over $Y$}) 
if,  for some open cover $Y = \bigcup_{i \in I} Y_i$, 
so is $(f^{-1}(Y_i), \Delta|_{f^{-1}(Y_i)})$ for any $i \in I$. 
\item The pair $(X, \Delta)$ is {\em strongly $F$-regular} 
(resp. {\em purely $F$-regular}, resp. {\em sharply $F$-pure}, resp. {\em $F$-pure}) if 
$(X, \Delta)$ is 
globally $F$-regular over $X$ 
(resp. globally purely $F$-regular over $X$, 
resp. globally sharply $F$-split over $X$, 
resp. globally $F$-split over $X$) 
with respect to the identity morphism ${\rm id}:X \to X$. 
\item We say that $X$ is strongly $F$-regular, globally $F$-regular,...etc if so is $(X, 0)$. 
\end{enumerate}
\end{dfn}

\begin{lem}\label{l-F-sing-QvsR}
Let $X$ be an integral normal noetherian $F$-finite $\F_p$-scheme 
and let $f : X \to Y$ be a morphism to a quasi-compact scheme $Y$. 
Let $E$ be an effective $\Z$-divisor and let $\Delta$ be an effective $\R$-divisor on $X$.  
Then the following hold. 
\begin{enumerate}
\item 
If $X$ is affine, 
then the following are equivalent. 
\begin{enumerate}
\item $\MO_X \to F_*^e\MO_X \hookrightarrow F_*^e\MO_X(E)$ splits. 
\item $\MO_{X, x} \to (F_*^e\MO_X)_x \hookrightarrow (F_*^e\MO_X(E))_x$ splits 
for any $x \in X$. 
\end{enumerate}
\item 
Assume that $(X, \Delta)$ is globally sharply $F$-split over $Y$. 
Then there exists a $\Q$-divisor $\Delta'$ such that 
$\Delta \leq \Delta'$, $\Supp\,\Delta = \Supp\,\Delta'$, 
and $(X, \Delta')$ is globally sharply $F$-split over $Y$. 
\item 
If $(X, \Delta)$ is globally $F$-regular over $Y$, 
then there exists an effective $\Q$-divisor $\Delta'$ such that $\Delta \leq \Delta'$ and $(X, \Delta')$ is globally $F$-regular over $Y$. 

\end{enumerate}
\end{lem}

\begin{proof}
The assertion (1) follows from the fact that the 
splitting of $F^e(E) : \MO_X \xrightarrow{F^e} F_*^e\MO_X \hookrightarrow F_*^e\MO_X(E)$ is equivalent to the surjectivity of 
\[
(F^e(E))^* : 
H^0(X, \mathcal Hom_{\MO_X}( F_*^e\MO_X(E), \MO_X)) 
\to 
H^0(X, \mathcal Hom_{\MO_X}( \MO_X, \MO_X)). 
\]

Let us show (2). 
We first treat the case when $Y = \Spec\,\F_p$, i.e. $X$ is globally sharply $F$-split. 
In this case, there exists $e \in \Z_{>0}$ such that 
$\MO_{X} \to F_*^e\MO_X \hookrightarrow F_*^e\MO_X( \ulcorner (p^e-1) \Delta \urcorner )$ splits.  
By definition, also $(X, \Delta')$ is globally sharply $F$-split for $\Delta' := \frac{1}{p^e-1} \ulcorner (p^e-1) \Delta \urcorner$. 
This completes the proof for the case when $Y = \Spec\,\F_p$. 

Let us go back to the general case. 
By Definition \ref{d-F-sing}, 
there exists an open cover $Y = \bigcup_{i \in I} Y_i$ such that 
$(f^{-1}(Y_i), \Delta_i|_{f^{-1}(Y_i)})$ is globally sharply $F$-split, 
i.e. globally sharply $F$-split over $\Spec\,\F_p$. 
Since $Y$ is quasi-compact, we may assume that $I$ is a finite set. 
As we have already proven the case when $Y = \Spec\,\F_p$, 
if $i \in I$, then there exists an effective $\Q$-divisor $\Gamma_i$ 
such that $\Delta|_{f^{-1}(Y_i)} \leq \Gamma_i$, $\Supp\,\Delta|_{f^{-1}(Y_i)}= \Supp\,\Gamma_i$, 
and $(f^{-1}(Y_i), \Gamma_i)$ is globally sharply $F$-split. 
Let $\Delta = \sum_k \delta_k \Delta_k$ and $\Gamma_i = \sum_k \gamma_{k, i} \Gamma_k|_{f^{-1}(Y_i)}$ be the irreducible composition. 
For any $k$ and $i \in I$, we define $\delta'_{k, i} \in \Q \cup \{\infty\}$ as follows: 
\[
\delta'_{k, i} := 
\begin{cases}
\gamma_{k, i} \qquad \text{if }\Delta_k \cap f^{-1}(Y_i) \neq \emptyset\\
\infty \qquad \text{if }\Delta_k \cap f^{-1}(Y_i) =\emptyset. 
\end{cases}
\]
We set 
\[
\delta'_k := \min_{i \in I} \delta'_{k, i} \in \Q_{>0}, 
\qquad 
\Delta' := \sum_k \delta'_k \Delta_k. 
\]
By $\Delta'|_{f^{-1}(Y_i)} \leq \Gamma_i$, 
$(f^{-1}(Y_i), \Delta'|_{f^{-1}(Y_i)})$ is globally sharply $F$-split for any $i \in I$, 
i.e. $(X, \Delta')$ is globally sharply $F$-split over $Y$. 
Thus, (2) holds. 

Let us show (3). 
By the same argument as in (2), 
we may assume that $Y= \Spec\,\F_p$, 
i.e. $X$ is globally $F$-regular. 
Fix  an effective $\Z$-divisor $C$ such that 
$\Supp\, \Delta \subset \Supp\,C$. 
Then there exists $e$ such that 
\[
\MO_X \to F_*^e \MO_X( \ulcorner (p^e-1)\Delta \urcorner + C)
\]
splits. 
Enlarging the irrational coefficients of $\Delta$, 
we can find an effective $\Q$-divisor $\Delta'$ such that $\Delta \leq \Delta'$, 
$\Supp\,\Delta = \Supp\,\Delta'$, 
and 
$\ulcorner (p^e-1)\Delta \urcorner = \ulcorner (p^e-1)\Delta' \urcorner$. 
Since $(X \setminus C, \Delta'|_{X \setminus C})$ is globally $F$-regular, 
it follows from \cite[Theorem 3.9]{SS10} that $(X, \Delta')$ is globally $F$-regular. 
\end{proof}

\subsection{Localisation of rings}

We collect some basic results on localisation of rings. 
Although the results in this subsection might be well known, 
we include the proofs, as the author 
{could not find} appropriate references. 


\begin{lem}\label{l-S-1A-inje}
Let $A$ be a ring. 
Let $S$ and $T$ be multiplicatively closed subsets of $A$ 
such that an $A$-algebra homomorphism $\varphi : S^{-1}A \to T^{-1}A$ exists: 
\[
\theta' : A \xrightarrow{\theta} S^{-1}A \xrightarrow{\varphi} T^{-1}A,   
\] 
where $\theta$ and $\theta'$ denote the induced ring homomorphisms. 
If $\theta' : A \to T^{-1}A$ is injective, 
then also $\varphi : S^{-1}A \to T^{-1}A$ is injective.  
\end{lem}

\begin{proof}
Take $a/s \in S^{-1}A$ such that $\varphi(a/s)=0$ in $T^{-1}A$, where $a \in A$ and $s \in S$. 
We also have 
\[
0 = \varphi(a/1) = \varphi(\theta(a)) = \theta'(a). 
\]
Since $\theta'$ is injective, we get $a=0$ in $A$. 
Hence, $a/s = 0$ holds in $S^{-1}A$, as required. 
\end{proof}

\begin{lem}\label{l-mcs-transitive}
Let $A$ be a ring. 
Let $S$ be a multiplicatively closed subset of $A$ and 
let $T$ be a multiplicatively closed subset of $S^{-1}A$. 
{Then there exists a multiplicatively closed subset $U$ of $A$ such that there is an isomorphim 
$T^{-1}(S^{-1}A) \simeq U^{-1}A$ of $A$-algebras.} 
\end{lem}

\begin{proof}
Set 
\begin{itemize}
\item 
$T' := \{ t' \in A\,|\, \frac{t'}{s'} \in T\text{ for some }s' \in S\}$ and 
\item 
$U := ST' := \{ st' \in A\,|\, s \in S, t' \in T'\}$. 
\end{itemize}
Then $T'$ is a multiplicatively closed subset of $A$. 
Hence, also $U$ is a multiplicatively closed subset of $A$. 

We now prove that the induced ring homomorphism 
\[
\theta : A \to T^{-1}(S^{-1}A)
\]
factors through the induced one $A \to U^{-1}A$. 
Pick $st' \in U$ with $s \in S, t' \in T'$. 
We have $t'/ s' \in T$ for some $s' \in S$. 
Then we obtain 
\[
\theta(st') = \frac{st'/1}{1/1} = \frac{(ss')/1}{1/1} \cdot \frac{t'/s'}{1/1}, 
\]
which is contained in $(T^{-1}(S^{-1}A))^{\times}$, as required. 
Therefore, we get the induced ring homomorphism: 
\[
\theta' : U^{-1}A \to T^{-1}(S^{-1}A), \qquad \frac{a}{u} \mapsto \theta(u)^{-1}\theta(a).
\]
for $a \in A$ and $u \in U$. 

Let us show that $\theta'$ is injective. 
Fix $a \in A$ and assume that $\theta'(a/1) =0$. 
Then there exists $t \in T$ such that $t \cdot (a/1) =0$ holds in $S^{-1}A$. 
We can write $t = t'/s'$ for some $t' \in T'$ and $s' \in S'$. 
Therefore, we obtain $s''t' a=0$ in $A$ for some $t' \in T'$ and $s'' \in S$. 
By $s''t' \in U$, we get $a/1=0/1$ in $U^{-1}A$. 
Thus $\theta'$ is injective. 

Let us show that $\theta'$ is surjective. 
Pick $\zeta \in T^{-1}(S^{-1}A)$. 
We can write $\zeta = \frac{a/s}{t}$ for some $a \in A, s \in S, t \in T$. 
We further have $t = t'/s'$ for some $t' \in T', s' \in S$. 
It holds that 
\[
\zeta =  \frac{a/s}{t} = \frac{a/s}{t'/s'}. 
\]
We can check that  $\zeta = \theta'( as'/ st')$ for $as' \in A, st' \in U$, as required. 
\end{proof}

\begin{lem}\label{l-S-1A-chard}
Let $A$ be a ring and let $S$ be a multiplicatively closed subset of $A$. 
Set 
\begin{itemize}
\item $M(S) := {\rm Im}(\theta^*: \Spm\,S^{-1}A \to \Spec\,A)$, and 
\item $\overline{S} := \bigcap_{\p \in M(S)} (A \setminus \p)$, 
\end{itemize}
where $\theta^* : \Spm\,S^{-1}A \to \Spec\,A$ is defined by $\theta^*(\n) := \theta^{-1}(\n)$ 
for the natural ring homomorphism $\theta: A \to S^{-1}A, a \mapsto a/1$. 
Then the following hold. 
\begin{enumerate}
\item 
The inclusion  $S \subset \overline{S}$ holds. 
Furthermore, the induce $A$-algebra homomorphism $S^{-1}A \to \overline{S}^{-1}A$ is an isomorphism. 
\item 
Let $S_1$ and $S_2$ be multiplicatively closed subsets of $A$. 
If $M(S_1) =M(S_2)$, then an $A$-algebra isomorphism 
$S_1^{-1}A \simeq S_2^{-1}A$ holds. 
\end{enumerate}
\end{lem}

\begin{proof}
Since (1) immediately implies (2), it suffices to show (1). 
Let us show the inclusion  $S \subset \overline{S}$. 
Fix $s \in S$ and $\p \in M(S)$. 
Note that $\p S^{-1}A$ is a maximal ideal of $S^{-1} A$. 
It follows from $s/1 \in (S^{-1}A)^{\times}$ that $s/1 \not\in \p S^{-1}A$. 
Therefore, we obtain $s \not\in \theta^{-1}( \p S^{-1}A) = \p$. 
Since $\p$ is chosen to be an arbitrary element of $M(S)$, 
we get $s \in \overline{S}$, i.e. $S \subset \overline{S}$.

We then obtain an $A$-algebra homomorphism  
$S^{-1}A \to \overline{S}^{-1}A$. 
Fix $t \in \overline{S}$. It is enough to show that $\theta(t) \in (S^{-1}A)^{\times}$.  
Take $\n \in \Spm\,(S^{-1}A)$. 
It holds that 
\[
t \in \overline{S} = \bigcap_{\p \in M(S)} (A \setminus \p) \subset A \setminus \theta^{-1}(\n).
\]
Hence we get $\theta(t) \not\in \n$. 
Therefore, we obtain $t \in (S^{-1}A)^{\times}$, as required. 
\end{proof}

\subsection{$\mathbb A^n$-localising morphisms}\label{ss-An-localise}

In this subsection, we first introduce notions of $\mathbb A^n$-localising morphisms 
and strongly $\mathbb A^n$-localising morphisms (Definition \ref{d-A-local}). 
We then establish some fundamental properties, e.g. we study the behaviour under base changes 
and compositions (Proposition \ref{p-A-local}). 
In a main theorem of this paper (Theorem \ref{t-Bertini1}), we shall prove that 
$X^{\gen}_{\varphi} \to X$ is an $\mathbb A^n$-localising morphism under suitable assumptions. 
In this sense, what we are mainly interested in is  not strongly $\mathbb A^n$-localising morphisms 
but $\mathbb A^n$-localising morphisms. 

\begin{dfn}\label{d-A-local}
Fix $n \in \Z$. 
\begin{enumerate}
\item 
We say that a ring homomorphism $\varphi:A \to B$ is $\mathbb A^n$-{\em localising} 
if $B$ is isomorphic to $S^{-1}(A[x_1, ..., x_n])$ as an $A$-algebra for some 
multiplicatively closed subset $S$ of $A[x_1, ..., x_n]$, 
where $A[x_1, ..., x_n]$ denotes the polynomial ring over $A$ with variables $x_1, ..., x_n$. 
Note that if $n<0$, then we have $A[x_1, ..., x_n] = \{0\}$ and $B=\{0\}$ under our terminology (Subsection \ref{ss-notation}(\ref{negative-An})). 
\item 
We say that a morphism $f:X \to Y$ of schemes is $\mathbb A^n$-{\em localising} if 
\begin{enumerate}
\item $f$ is affine and 
\item there exists an affine open cover $Y = \bigcup_{i \in I} Y_i$ of $Y$ such that 
the induced ring homomorphism $\MO_Y(Y_i) \to \MO_X(f^{-1}(Y_i))$ is 
$\mathbb A^n$-localising for any $i \in I$. 
\end{enumerate}
\item 
We say that a morphism $f:X \to Y$ of schemes is {\em strongly} $\mathbb A^n$-{\em localising} if 
\begin{enumerate}
\item $f$ is affine and 
\item for any affine open subset $Y'$ of $Y$, 
the induced ring homomorphism $\MO_Y(Y') \to \MO_X(f^{-1}(Y'))$ is 
$\mathbb A^n$-localising. 
\end{enumerate}
\end{enumerate}
\end{dfn}

\begin{prop}\label{p-A-local}
Fix $n, m \in \Z$. 
\begin{enumerate}
\item 
Let 
\[
\begin{CD}
X' @>f'>> Y'\\
@VVV @VVV\\
X @>f>> Y 
\end{CD}
\] 
be a cartesian diagram of schemes. 
\begin{enumerate}
\item If $f$ is an $\mathbb A^n$-localising morphism, then also $f'$ is an $\mathbb A^n$-localising morphism. 
\item If $f$ is a strongly $\mathbb A^n$-localising morphism and $Y$ is an affine scheme, then also $f'$ is a strongly $\mathbb A^n$-localising morphism. 
\end{enumerate} 
\item 
Let $\varphi:A \to B$ be a ring homomorphism. 
Then $\varphi$ is an $\mathbb A^n$-localising ring homomorphism 
if and only if the induced morphism $\Spec\,B \to \Spec\,A$ is strongly $\mathbb A^n$-localising. 
\item 
Let $A \xrightarrow{\varphi} B \xrightarrow{\psi} C$ be ring homomorphisms. 
If $\varphi$ is $\mathbb A^n$-localising and $\psi$ is $\mathbb A^m$-localising, then 
$\psi \circ \varphi$ is $\mathbb A^{n+m}$-localising. 
\item 
Let $X \xrightarrow{f} Y \xrightarrow{g} Z$ be morphisms of schemes. 
If $f$ is strongly $\mathbb A^n$-localising and $g$ is $\mathbb A^m$-localising, 
then $g \circ f$ is $\mathbb A^{n+m}$-localising. 
\item 
Let $f_0:X_0 \to Y_0$ be a morphism of affine schemes such that 
the induced ring homomorphism $\MO_{Y_0}(Y_0) \to \MO_{X_0}(X_0)$ is $\mathbb A^n$-localising. 
Let 
\[
\begin{CD}
X @>f>> Y\\
@VVV @VVV\\
X_0 @>f_0>> Y_0, 
\end{CD}
\]
be a cartesian diagram of schemes. 
If $g:Y \to Z$ is an $\mathbb A^m$-localising morphism of schemes, 
then the composite morphism $g \circ f:X \to Z$ is an $\mathbb A^{n+m}$-localising morphism. 
\end{enumerate}
\end{prop}

\begin{proof}
Let us show (1). 
The assertion (a) is clear. 
We now show (b). 
Replacing $Y'$ by an affine open subset of $Y'$, we may assume that $Y'$ is an affine scheme. 
Then all the schemes $X, Y, X',$ and $Y'$ are affine, and hence the assertion (b) holds, 
which completes the proof of (1). 
The assertion (2) directly follows from (b) of (1).

Let us show (3). 
We have 
\[
B = S^{-1}(A[x_1, ..., x_n]) =S^{-1}(A[\underline{x}]) 
\quad\text{and}\quad C= T^{-1}(B[y_1, ..., y_m]) = T^{-1}(B[\underline{y}]),
\]
where $\underline x$ and $\underline y$ denote the multi-variables. 
We then obtain 
\[
C = T^{-1}(B[\underline{y}]) = T^{-1}( (S^{-1}(A[\underline{x}]))[\underline y]). 
\]
We define a multiplicatively closed subset $\widetilde S$ of $A[\underline{x}, \underline{y}]$ by 
\[
\widetilde S := {\rm Im}(S \hookrightarrow A[\underline{x}] \to A[\underline{x}, \underline{y}]), 
\]
where each of $S \hookrightarrow A[\underline{x}]$ and $A[\underline{x}] \to A[\underline{x}, \underline{y}]$ denotes the induced map. 
Then it is easy to see that 
\[
(S^{-1}(A[\underline{x}]))[\underline y] \simeq \widetilde{S}^{-1}(A[\underline{x}, \underline y]). 
\]
Therefore, we obtain 
\[
C = T^{-1}( (S^{-1}(A[\underline{x}]))[\underline y]) \simeq T^{-1}(\widetilde{S}^{-1}(A[\underline{x}, \underline y])) \simeq U^{-1}(A[\underline{x}, \underline y])
\]
for some multiplicatively closed subset $U$ of $A[\underline{x}, \underline y]$ 
(Lemma \ref{l-mcs-transitive}). 
Thus (3) holds.

Let us show (4). 
Let $Z = \bigcup_{i \in I} Z_i$ be an affine open cover of $Z$ such that 
if we set $Y_i := g^{-1}(Z_i)$, 
then $\MO_Z(Z_i) \to \MO_Y(Y_i)$ is $\mathbb A^n$-localising for any $i \in I$. 
Set $X_i := f^{-1}(Y_i)$ for any $i \in I$. 
As $f$ is strongly $\mathbb A^n$-localising, 
$\MO_Y(Y_i) \to \MO_X(X_i)$ is $\mathbb A^n$-localising. 
Then (3) implies that $\MO_Z(Z_i) \to \MO_Y(Y_i) \to \MO_X(X_i)$ is $\mathbb A^{n+m}$-localising. 
Hence (4) holds.

Let us show (5). 
By (b) of (1), $f$ is strongly $\mathbb A^n$-localising. 
Then (4) implies that $g \circ f$ is $\mathbb A^{n+m}$-localising. 
\end{proof}

\section{Purely transcendental extensions for rings}\label{s-p-tr-ex}

\subsection{Purely transcendental extensions for rings}\label{ss-ts-ext1}

In this subsection, 
we first introduce  purely transcendental extensions for rings (Definition \ref{d-tr-ext}): 
\[
R \subset R(t_1, ..., t_n), 
\]
which generalise purely transcendental extensions for fields. 
We then establish some foundational properties as follows: 
\begin{enumerate}
    \item $R(t_1, ..., t_n)/IR(t_1, ..., t_n) \simeq (R/I)(t_1, ..., t_n)$ (Proposition \ref{p-tr-ext-residue}). 
    \item $R \to R(t_1, ..., t_n)$ is faithfully flat and essentially smooth (Proposition \ref{p-tr-ext-ff}). 
    \item The map $\Spm R(t_1, ..., t_n) \to \Spm R$ defined by the pullback is bijective (Theorem \ref{t-tr-ext-maximal}). 
    \item $R(t_1, ..., t_n) \simeq R(t_1, ..., t_{n-1})(t_n)$ (Proposition \ref{p-tr-ext-transitive}). 
\end{enumerate}

\vspace{4mm}

\begin{dfn}\label{d-tr-ext}
Let $R$ be a ring. 
Set 
\begin{enumerate}
\item $R[\underline{t}] := R[t_1, ..., t_n]$, 
\item 
$U_{R[\underline{t}]/R} := 
\bigcap_{\p \in \Spec\,R} (R [\underline{t}] \setminus \p R[ \underline{t}] ) 
= \bigcap_{\m \in \Spm\,R} (R [\underline{t}] \setminus \m R[ \underline{t}] )  $, and
\item 
$R(\underline{t}) := R(t_1, ..., t_n) := U_{R[\underline{t}]/R}^{-1}(R [\underline{t}])$. 
\end{enumerate}
It is obvious that $U_{R[\underline{t}]/R}$ is a multiplicatively closed subset of 
$R[\underline{t}]$. 
We call $R(t_1, ..., t_n)$ a {\em purely transcendental extension of $R$ (of degree $n$)}. 
\end{dfn}

\begin{rem}\label{r-tr-ext}
Let $k$ be a field. 
Let $k(t_1, ..., t_n)$ be as in Definition \ref{d-tr-ext}. 
Then $k(t_1, ..., t_n)$ is 
nothing but a purely transcendental extension of $k$ of degree $n$ in the usual sense. 
This follows from 
\[
U_{k[\underline{t}]/k} = \bigcap_{\p \in \Spec\,k} k[\underline{t}] \setminus \p k[\underline{t}] = k[t_1, ..., t_n] \setminus \{0\}. 
\]
\end{rem}

\begin{lem}\label{l-tr-ext-inje}
Let $R$ be a ring. Then the induced ring homomorphism 
\[
R[t_1, ..., t_n] \to R(t_1, ..., t_n)
\]
is injective. 
In particular, if $R \neq 0$, then $R(t_1, ..., t_n) \neq 0$. 
\end{lem}

\begin{proof}
In what follows, $\underline{t}$ denotes the multi-variable $(t_1, ..., t_n)$. 
Fix $u \in U_{R[\underline t]/R} \subset R[\underline t]$. 
It suffices to show that $u$ is not a zero-divisor of $R[\underline{t}]$. 
Suppose that $u$ is a zero-divisor of $R[\underline{t}]$. 
Then there exists $r \in R \setminus \{0\}$ such that $r u =0$ 
(cf. \cite[iii) of Exercise 2 in Chapter 1]{AM69}). 

Fix $\p \in \Spec\,R$. 
By $ru=0$ and 
\[
u \in  U_{R[\underline t]/R} \subset R[\underline t] \setminus \p R[\underline t], 
\]
we can find a coefficient $u'$ of $u$ 
such that $u' \in R \setminus \p$ and $ru' =0$, 
which implies that $r_{\p}=0$ in $R_{\p}$, where $r_{\p}$ denotes the image: $R \to R_{\p}, r \mapsto r_{\p}$. 
Then it follows from \cite[Proposition 3.8]{AM69} that $r=0$. This is a contradiction. 
\end{proof}

\begin{lem}\label{l-tr-ext-functor}
Let $\varphi : R \to R'$ be a ring homomorphism 
and let $\widetilde{\varphi} : R[t_1, ..., t_n] \to R'[t_1, ..., t_n]$ 
be the ring homomorphism such that 
$\widetilde{\varphi}(t_i)=t_i$ for any $1 \leq i \leq n$ 
and $\widetilde{\varphi}(r) = \varphi(r)$  for any $r \in R$. 
Then 
{$\widetilde{\varphi}$ induces a ring homomorphism}
\[
\psi : R(t_1, ..., t_n)  \to R'(t_1, ..., t_n), 
\]
{i.e., the follwoing diagram is commutative 
\[
\begin{CD}
R[t_1, ..., t_n] @> \widetilde{\varphi}>> R'[t_1, ..., t_n]\\
@VVV @VVV\\
R(t_1, ..., t_n)  @>\psi >> R'(t_1, ..., t_n), 
\end{CD}
\]
where the vertical arrows are the induced ring homomorphisms. }
\end{lem}

\begin{proof}
In what follows, $\underline{t}$ denotes the multi-variable $(t_1, ..., t_n)$. 
We have 
\[
\widetilde{\varphi} : R[\underline{t}] \to R' [\underline{t}]. 
\]
Fix $u(\underline{t}) \in U_{R[\underline{t}]/R}$. 
By $R(\underline{t}) = U_{R[\underline{t}]/R}^{-1}R[\underline{t}]$ and 
$R'(\underline{t}) = U_{R'[\underline{t}]/R'}^{-1}R'[\underline{t}]$, 
it suffices to show $\widetilde{\varphi}(u(\underline{t})) \in  U_{R'[\underline{t}]/R'}$. 
Take $\mathfrak q \in \Spec\,R'$. 
It is enough to prove $\widetilde{\varphi}(u(\underline{t})) \not\in \mathfrak q R'[\underline{t}]$. 
Set $\p := \varphi^{-1}(\mathfrak q)$. 
By $u(\underline{t}) \in U_{R[\underline{t}]/R} \subset R [\underline{t}] \setminus \p R [\underline{t}]$, 
we get $u(\underline{t}) \not\in \p R[\underline{t}] = \varphi^{-1}(\mfq) R[\underline{t}]$, which implies that $\widetilde{\varphi}(u(\underline{t})) \not\in \mathfrak q R'[\underline{t}]$.
\end{proof}

\begin{prop}\label{p-tr-ext-residue}
Let $R$ be a noetherian ring and let $I$ be an ideal of $R$. 
Then the induced ring homomorphism 
\[
\theta : R(t_1, ..., t_n)/I R(t_1, ..., t_n) \to 
(R/I)(t_1, ..., t_n) 
\]
is an isomorphim. 
\end{prop}

\begin{proof}
Let $\underline{t} := (t_1, ..., t_n)$ be the multi-variable. 
We first recall the construction of $\theta$. 
By the natural surjective ring homomorphisms $\pi : R \to R/I$ and 
$\pi' : R[\underline{t}] \to (R/I)[\underline{t}]$, 
we obtain 
$\pi'': R(\underline{t}) \to (R/I)(\underline{t})$ (Lemma \ref{l-tr-ext-functor}). 
It follows from  $\pi''(I)=0$ that 
\[
\theta : R(\underline{t})/IR(\underline{t}) \to (R/I)(\underline t)
\] 
is induced. 

Recall that 
\[
U_{R[\underline{t}]/R} = \bigcap_{\p \in \Spec\,R} (R[\underline{t}] \setminus 
\p R[\underline{t}])\quad\text{and}\quad 
U_{(R/I)[\underline{t}]/(R/I)} = \bigcap_{\mfq \in \Spec\,(R/I)} ((R/I)[\underline{t}] \setminus 
\mfq (R/I)[\underline{t}]).\] 
Set 
\[
\overline{U}_{R[\underline{t}]/R} := \pi'( U_{R[\underline{t}]/R}). 
\]
Note that we have 
\[
R(\underline{t})/IR(\underline{t}) \simeq \overline{U}_{R[\underline{t}]/R}^{-1}((R/I)[\underline{t}])
\]
and 
\[
(R/I)(\underline t) = U_{(R/I)[\underline{t}]/(R/I)}^{-1}((R/I)[\underline{t}]). 
\]
In order to prove that $\theta$ is an isomorphism, it suffices to show that 
there exists an $(R/I)[\underline{t}]$-algebra homomorphism 
\[
(R/I)(\underline t)
\to 
R(\underline{t})/IR(\underline{t}). 
\]
Hence it suffices to show $\overline{U}_{R[\underline{t}]/R} \supset 
U_{(R/I)[\underline{t}]/(R/I)}$. 

Take $\overline{f}(\underline{t}) \in U_{(R/I)[\underline{t}]/(R/I)} \subset (R/I)[\underline{t}]$. 
Pick an element $f(\underline{t}) \in R[\underline t]$ such that $\pi'( f(\underline{t})) = 
\overline{f}(\underline{t})$. 
Take a {set of} generator of $I$: 
\[
I = (\alpha_1, ..., \alpha_r). 
\]
Set $N := \deg f$. 
We then take another lift $f'(x)$ of $\overline{f}(x)$ as follows: 
\[
f'(\underline{t}) := f(\underline{t}) + \alpha_1 t_1^{N+1} + \alpha_2 t_1^{N+2} + \cdots + \alpha_r t_1^{N+r} \in R[\underline t]. 
\]
It suffices to show $f'(\underline{t}) \in U_{R[\underline{t}]/R}$. 

Fix $\p \in \Spec\,R$. 
It is enough to prove $f'(\underline{t}) \not\in \p R [\underline t]$. 
If $I\not\subset \p$, 
then we have $\alpha_i \not\in \p$ for some $i$, which implies 
$f'(\underline{t}) \not\in \p R [\underline t]$, as required. 
Therefore, we may assume that $I \subset \p$. 
In this case, $\mathfrak q := \pi(\p)$ is a prime ideal of $R/I$. 
By $\overline{f}(\underline{t}) \in U_{(R/I)[\underline{t}]/(R/I)}$, 
we have $\overline{f}(\underline{t}) \not\in \mathfrak q (R/I)[\underline{t}]$. 
We then obtain $f'(\underline{t}) \not\in \p R[\underline{t}]$ by $\p = \pi^{-1}(\mathfrak q)$. 
\end{proof}

\begin{prop}\label{p-tr-ext-ff}
Let $R$ be a noetherian ring and let $n \in \Z_{\geq 0}$. 
Then the induced morphism 
\[
f: \Spec\,R(t_1, ..., t_n) \to \Spec\,R
\]
is a strongly $\mathbb A^n$-localising surjective morphism. 
In particular, $f$ is essentially smooth and faithfully flat. 
\end{prop}

\begin{proof}
It is clear that $f$ is 
a strongly $\mathbb A^n$-localising surjective morphism (cf. 
Definition \ref{d-A-local}, Proposition \ref{p-A-local}(2)). 
In particular, $f$ is flat. 
Therefore, it suffices to show that 
$R(t_1, ..., t_n) \otimes_R (R/\m) \neq 0$ for any maximal ideal $\m$ of $R$. 
By Proposition \ref{p-tr-ext-residue}, we obtain
\[
 R(t_1, ..., t_n) \otimes_R (R/\m) \simeq 
(R/\m)(t_1, ..., t_n), 
\]
which is a field (Remark \ref{r-tr-ext}) and hence not the zero ring.  
\end{proof}

The following theorem describes the maximal ideals of $R(t_1, ..., t_n)$.

\begin{thm}\label{t-tr-ext-maximal}
Let $R$ be a noetherian ring. 
Then the following hold. 
\begin{enumerate}
\item 
Let $\mathfrak q$ be a prime ideal of $R[t_1, ..., t_n]$. 
Then $\mathfrak q \cap U_{R[t_1, ..., t_n]/R} = \emptyset$ if and only if 
$\mathfrak q \subset \m R[t_1, ..., t_n]$ for some maximal ideal $\m$ of $R$. 
\item For any maximal ideal $\m$ of $R$, $\m R(t_1, ..., t_n)$ is a maximal ideal of $R(t_1, ..., t_n)$. 
\item 
For the induced ring homomorphism $\epsilon : R \to R(t_1, ..., t_n)$, 
we have the bijection  
\begin{eqnarray*}
\Spm\,R(t_1, ..., t_n) &\to& \Spm\,R\\
\n &\mapsto& \epsilon^{-1}(\n)
\end{eqnarray*}
whose inverse map is given by $\Spm\,R \to \Spm\,R(t_1, ..., t_n),$ 
$\m \mapsto \m R(t_1, ..., t_n)$. 
\end{enumerate}
\end{thm}


\begin{proof}
In what follows, $\underline{t}$ denotes the multi-variable $(t_1, ..., t_n)$. 

Let us show (1).  
Assume $\mathfrak q \subset \m R[\underline{t}]$ 
for some maximal ideal $\m$ of $R$. 
We obtain $\mathfrak q \cap (R[\underline{t}] \setminus \m  R[\underline{t}]) = \emptyset$. 
Then we get $\mathfrak q \cap U_{R[\underline{t}]/R} = \emptyset$.

Take $\mfq \in \Spec\,R[\underline{t}]$ 
such that  $\mfq \not\subset \m R[t]$ for any $\m \in \Spm\,R$. 
It suffices to show $\mfq \cap  U_{R[\underline{t}]/R} \neq \emptyset$. 
For any $\m \in \Spm\,R$, there exists $f_{\m}(\underline{t}) \in \mfq$ such that $f_{\m}(\underline{t}) \not\in \m R[\underline{t}]$. 
For each $\m \in \Spm\,R$, 
let $I_{\m} \subset R$ be the ideal of $R$ generated by the coefficients of $f_{\m}(\underline t)$. 
By $I_{\m} \not\subset \m$, we have $\sum_{\m \in \Spm\,R} I_{\m} =R$. 
Since $R$ is a noetherian ring, there exist $\m_1, ..., \m_r \in \Spm\,R$ such that 
$I_{\m_1} + \cdots + I_{\m_r} =R$. 
For $0 \ll N_2 \ll N_3 \ll \cdots \ll N_r$, 
$I_{\m_1} + \cdots + I_{\m_r}$ is nothing but the ideal generated by 
the coefficients of 
\[
f(\underline{t}) := f_{\m_1}(\underline{t}) + t_1^{N_2} f_{\m_2}(\underline{t}) + 
\cdots + t_1^{N_r} f_{\m_2}(\underline{t}) \in \mfq.  
\]
Therefore, we obtain $f(\underline{t}) \not\in \m R[\underline t]$ for any $\m \in \Spm\,R$. 
Hence, we get $f(\underline{t}) \in \mfq \cap U_{R[\underline t]/R}$. 
Thus (1) holds.


Let us show (2) and (3). 
By (1), we have 
\[
\Spm\,R(\underline{t}) = \{ \m R(\underline{t})\,|\, \m \in \Spm\,R\}. 
\]
Thus (2) holds. 
We obtain $\epsilon^{-1}(\m R(\underline{t})) = \m$. 
Hence we have the following two maps
\[
\Spm\,R(\underline{t}) \to \Spm\,R, \qquad \n \mapsto \epsilon^{-1}(\n)
\]
and 
\[
\Spm\,R \to \Spm\,R(\underline{t}), \qquad \m \mapsto \m R(\underline{t}), 
\]
 which are inverse each other. 
Thus (3) holds. 
\end{proof}

\begin{prop}\label{p-tr-ext-transitive}
Let $R$ be a noetherian ring. 
Then an isomorphism of $R[t_1, ..., t_n]$-algebras 
\[
R(t_1, ..., t_n)\simeq 
R(t_1, ..., t_{n-1})(t_{n}) 
\]
holds. 
\end{prop}

\begin{proof}
In what follows, $\underline{t}$ denotes the multi-variable $(t_1, ..., t_n)$. 
Both hand sides can be written as $S^{-1}(R[\underline{t}])$ for some multiplicatively closed subset of 
$R[\underline{t}]$. 
By Lemma \ref{l-S-1A-chard}(2), 
it suffices to show that 
\[
{\rm Im}(\Spm\,R(\underline{t}) \to \Spec\,R[\underline{t}]) =
{\rm Im}(\Spm\,R(t_1, ..., t_{n-1})(t_{n})  \to \Spec\,R[\underline{t}]). 
\]
By Theorem \ref{t-tr-ext-maximal}, 
each of both sides is equal to $\{ \m R[\underline{t}] \,|\, \m \in \Spm\,R\}$. 
\end{proof}

\begin{prop}\label{p-tr-ext-localise}
Let $R$ be a noetherian ring and let $\m$ be a maximal ideal of $R$. 
Then an isomorphism
\[
R_{\m}(t_1, ..., t_n) \simeq R(t_1, ..., t_n)_{\m R(t_1, ..., t_n)}. 
\]
of $R(t_1, ..., t_n)$-algebras holds. 
\end{prop}

\begin{proof}
In what follows, $\underline{t}$ denotes the multi-variable $(t_1, ..., t_n)$. 
Note that $R(\underline{t})_{\m R(\underline{t})}$ is the local ring of 
$R(\underline{t})$ at the maximal ideal $\m R(\underline{t})$ (Theorem \ref{t-tr-ext-maximal}). 
We can write $R_{\m}(\underline{t}) \simeq S^{-1}R(\underline{t})$ 
for some multiplicatively closed subset $S$ of $R(\underline{t})$.
Since $R_{\m}(\underline{t})$ is a local ring whose maximal ideal is $\m R_{\m}(\underline{t})$, 
it is enough to prove 
the equation 
$\theta^{-1}(\m R_{\m}(\underline{t})) = \m R(\underline{t})$ for 
the induced ring homomorphism $\theta : R(\underline{t}) \to R_{\m}(\underline{t})$. 
The inclusion $\theta^{-1}(\m R_{\m}(\underline{t})) \supset \m R(\underline{t})$ is clear. 
Since $\m R(\underline{t})$ is a maximal ideal, we obtain the required equation
$\theta^{-1}(\m R_{\m}(\underline{t})) = \m R(\underline{t})$. 
\end{proof}

\begin{rem}\label{r-tr-ext-localise}
Let $R$ be a ring and let $S$ be a multiplicatively closed subset of $R$. 
Then the induced ring homomorphism $R(t_1, ,...., t_n) \to (S^{-1}R)(t_1, ..., t_n)$ 
factors through $R(t_1, ,...., t_n) \to S^{-1}(R(t_1, ..., t_n))$: 
\[
R(t_1, ,...., t_n) \to S^{-1}(R(t_1, ..., t_n)) \xrightarrow{\varphi} (S^{-1}R)(t_1, ..., t_n). 
\] 
Indeed, the image of $S$ by $R \to (S^{-1}R)(t_1, ..., t_n)$ is contained 
in the multiplicative group $((S^{-1}R)(t_1, ..., t_n))^{\times}$. 
Although it is tempting to hope that $\varphi$ is an isomorphism, 
this does not hold in general (Proposition \ref{p-local-tex-no}). 
\end{rem}

\subsection{Base changes}\label{ss-tr-ext-bc}

Given $n \in \Z_{>0}$, a noetherian ring $R$, and a morphism $f: X \to \Spec\,R$ of finite type, 
consider the following cartesian diagram 
\[
\begin{CD}
X \times_R R(t_1, ..., t_n) @>>> X\\
@VVV @VVV\\
\Spec\,R(t_1, ..., t_n) @>>> \Spec\,R.
\end{CD}
\]
Since $\Spec\,R(t_1, ..., t_n) \to \Spec\,R$ is essentially smooth and faithfully flat (Proposition \ref{p-tr-ext-ff}), the following hold. 
\begin{enumerate}
\item[(i)] $X$ is regular if and only if $X \times_R R(t_1, ..., t_n)$ is regular. 
\item[(ii)] $X$ is Cohen--Macaulay if anf only if 
$X \times_R R(t_1, ..., t_n)$ is Cohen--Macaulay. 
\end{enumerate}
The purpose of this subsection is to establish 
the following properties. 
\begin{enumerate}
\item[(iii)] $X$ is integral if and only if 
$X \times_R R(t_1, ..., t_n)$ is integral (Theorem \ref{t-tr-bc-conne}). 
\item[(iv)] $\dim X = \dim (X \times_R R(t_1, ..., t_n))$ if $R$ is excellent (Theorem \ref{t-tr-ext-dim}). 
\end{enumerate}

\subsubsection{Integrality}\label{sss-conne}

\begin{lem}\label{l-tr-fibre}
Let $R$ and $A$ be rings and let $\kappa$ be a field. 
Then the following hold. 
\begin{enumerate}
\item If $R \to A$ is a ring homomorphism, 
then the induced ring homomorphism 
\[
A \otimes_R R(t_1, ..., t_n) \to A(t_1, ..., t_n)
\]
is injective. 
\item 
If $R \to \kappa$ is a ring homomorphism, then an isomorphism 
as $\kappa [t_1, ..., t_n]$-algberas 
\[
\kappa \otimes_R R(t_1, ..., t_n) \simeq V^{-1}(\kappa[t_1, ..., t_n])
\]
holds for some multiplicatively closed subset $V$ of $\kappa[t_1, ..., t_n]$ 
satisfying $0 \not\in V$. 
In particular, $\kappa \otimes_R R(t_1, ..., t_n)$ is an integral domain 
with $\dim (\kappa \otimes_R R(t_1, ..., t_n)) \leq n$. 
\end{enumerate}
\end{lem}

\begin{proof}
{Let us show (1).} 
Note that $A[t_1, ..., t_n] \to A(t_1, ..., t_n)$ is injective (Lemma \ref{l-tr-ext-inje}). 
Since both sides are localisations of $A[t_1, ..., t_n]$, 
also $A \otimes_R R(t_1, ..., t_n) \to A(t_1, ..., t_n)$ is injective (Lemma \ref{l-S-1A-inje}). 
{Thus (1) holds. 
Let us show (2). 
By Definition \ref{d-tr-ext}, we obtain $\kappa \otimes_R R(t_1, ..., t_n) \simeq V^{-1}(\kappa[t_1, ..., t_n])$ 
for some multiplicatively closed subset $V$ of $\kappa[t_1, ..., t_n]$. 
It follows from (1) that $0 \not\in V$. 
Thus (2) holds.} 
\end{proof}

\begin{prop}\label{p-tr-bc-red}
Let $R$ be a noetherian ring, $X$ a scheme, 
and $f : X \to \Spec\,R$ a morphism of finite type. 
Consider the cartesian diagrams: 
\[
\begin{CD}
X_{\red} \times_R R(t_1, ..., t_n) @>\gamma >> X_{\red}\\
@VVj' V @VVj V\\
X \times_R R(t_1, ..., t_n) @>\beta >> X\\
@VVf' V @VVf V\\
\Spec\,R(t_1, ..., t_n) @>\alpha >> \Spec\,R,  
\end{CD}
\]
where $j : X_{\red} \hookrightarrow X$ denotes the surjective closed immersion 
from the reduced scheme $X_{\red}$. 
Then the following induced morphism is an isomorphism: 
\[
X_{\red} \times_R R(t_1, ..., t_n)
\xrightarrow{\simeq}
(X \times_R R(t_1, ..., t_n))_{\red}. 
\]
\end{prop}

\begin{proof}
Since $X_{\red}$ is reduced and $R \to R(t_1, ..., t_n)$ is essentially smooth (Definition \ref{d-tr-ext}), also $X_{\red} \times_R R(t_1, ..., t_n)$ is reduced. 
Then the assertion follows from the fact that $j'$ is a surjective closed immersion 
from a reduced scheme $X_{\red} \times_R R(t_1, ..., t_n)$. 
\end{proof}

\begin{thm}\label{t-tr-bc-conne}
Let $R$ be a noetherian ring, $X$ a scheme, 
and $f: X \to \Spec\,R$ a morphism of finite type. 
Then the following hold. 
\begin{enumerate}
\item $X$ is integral if and only if $X \times_R R(t_1, ..., t_n)$ is integral. 
\item $X$ is irreducible if and only if $X \times_R R(t_1, ..., t_n)$ is irreducible. 
\item $X$ is connected if and only if $X \times_R R(t_1, ..., t_n)$ is connected. 
\end{enumerate}
\end{thm}

\begin{proof}
The \lq\lq if" parts for (1)--(3) follow from the fact that 
$X \times_R R(t_1, ..., t_n) \to X$ is faithfully flat (Proposition \ref{p-tr-ext-ff}). 
In what follows, let us consider the \lq\lq only-if" parts. 
 
Let us show (1). 
If $X$ is affine, then the assertion follows from 
Lemma \ref{l-tr-ext-inje} and Lemma \ref{l-tr-fibre}(1). 
Set $Y := X \times_R R[t_1, ..., t_n]$ and $Z := X \times_R R(t_1, ..., t_n)$. 
We consider $Z$ as a subset of $Y$. 
Fix a non-empty affine open subset $X_0$ of $X$. 
Set $Y_0 \subset Y$ and $Z_0 \subset Z$ to be the inverse images of $X_0$. 
Let $\xi$ be the generic point of $Y$. 
We have $\xi \in Y_0$ and $\xi \in Z_0 \subset Z$. 
Therefore, 
$Z$ coincides with the closure of $\{\xi\}$. 
Thus (1) holds. 

The assertion (2) follows from (1) and Proposition \ref{p-tr-bc-red}. 

Let us show (3). 
Let $X = \bigcup_{i \in I} X_i$ be the irreducible decomposition. 
It follows from (2) that each $Z_i := X_i \times_R R(t_1, ..., t_n)$ is irreducible. 
Hence it is enough to show that $Z_i \cap Z_j \neq \emptyset$ 
if $X_i \cap X_j \neq \emptyset$. 
This follows from the surjecitivity of $\alpha : Z= X \times_R R(t_1, ..., t_n) \to X$ 
(Proposition \ref{p-tr-ext-ff}) and 
\[
Z_i \cap Z_j = \alpha^{-1}(X_i) \cap \alpha^{-1}(X_j) = \alpha^{-1}(X_i \cap X_j). 
\]
Thus (3) holds. 
\end{proof}


\begin{prop}\label{p-tr-bc-normal}
Let $R$ be an excellent ring, $X$ an integral scheme, 
and $f : X \to \Spec\,R$ a morphism of finite type. 
For the normalisation $\nu : X^N \to X$ of $X$, consider the cartesian diagrams: 
\[
\begin{CD}
X^N \times_R R(t_1, ..., t_n) @>\gamma >> X^N\\
@VV\nu' V @VV\nu V\\
X \times_R R(t_1, ..., t_n) @>\beta >> X\\
@VVf' V @VVf V\\
\Spec\,R(t_1, ..., t_n) @>\alpha >> \Spec\,R. 
\end{CD}
\]
Then we have the following isomorphism as $X \times_R R(t_1, ..., t_n)$-schemes: 
\[
X^N \times_R R(t_1, ..., t_n)
\simeq  
(X \times_R R(t_1, ..., t_n))^N. 
\]
\end{prop}

\begin{proof}
{Since $R$ is an excellent ring, 
all the schemes appearing above are excellent.} 
Note that both sides are integral normal schemes such that 
the induced morphisms to $X \times_R R(t_1, ..., t_n)$ are finite morphisms. 
Hence these are isomorphic over $X \times_R R(t_1, ..., t_n)$. 
\end{proof}

\subsubsection{Fibres}

\begin{lem}\label{l-dim-leq-1}
Let $R$ be a noetherian integral domain such that $K := \Frac\,R$ is a finitely generated $R$-algebra. 
Then $R$ is a semi-local ring with $\dim R \leq 1$. 
\end{lem}

\begin{proof}
We have 
\[
K = R\left[ \frac{r'_1}{r_1}, ... \frac{r'_m}{r_m}\right]
= R\left[ \frac{1}{r_1}, ... \frac{1}{r_m}\right] = R\left[ \frac{1}{f} \right]
\]
for some $r_1, ..., r_m, r'_1, ..., r'_m \in R$ and 
the product $f := r_1 \cdots r_m$. 
In other words, the equation $D_f = \{(0)\}$ holds in $\Spec\,R$. 
By taking the complements of both hands, we get 
\[
\Spec R \setminus \{(0)\} = V(f) = V(\p_1) \cup \cdots \cup V(\p_r),
\]
where $\p_1, ..., \p_r$ are the minimal prime ideals of $fR$, and hence of height one. 
In particular, $R$ has 
only finitely many primes ideals of height one.  
We then obtain $\dim R \leq 1$ by Ratliff's existence theorem \cite[Theorem 31.2]{Mat89}. 
\end{proof}

\begin{lem}\label{l-tr-ext-normaln}
Let $R \subset S$ be an integral extension of noetherian 
integral domains. 
Then the following induced ring homomorphism is an isomorphism: 
\[
S \otimes_R R(t_1, ..., t_n) \xrightarrow{\simeq} S(t_1, ..., t_n).
\]
\end{lem}

\begin{proof}
Let $\underline t :=(t_1, ..., t_n)$ be the multi-variable. 
Note that 
both sides are localisations of $S[\underline{t}]$. 
In particular, also $S \otimes_R R(\underline{t}) \to S(\underline{t})$ 
is a localisation. 

We claim that the following hold. 
\begin{enumerate}
\item If $\n'$ is a maximal ideal of $S \otimes_R R(\underline{t})$, then its pullback $\n' \cap S$ is a maximal ideal of $S$, i.e. we have the following map 
\begin{equation}\label{e1-tr-ext-normaln}
\Spm\,(S \otimes_R R(\underline{t})) \to \Spm\,S, \qquad \n' \mapsto \n' \cap S. 
\end{equation}
\item The map (\ref{e1-tr-ext-normaln}) is bijective. 
\end{enumerate}
For now, we complete the proof under assuming (1) and (2). 
Since $\Spm\,S(\underline{t}) \to \Spm\,S, \mathfrak l \mapsto \mathfrak l \cap R$ is bijective (Theorem \ref{t-tr-ext-maximal}), we get the bijection by (1) and (2): 
\[
\Spm\,S(\underline{t}) \to \Spm\,(S \otimes_R R(\underline{t})), \qquad \mathfrak l \mapsto \mathfrak l \cap S \otimes_R R(\underline{t}). 
\]
Since $S \otimes_R R(\underline{t}) \to S(\underline{t})$ 
is a localisation, this is faithfully flat, and hence an isomorphism. 

Therefore, it is enough to show (1) and (2). 
Let us show (1). 
Fix $\n' \in \Spm\,(S \otimes_R R(\underline{t}))$ and take its pullbacks: 
\begin{equation}\label{e2-tr-ext-normaln}
\begin{CD}
S @>>> S \otimes_R R(\underline{t})\\
@AAA @AAA\\
R @>>> R(\underline{t}) 
\end{CD}
\qquad\qquad  
\begin{CD}
\n @<<< \n'\\
@VVV @VVV\\
\m @<<< \m'.
\end{CD}
\end{equation}
All the arrows in the diagram (\ref{e2-tr-ext-normaln}) are injective. 
Indeed, since $R \to R(\underline{t})$ is faithfully flat, also $S \to S \otimes_R R(\underline{t})$ 
is faithfully flat, and hence injective. 
Since $R(\underline{t})$ is reduced and 
$\Spec\,(S \otimes_R R(\underline{t})) \to \Spec\,R(\underline{t})$ is surjective, 
$R(\underline{t}) \to S \otimes_R R(\underline{t})$ is injective. 
In what follows, we consider all the rings in the diagram (\ref{e2-tr-ext-normaln}) 
as subrings of $S \otimes_R R(\underline{t})$. 
Furthermore, these rings are integral domains since so is $S \otimes_R R(\underline{t})$ (Theorem \ref{t-tr-bc-conne}(1)). 
Since $R \subset S$ is an integral extension, 
so is $R(\underline{t}) \subset S \otimes_R R(\underline{t})$, and hence $\m' := \n' \cap R(\underline{t})$ is a maximal ideal. 
Therefore, $\m := \m' \cap R$ is a maximal ideal (Theorem \ref{t-tr-ext-maximal}). 
Since $R \to S$ is an integral extension of integral domains, also 
$\n := \n' \cap S$ is a maximal ideal. 
Thus (1) holds. 

Let us show (2). 
It suffices to give the inverse map of the map (\ref{e1-tr-ext-normaln}). 
Pick $\n \in \Spm\,S$. 
Then $\m := \n \cap R$ and $\m' := \m R(\underline{t})$ are maximal ideals (Theorem \ref{t-tr-ext-maximal}). 
We get a prime ideal $\n' \in \Spec\,(S \otimes_R R(\underline{t}))$ such that $\n' \cap S = \n$ and $\n' \cap R(\underline{t}) = \m'$. 
Since $R(\underline{t}) \to S \otimes_R R(\underline{t})$ is an integral extension of integral domains, $\n'$ is a maximal ideal. 
Then the map $\n \mapsto \n'$ given here is nothing but the inverse map of (\ref{e1-tr-ext-normaln}). 
Thus (2) holds. 
\end{proof}

\begin{lem}\label{l-pid-frac}
Let $R$ be a principal ideal domain.  
Set $K := \Frac\,R$. Then the following induced ring homomorphism 
is an isomorphism: 
\[
K \otimes_R R(t_1, ..., t_n) \xrightarrow{\simeq} K(t_1, ..., t_n)
\]
\end{lem}

\begin{proof}
Let $\underline t :=(t_1, ..., t_n)$ be the multi-variable. 
If $\dim R=0$, i.e. $R=K$, then the assertion is clear. 
Hence we may assume that $\dim R=1$. 
Since the induced ring homomorphism $K \otimes_R R(\underline{t}) \to K(\underline{t})$ 
is injective (Lemma \ref{l-tr-fibre}), 
we may consider $K \otimes_R R(\underline{t})$ is a subring of $K(\underline{t})$. 
Inside the field $K(\underline{t})$, we have that 
\[
R(\underline{t}) \otimes_{R} K =
\left\{ \frac{f(\underline{t})}{a g(\underline{t})} \in K(\underline{t}) \,\middle|\, 
a \in R \setminus \{ 0\}, 
f(\underline{t}) \in R[\underline{t}], 
g(\underline{t}) \in \bigcap_{\m \in \Spm\,R} R[\underline{t}] \setminus \m R[t] 
\right\}.
\]
Pick $h(\underline{t}) \in R[\underline{t}] \setminus \{0\}$. 
For a greatest common divisor $a \in R \setminus \{0\}$ of the coefficients of $h(\underline{t})$, 
we have that $h(\underline{t}) = ag(\underline{t})$ for some $g(\underline{t}) \in \bigcap_{\m \in \Spm\,R} R[\underline{t}] \setminus \m R[\underline{t}]$, which implies $R(\underline{t}) \otimes_R K = K(\underline{t})$. 
\end{proof}

\begin{thm}\label{t-tr-fibre}
Let $R$ be an excellent ring and let $\kappa$ be a field which is a finitely generated $R$-algebra. 
Then the following induced ring homomorphism is an isomorphism: 
\[
\kappa \otimes_R R(t_1, ..., t_n) \xrightarrow{\simeq} \kappa(t_1, ..., t_n). 
\]
\end{thm}

\begin{proof}
Let $\underline t :=(t_1, ..., t_n)$ be the multi-variable. 
By Proposition \ref{p-tr-ext-residue}, 
we have $\kappa \otimes_R R(\underline{t}) \simeq \kappa \otimes_{R'} R'(\underline{t})$ 
for the image $R' := {\rm Im}(R \to \kappa)$. 
Hence, after replacing $R$ by $R'$, 
we may assume that $R$ is a subring of $\kappa$. 
Furthermore, by Lemma \ref{l-tr-ext-normaln} {and the assumption that $R$ is excellent}, 
the problem is reduced to the case when $R$ is integrally closed. 
For $K := \Frac\,R$, we have 
\[
R  \subset K \subset \kappa.  
\]
Since $\kappa$ is a finitely generated $R$-algebra, $\kappa$ is also a finitely generated $K$-algebra. 
Then Hilbert's Nullstellensatz implies that $K \subset \kappa$ is a finite extension. 
By the Artin--Tate lemma \cite[Proposition 7.8]{AM69}, also $K$ is a finitely generated $R$-algebra. 
Then Lemma \ref{l-dim-leq-1} implies that 
$R$ is a semi-local ring with $\dim R \leq 1$. 
Therefore, $R$ is 
either a field or a Dedekind domain with finitely many prime ideals, and hence $R$ is a principal ideal domain. 
We then obtain  
\begin{align*}
\kappa \otimes_R R(\underline{t}) 
\simeq \kappa \otimes_K (K \otimes_R R(\underline{t}))
\simeq \kappa \otimes_K K(\underline{t})
\simeq \kappa(\underline{t}), 
\end{align*}
where the second isomorphism holds by Lemma \ref{l-pid-frac} 
and the last one follows from Lemma \ref{l-tr-ext-normaln}. 
\end{proof}


\begin{cor}\label{c-tr-fibre}
Let $R$ be a noetherian ring, $X$ a scheme, and 
$f: X \to \Spec\,R$ a morphism of finite type.  
Let $x \in X$ be a point. 
For the induced morphism 
\[
\alpha : X \times_R R(t_1, ..., t_n) \to X, 
\]
the following hold. 
\begin{enumerate}
\item The fibre $\alpha^{-1}(x)$ is an affine integral scheme. 
\item If $R$ is excellent and $x$ is a closed point, 
then $\alpha^{-1}(x)$ consists of one point, i.e. 
$\alpha^{-1}(x)$ is an affine spectrum of a field. 
\end{enumerate}
\end{cor}

\begin{proof}
The assertion (1) follows from Lemma \ref{l-tr-fibre}(2). 
The assertion (2) holds by Theorem \ref{t-tr-fibre}. 
\end{proof}

\subsubsection{Dimension}\label{sss-dim}

\begin{lem}\label{l-highest-pullback}
Let $A$ be a noetherian ring with $d := \dim A <\infty$. 
Let $\n$ be a maximal ideal of $A[x]$ with ${\rm ht}(\n)=d+1$. 
Then $\m := \n \cap A$ is a maximal ideal of $A$ with ${\rm ht}(\m)=d$. 
\end{lem}

\begin{proof}
By $d = \dim A$, it suffices to prove that the prime ideal $\m = \n \cap A$ satisfies ${\rm ht}(\m) = d$. 
We have ${\rm ht}(\m) \leq \dim A = d$. 
Suppose that ${\rm ht}(\m) < d$. 
It suffices to derive a contradiction. 
By $\m = \n \cap A$, we have that $\n \cap S_{\m}=\emptyset$ for $S_{\m} := A \setminus \m$. 
Therefore, $\n$ comes from $S^{-1}_{\m}(A[x]) = A_{\m}[x]$, i.e. 
$\n = \n S^{-1}_{\m}(A[x]) \cap A[x]$. 
We have 
\[
{\rm ht}(\n) =  {\rm ht}(\n S^{-1}_{\m}(A[x])) \leq \dim S^{-1}_{\m}(A[x])
= \dim A_{\m}[x] =\dim A_{\m} +1 ={\rm ht}(\m)+1 <d+1, 
\]
which contradicts ${\rm ht}(\n)=d+1$. 
\end{proof}

\begin{lem}\label{l-ht-bound}
Let $R$ be an excellent ring 
and let $A$ be a finitely generated $R$-algebra. 
Then $\dim A \otimes_R R(t) \leq \dim A$. 
\end{lem}

\begin{proof}
If $\dim A = \infty$, then there is nothing to show. 
Hence we may assume that $\dim A < \infty$. 
Fix a maximal ideal $\n$ of $A \otimes_R R(t)$. 
It suffices to show ${\rm ht}(\n) \leq \dim A$. 
Set $\mathfrak l := \n \cap A[t]$ and $\m := \n \cap A$: 
\[
A \to A[t] \to A \otimes_R R(t), \qquad \m \leftarrow \mathfrak l \leftarrow \n.  
\]

We now treat the case when $\dim A_{\m} < \dim A$. 
We have 
\[
{\rm ht}(\n) = {\rm ht}(\mathfrak l) \leq \dim A_{\m}[t] = \dim A_{\m} +1 <\dim A+1, 
\]
where the first equality holds by 
$A[t]_{\mathfrak l} \simeq (A \otimes_R R(t))_{\n}$. 
This completes the proof for the case when $\dim A_{\m} < \dim A$. 

We may assume that $\dim A_{\m} = \dim A$. 
In particular, $\m$ is a maximal ideal of $A$. 
Since $A_{\m} \to (A \otimes_R R(t))_{\n}$ is faithfually flat, 
it holds by \cite[Theorem 15.1]{Mat89} that 
\[
{\rm ht}(\n) = 
\dim (A \otimes_R R(t))_{\n}
= \dim A_{\m} + \dim (A \otimes_R R(t) / \m (A \otimes_R R(t)))_{\n}
\]
\[
\leq \dim A + \dim ( (A/\m) \otimes_R R(t)) = \dim A, 
\]
where 
the last equality 
follows from Theorem \ref{t-tr-fibre}. 
\end{proof}

\begin{thm}\label{t-tr-ext-dim}
Let $R$ be an excellent ring, 
$X$ a scheme, and $f: X \to \Spec\,R$ a morphism of finite type. 
Then it holds that 
\[
\dim X = \dim (X \times_R R(t_1, ..., t_n)).
\] 
\end{thm}

\begin{proof}
By Proposition \ref{p-tr-ext-transitive}, we may assume that $n=1$. 
Set $t:=t_1$. 
Since $R \to R(t)$ is faithfully flat (Proposition \ref{p-tr-ext-ff}), 
we obtain $\dim X \leq \dim (X \times_R R(t))$ \cite[Theorem 15.1(2)]{Mat89}. 
The opposite inequality $\dim X \geq \dim (X \times_R R(t))$ 
follows from Lemma \ref{l-ht-bound}. 
\end{proof}

\subsection{Examples}\label{ss-tr-example}

\begin{prop}\label{p-PID}
Let $R$ be an excellent one-dimensional principal ideal domain. 
Set $K := \Frac\,R$. 
Then the following hold. 
\begin{enumerate}
\item $R(t)$ is an excellent one-dimensional principal ideal domain. 
\item $R(t) \otimes_R (R/\m) \simeq (R/\m)(t)$ for any $\m \in \Spm\,R$. 
\item $R(t) \otimes_{R} K \simeq K(t)$. 
\item The induced morphism $\Spec\,R(t) \to \Spec\,R$ is bijective. 
\end{enumerate}
In particular, $\Z(t)$ is a one-dimensional principal ideal domain, 
$\Spec\,\Z(t) \to \Spec\,\Z$ is bijective, $\Z(t) \otimes_{\Z} \Q \simeq \Q(t)$, 
and $\Z(t) \otimes_{\Z} \F_p \simeq \F_p(t)$ for any prime number $p$. 
\end{prop}

\begin{proof}
We have the induced ring homomorphisms
\[
R \to R[t] \to U_{R[t]/R}^{-1}(R[t]) = R(t). 
\]
Note that $\dim R(t)= \dim\,R = 1$ (Theorem \ref{t-tr-ext-dim}). 
Since $R[t]$ is a unique factorisation domain, 
$R(t)$ is a one-dimensional unique factorisation domain, 
and hence a principal ideal domain. 
Thus (1) holds. 
It follows from Proposition \ref{p-tr-ext-residue} and 
Lemma \ref{l-pid-frac} that  (2) and (3) hold, respectively. 
Then (4) holds by (2) and (3). 
\end{proof}

As the following example shows, 
the localisation functor $A \mapsto S^{-1}A$ does not commute with the purely transcendental extension functor $R \mapsto R(t)$. 

\begin{prop}\label{p-local-tex-no}
Let $k$ be a field. 
Set 
\[
R := k[x, y] \quad  {\rm and} \quad S := R \setminus \{0\}.
\]
Then $S^{-1}(R (t))$ and $(S^{-1}R)(t)$ are not isomorphic as rings. 
\end{prop}

\begin{proof}
Set $K := S^{-1}R = k(x, y)$. 
We have  $(S^{-1}R)(t) = (k(x, y))(t) = k(x, y, t)$. 
In particular, it suffices to show that $S^{-1}(R (t))$  is not a field. 
As a subring of ${\rm Frac}(R[t]) = k(x, y, t)$, 
$S^{-1}(R(t))$ can be written as follows: 
\[
S^{-1}(R(t)) = 
S^{-1}( U^{-1}_{k[x, y][t]/k[x,y]} k[x, y, t]) 
\]
\[
=
\left\{
\frac{f(x, y, t)}{s(x, y)u(x, y, t)} \,\middle|\, 
f(x, y, t) \in k[x, y, t], s(x, y) \in S, u(x, y, t) \in U_{k[x, y][t]/k[x,y]}
\right\}. 
\]
We then have $S^{-1}(R(t)) = V^{-1}k[x, y, t]$ for 
\[
V := \{ s(x, y)u(x, y, t) \,|\, s(x, y) \in S, u(x, y, t) \in U_{k[x, y][t]/k[x,y]} \}. 
\]
Set $\p :=(x+ty)k[x, y, t]$, which is a nonzero prime ideal of $k[x, y, t]$. 
It suffices to show that $\p \cap V = \emptyset$. 

Suppose $\p \cap V \neq \emptyset$. 
Pick $\zeta \in \p \cap V$. 
We have 
\[
\zeta = (x+ty)f(x, y, t) = s(x, y)u(x, y, t) \neq 0
\]
for some $f(x, y, t) \in k[x, y, t], s(x, y) \in S, u(x, y, t) \in U_{k[x, y][t]/k[x,y]}$. 
Since $x+ty$ is a prime element of $k[x, y, t]$ and $s(x, y) \not\in (x+ty) k[x, y, t]$, 
there exists $g(x, y, t) \in k[x, y, t]$ such that 
\[
u(x, y, t) = (x+ty) g(x, y, t) \in \m k[x, y, t], 
\]
where $\m := (x, y)k[x, y]$. 
This contradicts 
\[
u(x, y, t) \in U_{k[x, y][t]/k[x,y]} = \bigcap_{\n \in \Spm\,k[x, y]} k[x, y, t] \setminus \n k[x, y, t] \subset k[x, y, t] \setminus \m k[x, y, t]. 
\]
\end{proof}


The following proposition shows that the equation 
$\dim X = \dim (X \times_R R(t_1, ..., t_n))$ in Theorem \ref{t-tr-ext-dim} 
no longer holds if $X$ is essentially of finite type.

\begin{prop}
Let $k$ be a field. 
Then $\dim (B \otimes_k k(t)) > \dim B$ for $B :=k(x)$. 
\end{prop}

\begin{proof}
We have 
\begin{eqnarray*}
B \otimes_k k(t) &=& k(x) \otimes_k k(t)\\
& \simeq &
\left\{
\frac{f(x, t)}{g(x)h(t)} \,\middle|\,
f(x, t) \in k[x, t], g(x) \in k[x] \setminus \{0\}, h[t] \in k[t] \setminus \{0\} 
\right\}\\
&\simeq& V^{-1}k[x, t]
\end{eqnarray*}
for 
\[
V := \{ g(x)h(t) \in k[x, t]\,|\, 
g(x) \in k[x] \setminus \{0\}, h[t] \in k[t] \setminus \{0\} \}. 
\]
By $(x-t) k[x, t] \cap V = \emptyset$, we have that 
$\dim B \otimes_k k(t)  = \dim V^{-1}k[x, t] \neq 0=\dim B$. 
\end{proof}

\begin{prop}\label{e-not-fin-type}
Let $k$ be a field. Set $R := k[x]$. 
Then 
the following hold. 
\begin{enumerate}
\item $R(t)$ is an essentially finitely generated $k(t)$-algebra. 
\item $R(t)$ is not a finitely generated $k(t)$-algebra.
\end{enumerate}
\end{prop}

\begin{proof}
Since (1) is obvious, let us show (2). 
Suppose that $R(t) =(k[x])(t)$ is a finitely generated $k(t)$-algebra. 
We have the induced injective $k(t)$-algebra homomorphism: 
\[
\varphi : k(t)[x] \hookrightarrow (k[x])(t), 
\]
which induces a $k(t)$-morphism: 
\[
\varphi^*: \Spec\,( (k[x])(t)) \to \Spec\,k(t)[x]. 
\]
By Chevalley's theorem, ${\rm Im}(\varphi^*)$ is a constructible set. 
Since $\varphi$ is injective, $\varphi^*$ is dominant. 
Therefore, ${\rm Im}(\varphi^*)$ is a cofinite set, i.e. its complement is a finite set. 

We have morphisms. 
\[
\psi^* \circ \varphi^* : \Spec\,( (k[x])(t)) \xrightarrow{\varphi^*} \Spec\,k(t)[x] \xrightarrow{\psi^*} \Spec\,k[x]. 
\]
By Proposition \ref{p-PID}, the composition $\psi^* \circ \varphi^*$ is bijective and 
\[
{\rm Im}(\varphi^*) = \{ \p k(t)[x]  \,|\, \p \in \Spec\,k[x]\}. 
\]
In particular, we obtain
\[
(x-t^m)k(t)[x]  \not\in {\rm Im}(\varphi^*)
\]
for any $m \in \Z_{>0}$. 
Hence ${\rm Im}(\varphi^*)$ is not a cofinite set. 
\end{proof}

\subsection{A generalisation}\label{ss-tr-ex-general}

In order to treat a gluing problem of our construction (cf. Remark \ref{r-R(n)}), 
we give a slight generalisation of our purely transcendental extensions (Definition \ref{d-gene-well-def}) and establish a auxiliary result (Proposition \ref{p-gene-well-def}). 

\begin{dfn}\label{d-gene-well-def}
Let $R$ be a ring and let $S$ be an $R$-algebra 
such that $\m S$ is a prime ideal of $S$ for any $\m \in \Spm\,R$. 
We set 
\[
U_{S/R} := \bigcap_{\p \in \Spec\,R} (S \setminus \p S) = 
\bigcap_{\m \in \Spm\,R} (S \setminus \m S) 
\]
and 
\[
R(S) := U_{S/R}^{-1}S. 
\]
In particular, if $S = R[t_1, ..., t_n]$, then $R(S) = R(t_1, ..., t_n)$ (Definition \ref{d-tr-ext}). 
\end{dfn}

\begin{prop}\label{p-gene-well-def}
Let $R$ be a ring and let $S$ be an $R$-algebra 
such that $\m S$ is a prime ideal of $S$ for any $\m \in \Spm\,R$. 
Take $f \in U_{S/R}$. 
Then the following hold. 
\begin{enumerate}
\item 
For any $\m \in \Spm\,R$, $\m S_f$ is a prime ideal of $S_f$. 
\item 
The composition $S \xrightarrow{\theta} S_f \xrightarrow{\beta} R(S_f)$ of the induced ring homomorphisms factors through the induced ring homomorphism $\alpha : S \to R(S)$: 
\[
\begin{CD}
S @>\theta >> S_f\\
@VV\alpha V @VV\beta V\\
R(S) @>\theta' >> R(S_f). 
\end{CD}
\]
Furthermore, $\theta' : R(S) \to R(S_f)$ is an isomorphism. 
\end{enumerate}
\end{prop}

\begin{proof}
The assertion (1) follows from $f \in U_{S/R} \subset S \setminus \m S$. 

Let us show (2). 
We first show that 
there exists a ring homomorphism $\theta' : R(S) \to R(S_f)$ 
such that $\theta' \circ \alpha = \beta \circ \theta$. 
Take $u \in U_{S/R} = \bigcap_{\m \in \Spm\,R} (S \setminus \m S)$. 
It suffices to prove $\theta(u) \in U_{S_f/R} = \bigcap_{\m \in \Spm\,R} (S_f \setminus \m S_f)$. 
Fix $\m \in \Spm\,R$. 
Suppose $\theta(u) \in \m S_f$. 
Let us derive a contradiction. 
We have $f^r u \in \m S$ for some $r \in \Z_{>0}$. 
Since $\m S$ is a prime ideal of $S$, we have $f^r \in \m S$ or $u \in \m S$. 
By $f \not\in \m S$, we obtain $u \in \m S$, which contradictis $u \in U_{S/R} = \bigcap_{\m \in \Spm\,R} (S \setminus \m S)$. 
Therefore, we get a ring homomorphism $\theta' : R(S) \to R(S_f)$ such that $\theta' \circ \alpha = \beta \circ \theta$.

We now prove that $\theta' : R(S) \to R(S_f)$ is injective. 
Take $s/u \in R(S)$ with $\theta'(s/u)=0$, where $s \in S$ and $u \in U_{S/R}$. 
Then we have $\theta'(s/1)=0$. 
Hence we obtain $v \cdot  (s/1) =0$ in $S_f$ for some $v \in U_{S_f/R}$. 
We can write $v = t/f^a$ for some $t \in S$ and $a \in \Z_{>0}$. 
Hence, it holds that $f^b ts = 0$ in $S$ for some $b \in \Z_{>0}$. 
By $v \in U_{S_f/R} = \bigcap_{\m  \in \Spm\,R} (S_f \setminus \m S_f)$, 
we get $t \in \bigcap_{\m  \in \Spm\,R} (S \setminus \m S) =U_{S/R}$. 
By $f \in U_{S/R}$, we have that 
\[
\frac{s}{1} = \frac{f^b ts}{f^bt}=0 \qquad {\rm in} \qquad R(S) =U_{S/R}^{-1}S. 
\]
Thus $\theta'$ is injective.

Let us prove that $\theta'$ is surjective. 
Take $\zeta \in R(S_f)$. 
We can write $\zeta = \frac{s/f^a}{v}$ for some $s \in S, a \in \Z_{\geq 0}, v \in U_{S_f/R}$. 
Furthermore, we have $v = t/f^b$ for some $t \in U_{S/R}$ and $b \in \Z_{>0}$. 
Hence it holds that 
\[
\zeta = \frac{s/f^a}{t/f^b}. 
\]
By $f^bs \in S$ and $f^a t \in U_{S/R}$, 
$\zeta$ is the image of $f^bs/f^at \in U_{S/R}^{-1}R = R(S)$, as required. 
\end{proof}

\section{Bertini theorems for generic members}\label{s-main}


\subsection{Universal members and generic members}

The purpose of this subsection is to introduce the definitions and to establish some fundamental properties of the universal members $X^{\univ}_{\varphi}$ and 
the generic members $X^{\gen}_{\varphi}$ of $\varphi :X \to \mathbb P^n_R$. 
We use the word \lq\lq members", 
since the inverse image of a hyperplane of a projective space $\mathbb P^n_k$ over an algebraically closed field $k$ is usually called a member of the corresponding linear system.  

\begin{dfn}\label{d-R(n)}
Let $R$ be a noetherian ring and take $n \in \Z_{\geq 0}$. 
Fix 
\begin{enumerate}
\item an isomorphism $\delta: \mathbb P^n_R \xrightarrow{\simeq}  (\mathbb P^n_R)^*$,  
\item an open immersion $j : \mathbb A^n_R \hookrightarrow \mathbb P^n_R$ 
which is obtained as the complement of a hyperplane of $\mathbb P^n_R$ over $R$, and 
\item a coordinate: $\iota: \Spec\,R[t_1, ..., t_n] \xrightarrow{\simeq} \mathbb A^n_R $. 
\end{enumerate}
Consider  $R$-morphisms 
\[
\Theta: \Spec\,R(n) := \Spec\,R(t_1, ..., t_n) \xrightarrow{\theta^*} 
 \Spec\,R[t_1, ..., t_n] \xrightarrow{\iota, \simeq}
\mathbb A^n_R \overset{j}{\hookrightarrow} \mathbb P^n_R \xrightarrow{\delta, \simeq}  (\mathbb P^n_R)^*, 
\]
where $\theta^*$ is the morphism induced by 
the natural ring homomorphism $\theta : R[t_1, ..., t_n] \to R(t_1, ..., t_n)$. 
The $R$-algebra $R(n)$ is independent of the choices of (1)--(3) up to isomorphisms (Remark \ref{r-R(n)}). 
\end{dfn}

\begin{rem}\label{r-R(n)}
Let us prove that the $R$-algebra $R(n)$ is independent of (1)--(3) up to $R$-algebra isomorphisms. 
It is clear that the $R$-algbera $R(n)$ is independent of the choice of an isomorphism 
$\delta: \P^n_R \xrightarrow{\simeq} (\P^n_R)^*$ in (1). 
By $R(t_1, ..., t_n) = R(R[t_1, ..., t_n])$ (Definition \ref{d-gene-well-def}), $R(n)$ is independent of the choice of a coordinate in (3). 

Let $R(n)'$ be another choice. 
Let $H$ and $H'$ be the hyperplanes of $\P^n_R$ over $R$ such that their complements induce $R(n)$ and $R(n)'$, respectively. 
For $\P^n_R = \Proj\,R[t_0, ..., t_n]$, we may assume that $H$ is defined as $\{ t_0 =0 \}$. 
We can write 
\[
H' = \{a_0t_0+ a_1t_1 + \cdots +a_nt_n=0\} \subset \P^n_R = \Proj\,R[t_0, t_1, ..., t_n]. 
\]
We then have 
\[
\Gamma(\P^n_R \setminus (H \cup H'), \MO_{\P^n_R }) = 
\Gamma(\mathbb A^n_R \setminus H', \MO_{\mathbb A^n_R }) = R[t_1, ..., t_n]_f
\]
for $f :=a_0 + a_1 t_1+ \cdots + a_nt_n \in R[t_1, ..., t_n]$. 
Given $\m \in \Spm\,R$, 
there exists $0 \leq i \leq n$ such that $a_i \not\in \m$, 
because $\{ a_0t_0+ a_1t_1 + \cdots + a_nt_n =0 \} \subset \P^n_R$ is flat over $R$. 
This implies $f \in U_{R[t_1, ..., t_n]/R}$. 
It holds that 
\[
R(n) = R(R[t_1, ..., t_n]) =R(R[t_1, ..., t_n]_f) 
=R(\Gamma(\P^n_R \setminus (H \cup H'), \MO_{\P^n_R })), 
\]
where the second equality follows from 
Proposition \ref{p-gene-well-def}.  
By symmetry, we obtain $R(n)' =R(\Gamma(\P^n_R \setminus (H \cup H'), \MO_{\P^n_R }))$. 
\end{rem}

\begin{dfn}\label{d-gen-member}
Let $R$ be a noetherian ring and fix $n \in \Z_{\geq 0}$. 
Let $X$ be an $R$-scheme of finite type and 
let $\varphi: X \to \mathbb P^n_R$ be an $R$-morphism. 
\begin{enumerate}
\item 
Let $(\mathbb P^n_R)^{\univ}_{{\rm id}}$ be the universal member of the hyperplanes, which is a closed subscheme 
of $\mathbb P^n_R \times_R (\mathbb P^n_R)^*$. 
We set $(\mathbb P^n_R)^{\gen}_{{\rm id}} := (\mathbb P^n_R)^{\univ}_{{\rm id}} \times_{(\mathbb P^n_R)^*} \Spec\,R(n)$. 
In particular, we have the following commutative diagram in which each square is cartesian. 
\[
\begin{tikzcd}
	(\mathbb P^n_R)^{\gen}_{{\rm id}} & (\mathbb P^n_R)^{\univ}_{{\rm id}} \\
	\mathbb P^n_{R(n)} & \mathbb P^n_R \times_R (\mathbb P^n_R)^* & \mathbb P^n_R \\
	\Spec\,R(n) & (\mathbb P^n_R)^* & \Spec\,R
	\arrow[from=3-2, to=3-3]
	\arrow[from=3-1, to=3-2, "\Theta"]
	\arrow[from=2-3, to=3-3]
	\arrow[from=2-2, to=3-2, "{\rm pr}_2"]
	\arrow[from=2-1, to=3-1]
	\arrow[from=2-1, to=2-2]
	\arrow[hook, from=1-2, to=2-2]
	\arrow[hook, from=1-1, to=2-1]
	\arrow[from=1-1, to=1-2]
	\arrow[from=1-2, to=2-3]
	\arrow[from=2-2, to=2-3, "{\rm pr}_1"]
\end{tikzcd}
\]
\item 
We define 
$X^{\univ}_{\varphi} := (\mathbb P^n_R)^{\univ}_{{\rm id}} \times_{\mathbb P^n_R} X$ 
and 
$X^{\gen}_{\varphi} := (\mathbb P^n_R)^{\gen}_{{\rm id}} \times_{\mathbb P^n_R} X$. 
We call $X^{\univ}_{\varphi}$ and $X^{\gen}_{\varphi}$ {\em the universal member} and {\em the generic member} of $\varphi$, respectively. 
In particular, taking the base changes of the above commutative diagram by $\varphi :X \to \mathbb P^n_R$, 
we obtain the following commutative diagram in which each square is cartesian. 
\[
\begin{tikzcd}
	X^{\gen}_{\varphi}  & X^{\univ}_{\varphi} \\
	X \times_R R(n) & X \times_R (\mathbb P^n_R)^* & X\\
	\Spec\,R(n) & (\mathbb P^n_R)^* & \Spec\,R
	\arrow[from=3-2, to=3-3]
	\arrow[from=3-1, to=3-2, "\Theta"]
	\arrow[from=2-3, to=3-3]
	\arrow[from=2-2, to=3-2, "{\rm pr}_2"]
	\arrow[from=2-1, to=3-1]
	\arrow[from=2-1, to=2-2]
	\arrow[hook, from=1-2, to=2-2]
	\arrow[hook, from=1-1, to=2-1]
	\arrow[from=1-1, to=1-2]
	\arrow[from=1-2, to=2-3]
	\arrow[from=2-2, to=2-3, "{\rm pr}_1"]
\end{tikzcd}
\]
\end{enumerate}
\end{dfn}

\begin{rem}\label{r-gene-hyperplane}
We use the same notation as in Definition \ref{d-gen-member}. 
We have the following cartesian diagram in which each vertical arrow is a closed immersion. 
\[
\begin{tikzcd}
X^{\gen}_{\varphi} \arrow[r]\arrow[hook, d] & (\mathbb P^n_R)^{\gen} \arrow[hook, d]\\
X \times_R R(n) \arrow[r,"\varphi \times {\rm id}"]& \mathbb P^n_R \times_R R(n) =\mathbb P^n_{R(n)}
\end{tikzcd}
\]
Then the following hold. 
\begin{enumerate}
\item $(\mathbb P^n_R)^{\gen}$ is a hyperplane of $\mathbb P^n_{R(n)}$ over $R(n)$ by Definition \ref{d-gen-member}(1). 
\item By (1), the ideal sheaf on $X \times_R R(n)$ corresponding to $X^{\gen}_{\varphi}$ 
is locally principal. 
\end{enumerate}
\end{rem}

\begin{rem}\label{r-functorial}
Let $R$ be a noetherian ring and fix $n \in \Z_{\geq 0}$. 
Let $X$ and $Y$ be $R$-schemes of finite type and 
let 
\[
\varphi : X \xrightarrow{f} Y \xrightarrow{\psi} \mathbb P^n_R
\]
be $R$-morphisms. 
By Definition \ref{d-gen-member}, it holds that 
\begin{equation}\label{e1-functorial}
X^{\univ}_{\varphi} \simeq Y^{\univ}_{\psi} \times_Y X\qquad 
{\rm and}\qquad X^{\gen}_{\varphi} \simeq Y^{\gen}_{\psi} \times_Y X.
\end{equation}
\end{rem}

\subsection{Bertini theorems for schemes}

In this subsection, we prove a main theorem of this paper (Theorem \ref{t-Bertini1}).

\begin{prop}\label{p-Pn-compute}
Let $R$ be a noetherian ring and fix $n \in \Z_{>0}$. 
Set $X := \mathbb P^n_R = {\rm Proj}\,R[x_0, ..., x_n]$ 
and let $X = \bigcup_{i=0}^n X_i$ be the open cover with $X_i := D_+(x_i)$ (for the definition of $D_+(x_i)$, 
see \cite[Ch. II, Proposition 2.5]{Har77}). 
For each $i \in \{0, ..., n\}$, 
let $j_i : X_i \hookrightarrow X$ be the induced open immersion. 
Then the following hold. 
\begin{enumerate}
\item For each $i \in \{0, ..., n\}$, an isomorphism of $X_i$-schemes 
\[
(X_i)^{\univ}_{j_i} \simeq X_i \times_R \mathbb P^{n-1}_R 
\]
holds. 
\item 
The induced morphism $\beta:(X)^{\gen}_{{\rm id}} \to X$ is an $\mathbb A^{n-1}$-localising morphism. 
\end{enumerate}
\end{prop}

\begin{proof}
Let us show (1). 
By symmetry, we may assume that $i=0$. 
For $X = \mathbb P^n_{R, x} = \Proj\,R[x_0, ..., x_n]$ 
and $(\mathbb P^n_R)^* \simeq \mathbb P^n_{R, s} = \Proj\,R[s_0, ..., s_n]$, 
$X^{\univ}_{{\rm id}}$ is the effective Cartier divisor on $\P^n_{R, x} \times_R \P^n_{R, s}$ 
whose defining equation is $\{s_0x_0+ \cdots + s_nx_n = 0\}$. 
In particular, for $X_0 = \mathbb A^n_{R, y} = \Spec\,R[y_1, ..., y_n]$ 
with $y_i = x_i/x_0$, 
the assertion (1) holds by the following computation: 
\[
(X_0)^{\univ}_{j_0} \simeq 
X^{\univ}_{{\rm id}} \times_X X_0 \simeq 
\Proj\,\frac{(R[y_1, ..., y_n])[s_0, ..., s_n] }{(s_0+s_1y_1 \cdots + s_ny_n)} 
\]
\[
\simeq (\Spec\,R[y_1, ..., y_n]) \times_R \mathbb P^{n-1}_R = X_0 \times_R \mathbb P^{n-1}_R. 
\]
where the first isomorphism holds by Remark \ref{r-functorial}. 

Let us show (2). 
Take the base change: 
\[
\beta_0 := \beta \times_X X_0 : (X_0)^{\gen}_{j_0} \to X_0. 
\]
It suffices to show that $\beta_0$ is an $\mathbb A^{n-1}$-localising morphism. 
For $t_i := s_i/s_n$ and $\mathbb A^n_{R, t} = \Spec\,R[t_0, ..., t_{n-1}]$, 
we obtain the following diagram in which all the squares are cartesian: 
\[
\begin{tikzcd}
	(X_0)^{\gen}_{j_0} & W & (X_0)^{\univ}_{j_0}  \\
	\mathbb A^n_{R', y} & \mathbb A^n_{R, y} \times_R  \mathbb A^n_{R, t} & \mathbb A^n_{R, y} \times_R \mathbb P^n_{R, s} & \mathbb A^n_{R, y} =X_0\\
	\Spec\,R' & \mathbb A^n_{R, t} & \mathbb P^n_{R, s} & \Spec\,R
	\arrow["\iota_3"', hook, from=1-1, to=2-1]
	\arrow["{\rm pr}"', from=2-1, to=3-1]
	\arrow[from=1-1, to=1-2]
	\arrow[hook, from=1-2, to=1-3]
	\arrow["\iota_2"', hook, from=1-2, to=2-2]
	\arrow[from=2-1, to=2-2]
	\arrow["{\rm pr}_2"', from=2-2, to=3-2]
	\arrow[from=3-1, to=3-2]
	\arrow[hook, from=2-2, to=2-3]
	\arrow["\iota_1"', hook, from=1-3, to=2-3]
	\arrow[hook, from=3-2, to=3-3]
	\arrow[from=2-3, to=2-4]
	\arrow["{\rm pr}_2"', from=2-3, to=3-3]
	\arrow["{\rm pr}"', from=2-4, to=3-4]
	\arrow[from=3-3, to=3-4]
	\arrow["h", from=1-3, to=2-4]
\end{tikzcd}
\]
where $R' := R(t_0, ..., t_{n-1})$ and $\iota_1, \iota_2$, and $\iota_3$ are closed immersions. 
We have that 
\[
W = \Spec\,\frac{R[y_1, ..., y_n, t_0, ..., t_{n-1}]}{(t_0+t_1y_1 \cdots + t_{n-1}y_{n-1} + y_n)} 
\subset \mathbb A^n_{R, y} \times_R  \mathbb A^n_{R, t}. 
\]
It holds that 
\begin{eqnarray*}
&&\Gamma( (X_0)^{\gen}_{j_0}, \MO_{(X_0)^{\gen}_{j_0}}) \\
&\simeq & \frac{R[y_1, ..., y_n, t_0, ..., t_{n-1}]}{(t_0+t_1y_1 \cdots + t_{n-1}y_{n-1} + y_n)} 
\otimes_{R[t_0, ..., t_{n-1}]} R(t_0, ..., t_{n-1})\\
&\simeq & \frac{R[y_1, ..., y_n, t_0, ..., t_{n-1}]}{(t_0+t_1y_1 \cdots + t_{n-1}y_{n-1} + y_n)} 
\otimes_{R[t_0, ..., t_{n-1}]} S_1^{-1}(R[t_0, ..., t_{n-1}])\\
&\simeq & S_2^{-1}
\left( \frac{R[y_1, ..., y_n, t_0, ..., t_{n-1}]}{(t_0+t_1y_1 \cdots + t_{n-1}y_{n-1} + y_n)} \right)\\
&\simeq & S_3^{-1}
\left( R[y_1, ..., y_n, t_1, ..., t_{n-1}] \right)\\
\end{eqnarray*}
where each $S_i$ is a suitable multiplicatively closed subset. 
Thus (2) holds. 
\end{proof}

\begin{thm}\label{t-Bertini1}
Let $R$ be a noetherian ring and fix $n \in \Z_{\geq 0}$. 
Let $X$ be an $R$-scheme of finite type and 
let $\varphi: X \to \mathbb P^n_R$ be an $R$-morphism. 
Then the following hold. 
\begin{enumerate}
\item 
There exists an open cover $X = \bigcup_{i=0}^n X_i$ such that 
an $X_i$-isomorphism 
\[
X^{\univ}_{\varphi} \times_X X_i \simeq \mathbb P^{n-1}_R \times_R X_i 
\]
holds for each $i \in \{0, ..., n\}$. 
\item 
The induced morphism $\beta:X^{\gen}_{\varphi}  \to X$ is an $\mathbb A^{n-1}$-localising morphism. 
\item 
Let (P) be a property for noetherian schemes which satisfies the following properties (I)--(III). 
\begin{enumerate}
\renewcommand{\labelenumii}{(\Roman{enumii})}
\item 
Let $Y$ be a noetherian scheme and let $Y = \bigcup_{j \in J} Y_j$ be an open cover. 
Then $Y$ is (P) if and only if $Y_j$ is (P) for 
{every} $j \in J$. 
\item 
For a noetherian ring $A$, 
if $\Spec\,A$ is (P), then also $\Spec\,A[t]$ is (P), 
where $A[t]$ denotes the polynomial ring over $A$ with one variable. 
\item 
For a noetherian ring $A$ and a multiplicatively closed subset $S$ of $A$, 
if $\Spec\,A$ is (P), then also $\Spec\,S^{-1}A$ is (P). 
\end{enumerate}
If $X$ is (P), then also $X^{\gen}_{\varphi}$ is (P). 
\end{enumerate}
\end{thm}

\begin{proof}
The assertion (1) follows from Remark \ref{r-functorial} and Proposition \ref{p-Pn-compute}(1). 

Let us show (2). 
By Remark \ref{r-functorial}, we have the following cartesian diagram 
\[
\begin{CD}
X^{\gen}_{\varphi}  @>\beta >> X\\
@VVV @VV \varphi  V\\
(\mathbb P^n_R)^{\gen}_{{\rm id}} @>\beta' >> \mathbb P^n_R. 
\end{CD}
\]
It holds by Proposition \ref{p-Pn-compute}(2) that $\beta'$ is an $\mathbb A^{n-1}$-localising morphism. 
Then it follows from Proposition \ref{p-A-local}(1) that also $\beta$ is an $\mathbb A^{n-1}$-localising morphism. 
Thus (2) holds. 
The assertion (3) directly follows from (2). 
\end{proof}

\begin{rem}\label{r-tangent-bdl}
{We use the same notation as in Theorem \ref{t-Bertini1}. 
As a stronger result of Theorem \ref{t-Bertini1}(1), it holds that the induced morphism 
\[
\alpha : X^{\univ}_{\varphi} \to X
\]
is a $\P^{n-1}$-bundle, i.e., there exists a locally free sheaf $E$ on $X$ of rank $n$ such that 
$X^{\univ}_{\varphi} \simeq \P_X(E)$. 
In order to prove this, it is enough to assume that $X = \P^n_R$ and $\varphi : X \to \P^n_R$ is the identity morphism (Remark \ref{r-functorial}). 
By the proof of Proposition \ref{p-Pn-compute}(1), 
we obtain 
\[
X^{\univ}_{{\rm id}} = \{ s_0x_0+ \cdots + s_nx_n=0\} \subset \P^n_{R, x} \times \P^n_{R, s}
\]
for $X = \P^n_{R, x} = \Proj\,R[x_0, ..., x_n]$ and $(\P^n_R)^* \simeq \P^n_{R, s} = \Proj\,R[s_0, ..., s_n]$. 
By using the Euler sequence 
\[
0 \to \MO_{\P^n_R} \to \MO_{\P^n_{R}}(1)^{\oplus n+1} \to T_{\P^n_R/R} \to 0
\]
we can check $X \simeq \P_{\P^n_{R, s}}(T_{\P^n_{R, s}/R})$. 
}
\end{rem}

\begin{thm}\label{t-Bertini2}
Let $R$ be a noetherian ring and fix $n, m \in \Z_{\geq 0}$. 
Let $X$ be an $R$-scheme of finite type and 
let $\varphi: X \to \mathbb P^n_R$ be an $R$-morphism. 
Let (P) be one of the following properties. 
\begin{enumerate}
\item $R_m$. 
\item $S_m$. 
\item $G_m$. 
\item Regular. 
\item Cohen--Macaulay. 
\item Gorenstein. 
\item Reduced. 
\item Normal. 
\item Seminormal. 
\item Weakly normal. 
\item Irreducible or empty. 
\item Integral or empty. 
\item Resolution-rational. 
\item Du Bois (assume $R$ to be a field of characteristic zero).  
\item $F$-rational (assume $R$ to be an excellent $\F_p$-algebra). 
\item $F$-injective (assume $R$ to be an  $\F_p$-algebra). 
\end{enumerate}
If $X$ is (P), then 
$X_{\varphi}^{\gen}$ is (P). 
\end{thm}

\begin{proof}
By Theorem \ref{t-Bertini1}(1), the induced morphism $\beta: X_{\varphi}^{\gen}\to X$ is an $\mathbb A^{n-1}$-localising morphism. 
In particular, $\beta: X_{\varphi}^{\gen}\to X$ is flat and any fibre of $\beta$ is geometrically regular.

It follows from \cite[Theorem 23.9]{Mat89} that the cases (1) and (2) hold true. 
We will handle the case (3) in the next paragraph. 
Then each of the conditions (4)--(8) can be written as a combination of (1)--(3).  
The cases (9) and (10) follow from \cite[Proposition 5.1(a)]{GT80} and \cite[ii) of Proposition (III.3)]{Man80}, respectively.

Let us show (3). 
It suffices to prove that if $\Spec\,A$ is $G_m$, then also $\Spec\,A[x]$ is $G_m$ 
(although this is probably a well-known fact). 
Fix $\mathfrak q \in \Spec\,A[x]$ with $\dim A[x]_{\mathfrak q} \leq m$ and set $\mathfrak p$ to be the pullback of $\mathfrak q$ on $A$. 
Since $A_{\mathfrak p} \to A[x]_{\mathfrak q}$ is a faithfully flat ring homomorphism, we have that 
$\dim A_{\mathfrak p} \leq \dim A[x]_{\mathfrak q} \leq m$. 
As $\Spec\,A$ is $G_m$, $A_{\mathfrak p}$ is Gorenstein. 
We have a ring isomorphism $A[x]_{\mathfrak q} \simeq (A_{\mathfrak p}[x])_{\mathfrak q'}$ 
for some $\mathfrak q' \in \Spec\,A_{\mathfrak p}[x]$. 
Therefore, also $(A_{\mathfrak p}[x])_{\mathfrak q'}$ is Gorenstein, hence so is $A[x]_{\mathfrak q}$. 
Thus (3) holds true. 

Let us treat the case (12). 
Assume that $X$ is an integral scheme. 
We have the following diagram whose square is cartesian: 
\[
\begin{CD}
X @<\alpha << X^{\univ}_{\varphi}  
@<\gamma << 
X^{\gen}_{\varphi} \\
@. @VVV @VVV\\
@. (\mathbb P^n_R)^* @<<< \Spec\,R(n). 
\end{CD}
\]
We have a factorisation: 
\[
(\mathbb P^n_R)^* \overset{\iota}{\hookleftarrow} \mathbb A^n_R \xleftarrow{\theta^*} \Spec\,R(n)
\]
where $\iota$ is an open immersion and $\theta^*$ is 
the morphism induced by $\theta : R[t_1, ..., t_n] \to R(t_1, ..., t_n)$ up to isomorphisms (Definition \ref{d-R(n)}). 
Hence it suffices to show that $X^{\univ}_{\varphi}$ is either integral or empty.
Recall that $\alpha : X^{\univ}_{\varphi} \to X$ is a flat morphism of finite type whose fibres are $\mathbb P^{n-1}$ (Theorem \ref{t-Bertini1}(1)). 
If $n=0$, then $X^{\univ}_{\varphi}$ is empty. 
We may assume that $n>0$. 
Then $\alpha$ is an open map whose fibres are connected. 
Since $X$ is connected, $X^{\univ}_{\varphi}$ is connected. 
Thus (12) holds true. 
Since the induced closed immersion $X_{\red} \hookrightarrow X$ is a universal homeomorphism, 
it follows from Proposition \ref{r-functorial}(2) that the case (12) implies the case (11).

The remaining cases are (13)--(16). 
Let us treat (13). 
Let $f:W \to X$ be a resolution. 
We have a cartesian diagram 
\[
\begin{CD}
W^{\gen}_{f \circ \varphi} \simeq X^{\gen}_{\varphi} \times_X W @>\beta_W>> W\\
@VVf^{\gen} V @VVf V\\
X^{\gen}_{\varphi} @>\beta >> X. 
\end{CD}
\]
Then (13) holds true by the flat base change theorem. 
Concerning (14), 
we may apply a similar argument to (13) by using \cite[Theorem 4.6]{Sch07}.

As for (15), 
the assertion follows from \cite[Theorem 3.1]{Vel95} 
(note that a smooth homomorphism in \cite{Vel95} means a flat ring homomorphism with geometrically regular fibres: \cite[the paragraph immdiately after Theorem 3.1]{Vel95}). 
Finally, (16) holds true by \cite[Corollary 4.2]{DM}. 
\end{proof}

We prove the following lemma for a later usage. 

\begin{lem}\label{l-gene-dominant}
Let $R$ be a noetherian ring and fix $n \in \Z_{\geq 0}$. 
Let $X$ be an irreducible scheme of finite type over $R$ and 
let $\varphi: X \to \mathbb P^n_R$ be an $R$-morphism. 
Let $\beta : X^{\gen}_{\varphi} \to X$ be the induced morphism. 
Then either $X^{\gen}_{\varphi} = \emptyset$  
or $\beta(X^{\gen}_{\varphi})$ contains the generic point $\xi_X$ of $X$. 
\end{lem}

\begin{proof}
We may assume that $X^{\gen}_{\varphi} \neq \emptyset$ and $X$ is an integral scheme (Remark \ref{r-functorial}). 
Note that $X^{\univ}_{\varphi}$ is an integral scheme by the proof of Theorem \ref{t-Bertini2}(12). 
We consider $X^{\gen}_{\varphi}$ as a subset of $X^{\univ}_{\varphi}$ via the induced injection $X^{\gen}_{\varphi} \hookrightarrow X^{\univ}_{\varphi}$. 
Then we can find an affine open subset $\Spec\,A \subset X^{\univ}_{\varphi}$ 
whose pullback to $X^{\gen}_{\varphi}$ can be written as $\Spec\,(S^{-1}A)$ 
for some multiplicatively closed subset $S$ of $A$ with $S^{-1}A \neq 0$. 
Since $X^{\univ}_{\varphi}$ is an integral scheme, 
the image $X^{\gen}_{\varphi}$ contains the generic point 
$\xi_{X^{\univ}_{\varphi}}$ of $X^{\univ}_{\varphi}$. 
Since $\xi_{X^{\univ}_{\varphi}}$ is mapped to the generic point $\xi_X$ of $X$, 
$\beta(X^{\gen}_{\varphi})$ contains $\xi_X$, as required.  
\end{proof}

\subsection{Avoidance property for generic members}

\begin{prop}\label{p-gene-avoid}
Let $R$ be a noetherian ring, $X$ a non-empty $R$-scheme of finite type, and $\varphi : X \to \mathbb P^n_R$ an $R$-morphism with $n\in \Z_{\geq 0}$. 
Then $X \times_R R(n) \neq X^{\gen}_{\varphi}$ holds as sets. 
\end{prop}

\begin{proof}
If $X^{\gen}_{\varphi} = \emptyset$, then there is nothing to show. 
In what follows, we assume $X^{\gen}_{\varphi} \neq \emptyset$. 
In particular, it holds that $n \geq 1$. 

\setcounter{step}{0}

\begin{step}\label{s1-gene-Cartier}
Proposition \ref{p-gene-avoid} holds for the case when $R$ is an algebraically closed field. 
\end{step}

\begin{proof}(of Step \ref{s1-gene-Cartier}) 
Set $\kappa := R$. 
Fix a closed point $x$ of $X$. 
We consider $x$ as a reduced scheme. 
Recall that $X^{\univ}_{\varphi}$ is a closed subscheme of $X \times_{\kappa} (\P^n_{\kappa})^*$ and 
$X^{\gen}_{\varphi}$ is its generic fibre over $K((\P^n_{\kappa})^*)$ (Definition \ref{d-gen-member}). 
Take the scheme-theoretic inverse image of $x$ by the first projection 
$b:X \times_{\kappa} (\P^n_{\kappa})^* \to X$: 
\[
\Gamma := b^{-1}(x) = \{x\} \times_{\kappa} (\P^n_{\kappa})^* \subset X \times_{\kappa} (\P^n_{\kappa})^*. 
\]
We have morphisms 
\[
\pi : \Gamma \cap X^{\univ}_{\varphi} \overset{j}{\hookrightarrow}
 X \times_{\kappa} (\P^n_{\kappa})^* \xrightarrow{{\rm pr}_2} (\P^n_{\kappa})^*, 
\]
where $j$ is the induced closed immersion. 
Since $x \not\in \varphi^{-1}(H)$ holds for any general hyperplane $H$ of $\P^n_{\kappa}$, 
it holds that 
\[
\pi^{-1}(y) = \left( \Gamma \cap X^{\univ}_{\varphi}\right) \cap (X \times_{\kappa} \{y\})  =\emptyset
\]
for any general closed point $y \in (\P^n_{\kappa})^*$. 
By the generic flatness, $\Gamma \cap X^{\univ}_{\varphi}$ does not dominate $(\P^n_{\kappa})^*$, 
and 
hence the generic fibre of $\pi$ 
is empty, i.e. $\{x'\} \cap X^{\gen}_{\varphi} = \emptyset$ 
holds for $x' := x \times_{\kappa} K((\P^n_{\kappa})^*)$.  
This completes the proof of Step \ref{s1-gene-Cartier}. 
\end{proof}

\begin{step}\label{s2-gene-Cartier}
Proposition \ref{p-gene-avoid} holds for the case when $X$ is proper over $R$. 
\end{step}

\begin{proof}(of Step \ref{s2-gene-Cartier}) 
Fix a maximal ideal $\m$ of $R$ that is contained in the image of $X \to \Spec\,R$, 
whose existence is guaranteed by the properness of $X \to \Spec\,R$. 
Set $\kappa := R/\m$ and $X_{\kappa} := X \times_R \kappa$. 
For the induced morphism $\psi := \varphi \times_R \kappa : X_{\kappa} \to \P^n_{\kappa}$, 
we have the following cartesian diagram (Proposition \ref{p-tr-ext-residue}, Definition \ref{d-gen-member}):
\[
\begin{tikzcd}
(X_{\kappa})^{\gen}_{\psi} \arrow[r, hook, "\text{closed immersion}"] 
 \arrow[d, hook, "\text{closed immersion}"'] 
& [20mm]X^{\gen}_{\varphi} \arrow[d, hook, "\text{closed immersion}"]\\
X_{\kappa} \times_{\kappa} \kappa(n) \arrow[r, hook, "\text{closed immersion}"]
& X \times_R R(n).
\end{tikzcd}
\]
In order to show $X^{\gen}_{\varphi} \neq X \times_R R(n)$, it suffices to prove 
$(X_{\kappa})^{\gen}_{\psi} \neq X_{\kappa} \times_{\kappa} \kappa(n)$. 
Replacing $(R, X)$ by $(\kappa, X_{\kappa})$, 
we may assume that $\kappa := R$ is a field. 
Furthermore, taking the base change to the algebraic clsoure $\overline{\kappa}$ of $\kappa$, 
we may assume that $R$ is an algebraically closed field (Lemma \ref{l-tr-ext-normaln}, Definition \ref{d-gen-member}). 
By Step \ref{s1-gene-Cartier}, we completes the proof of Step \ref{s2-gene-Cartier}. 
\end{proof}

\begin{step}\label{s3-gene-Cartier}
Proposition \ref{p-gene-avoid} holds for the case when $X$ is an integral scheme. 
\end{step}

\begin{proof}(of Step \ref{s3-gene-Cartier}) 
Fix a non-empty affine open subset $X'$ of $X$. 
For the induced morphism $\varphi' : X' \hookrightarrow X \xrightarrow{\varphi} \P^n_R$, 
we have the following cartesian diagram (Remark \ref{r-functorial}): 
\[
\begin{tikzcd}
X'^{\gen}_{\varphi'} \arrow[r, hook, "\text{open immersion}"] 
 \arrow[d, hook, "\text{closed immersion}"']
& [20mm]X^{\gen}_{\varphi} \arrow[d, hook, "\text{closed immersion}"]\\
X' \times_R R(n) \arrow[r, hook, "\text{open immersion}"] 
& X \times_R R(n). 
\end{tikzcd}
\]
Note that all the schemes in this diagram are integral (Theorem \ref{t-tr-bc-conne}, Theorem \ref{t-Bertini2}). 
Therefore, $X^{\gen}_{\varphi} = X \times_R R(n)$ if and only if 
$X'^{\gen}_{\varphi'} = X' \times_R R(n)$. 
Hence the problem is reduced to the case when $X$ is affine. 
Pick an open immersion $X \hookrightarrow X''$ to an integral projective $R$-scheme $X''$. 
Applying the above argument again, the problem is reduced to the case when  $X$ is projective over $R$, 
which has been settled in Step \ref{s2-gene-Cartier}. 
This completes the proof of Step \ref{s3-gene-Cartier}.
\end{proof}

\begin{step}\label{s4-gene-Cartier}
Proposition \ref{p-gene-avoid} holds without any additional assumptions. 
\end{step}

\begin{proof}(of Step \ref{s4-gene-Cartier}) 
Fix an irreducible component $Y$ of $X$, 
which we equip with the reduced scheme structure. 
We have the following cartesian diagram (Remark \ref{r-functorial}): 
\begin{equation}\label{e-gene-Cartier}
\begin{tikzcd}
Y^{\gen}_{\psi} \arrow[r, hook] \arrow[d, hook, "\iota_Y"] & X^{\gen}_{\varphi} \arrow[d, hook, "\iota_X"]\\
Y \times_R R(n) \arrow[r, hook] & X \times_R R(n), 
\end{tikzcd}
\end{equation}
where $\psi : Y \hookrightarrow X \xrightarrow{\varphi} \P^n_R$ denotes the induced morphism. 
Suppose the set-theoretic equation $X \times_R R(n) = X^{\gen}_{\varphi}$ holds, 
i.e. the closed immersion $\iota_X$ is surjective. 
Then $\iota_Y$ is surjective, and hence we obtain 
the set-theoretic equation $Y \times_R R(n) = Y^{\gen}_{\psi}$, 
which contradicts Step \ref{s3-gene-Cartier}. 
This completes the proof of Step \ref{s4-gene-Cartier}. 
\end{proof}
Step \ref{s4-gene-Cartier} completes the proof of Proposition \ref{p-gene-avoid}. 
\qedhere

\end{proof}

\begin{thm}\label{t-gene-avoid}
Let $R$ be a noetherian ring, $X$ a non-empty $R$-scheme of finite type, and $\varphi : X \to \mathbb P^n_R$ an $R$-morphism with $n\in \Z_{\geq 0}$. 
Then the following hold. 
\begin{enumerate}
\item If $Y$ is a closed subscheme of $X$, then $Y \times_R R(n) \not\subset X^{\gen}_{\varphi}$ holds as sets. 
\item $X^{\gen}_{\varphi}$ is an effective Cartier divisor on $X \times_R R(n)$. 
\item If $R$ is an excellent ring, then $\dim X^{\gen}_{\varphi} \leq \dim X-1$. 
\end{enumerate}
\end{thm}

\begin{proof}
The assertion (1) follows from the cartesian diagram (\ref{e-gene-Cartier}) and Proposition \ref{p-gene-avoid}. 

Let us show (2). 
By (1), $X^{\gen}_{\varphi}$ contains no irreducible components of $X \times_R R(n)$. 
Then $X^{\gen}_{\varphi}$ is an effective Cartier divisor on $X \times_R R(n)$, 
since $X^{\gen}_{\varphi}$ is a closed subscheme on $X \times_R {R(n)}$ 
whose defining ideal sheaf is locally principal (Remark \ref{r-gene-hyperplane}). 
Thus (2) holds. 

The assertion (3) follows from (2) and $\dim X = \dim X \times_R R(n)$ 
(Theorem \ref{t-tr-ext-dim}). 
\end{proof}

\subsection{Restriction to generic members}\label{ss-rest-gene}

Let $R$ be a noetherian ring and fix $n \in \Z_{\geq 0}$. 
Let $X$ be an integral normal scheme of finite type over $R$ and 
let $\varphi: X \to \mathbb P^n_R$ be an $R$-morphism. 
Then $X_{R(n)} := X \times_R R(n)$ is an integral normal scheme of finite type over $R(n)$ (Proposition \ref{p-tr-ext-ff}, Theorem \ref{t-tr-bc-conne}(1)). 
Recall that the generic member $X^{\gen}_{\varphi}$ is an effective Cartier divisor on $X_{R(n)}$ (Theorem \ref{t-gene-avoid}(2)). 
In particular, we have an exact sequence 
\begin{equation}\label{e1-rest-gene}
0 \to \MO_{X_{R(n)}}(-X^{\gen}_{\varphi}) \to  \MO_{X_{R(n)}} \to \MO_{X^{\gen}_{\varphi}} \to 0. 
\end{equation}
Note that $X^{\gen}_{\varphi}$ might be the empty set, although we are mainly interested in the case when $X^{\gen}_{\varphi} \neq \emptyset$. 
Furthermore, we have $\dim X_{R(n)} = \dim X$ when $R$ is excellent 
(Theorem \ref{t-tr-ext-dim}).

\begin{dfn}\label{d-pull-div-to-gen}
Let $R$ be a noetherian ring and fix $n \in \Z_{\geq 0}$. 
Let $X$ be an integral normal scheme of finite type over $R$ and 
let $\varphi: X \to \mathbb P^n_R$ be an $R$-morphism. 
Let $\beta: X^{\gen}_{\varphi} \to  X$ be the induced morphism. 
Recall that 
\begin{enumerate}
\item[(i)] $X^{\gen}_{\varphi}$ is empty or 
an integral normal scheme of finite type over $R(n)$ (Theorem \ref{t-Bertini2}), and  
\item[(ii)] 
if $P$ is a prime divisor on $X$, 
then $\beta^{-1}(P) \simeq P^{\gen}_{j \circ \varphi}$ is either empty or a prime divisor on $X^{\gen}_{\varphi}$ for the induced closed immersion $j : P \hookrightarrow X$  
(Remark \ref{r-functorial} and Theorem \ref{t-Bertini2}). 
\end{enumerate}
\begin{enumerate}
\item 
For a prime divisor $P$ on $X$, we set 
\[
\beta^*P := 
\begin{cases}
\beta^{-1}(P)  \qquad & \text{if } \beta^{-1}(P)\text{ is a prime divisor}\\
0 \qquad & \text{otherwise.}\\
\end{cases}
\]
\item 
For an $\R$-divisor $D$ on $X$ and its irreducible decomposition $D= \sum_{i=1}^r a_iP_i$, we set 
\[
\beta^*D := \sum_{i=1}^r a_i\beta^*P_i, 
\]
which is an $\R$-divisor on $X^{\gen}_{\varphi}$. 
\end{enumerate}
\end{dfn}

\begin{prop}
Let $R$ be an excellent ring and fix $n \in \Z_{\geq 0}$. 
Let $X$ be an integral normal scheme of finite type over $R$ and 
let $\varphi: X \to \mathbb P^n_R$ be an $R$-morphism. 
Let $D$ be an $\R$-divisor on $X$. 
Set $X_{R(n)} := X \times_R R(n)$ and let $D_{R(n)}$ be the pullback of $D$ to $X_{R(n)}$. 
\begin{enumerate}
\item 
If $P$ and $Q$ are distinct prime divisors on $X$ such that $\beta^{-1}(P) \neq \emptyset$, 
then $\beta^*(P) \neq \beta^*(Q)$. 
\item 
The following hold: 
\[
\llcorner \beta^*D \lrcorner = \beta^*(\llcorner D \lrcorner), \quad
\ulcorner \beta^*D \urcorner = \beta^*(\ulcorner D \urcorner), \quad \text{and} \quad
\{ \beta^*D \} = \beta^*( \{ D \}). 
\]
\item 
There exists an exact sequence of $\MO_{X_{R(n)}}$-modules 
\[
0 \to \MO_{X_{R(n)}}(D_{R(n)} -X^{\gen}_{\varphi}) \to 
\MO_{X_{R(n)}}(D_{R(n)}) \to \MO_{X^{\gen}_{\varphi}}(\beta^*D) \to 0
\]
whose restriction to the regular locus of $X_{R(n)}$ coincides with the natural exact sequence obtained from (\ref{e1-rest-gene}). 
\item 
$\beta^*\MO_X(D) \simeq \MO_{X^{\gen}_{\varphi}}(\beta^*D)$. 
\end{enumerate}
\end{prop}

\begin{proof}
If $X^{\gen}_{\varphi} = \emptyset$, then all the assertions are obvious. 
Hence we may assume that $X^{\gen}_{\varphi} \neq \emptyset$.

Let us show (1). 
Assume that $\beta^{-1}(P) \neq \emptyset$ and $\beta^{-1}(P) = \beta^{-1}(Q)$ (automatically, we have $\beta^{-1}(Q) \neq \emptyset$). 
It suffices to show that $P=Q$. 
This follows from 
\[
P = \overline{\beta( P^{\gen}_{\varphi|_P}) } = \overline{\beta( \beta^{-1}(P))} 
=  \overline{\beta( \beta^{-1}(Q))} 
= \overline{\beta( Q^{\gen}_{\varphi|_Q}) }  = Q, 
\]
where the second and fourth equalities hold by the surjectivity of $\beta$ 
and the first and last equalities follow from Lemma \ref{l-gene-dominant}. 
Thus (1) holds. 
The assertion (2) follows from (1) and the fact that $\beta^{-1}(P)$ is either empty or a prime divisor 
(cf. Definition \ref{d-pull-div-to-gen}).

Let us show (3). 
By (2), we may assume that $D$ is a $\Z$-divisor. 
For the regular locus $U$ of $X$ and the induced open immersion $i: U_{R(n)} \hookrightarrow X_{R(n)}$, 
we obtain $i_*i^* \mathcal F = \mathcal F$ for 
\[
\mathcal F \in \{\MO_{X_{R(n)}}(D_{R(n)} -X^{\gen}_{\varphi}), \MO_{X_{R(n)}}(D_{R(n)}), \MO_{X^{\gen}_{\varphi}}(\beta^*D)\}.
\] 
Therefore, we may assume that $X$ is regular, and hence $D$ is Cartier. 
Then the exact sequence (\ref{e1-rest-gene}) induces the required exact sequence. 
Thus (3) holds. 
The assertion (4) follows from (3). 
\end{proof}

\subsection{Bertini theorems for pairs}

\begin{thm}\label{t-Bertini3.1}
Let $R$ be an excellent ring admitting a dualising complex. 
Let $X$ be an integral normal scheme of finite type over $R$, 
$\Delta$ an 
$\R$-divisor such that $K_X+\Delta$ is $\R$-Cartier, 
and 
$\varphi : X \to \mathbb P^n_R$ an $R$-morphism. 
Let $(P)$ be one of the following properties.
\begin{enumerate}
\setcounter{enumi}{16}
\item Terminal. 
\item Canonical. 
\item Klt. 
\item Plt. 
\item Dlt. 
\item Lc. 
\end{enumerate}
If $(X, \Delta)$ is (P), then 
$(X_{\varphi}^{\gen}, \beta^*\Delta)$ is (P). 
\end{thm}

\begin{proof}
We apply the same argument as in Theorem \ref{t-Bertini2}. 
Take an $\mathbb A^n$-localising morphism $f:Y \to X$. 
It is enough to show that if $(X, \Delta)$ satisfies (P), then $(Y, \Delta_Y := f^*\Delta)$ satisfies (P). 
Since taking blowups commutes with flat base changes, 
if  $(X, \Delta)$ satisfies (P), then $(Y, \Delta_Y)$ satisfies (P) 
(cf. \cite[Proposition 2.15 and Warning immidiately after Proposition 2.15]{Kol13}). 
\end{proof}

\begin{thm}\label{t-Bertini-GFR}
Let $R$ be an $F$-finite noetherian $\F_p$-algebra. 
Let $f:X \to Y$ be an $R$-morphism of schemes which are of finite type over $R$. 
Assume that $X$ is integral and normal. 
Let 
$\Delta$ an effective $\R$-divisor on $X$. 
Let $\psi : Y \to \mathbb P^n_R$ be an $R$-morphism and  set $\varphi := \psi \circ f$: 
\[
\varphi : X \xrightarrow{f} Y \xrightarrow{\psi} \mathbb P^n_R. 
\]
Recall that we have the following cartesian diagram (Remark \ref{r-functorial}): 
\[
\begin{CD}
X^{\gen}_{\varphi}  @>\alpha >> X\\
@VV g V @VVf V\\
Y^{\gen}_{\psi} @>\beta >> Y.\\
\end{CD}
\]
Let $(P)$ be one of the following properties.  
\begin{enumerate}
\renewcommand{\labelenumi}{(\roman{enumi})}
\item Globally $F$-regular. 
\item Globally sharply $F$-split. 
\item Globally $F$-split. 
\end{enumerate}
If $(X, \Delta)$ is (P) over $Y$, 
then 
$X^{\gen}_{\varphi}$ is (P) over $Y^{\gen}_{\psi}$. 
\end{thm}

\begin{proof}
We may replace $Y$ by a piece of an open cover of $Y$. 
In particular, we may assume that 
\begin{enumerate}
\item $X$ is (P) (over $\Spec\,\F_p$), and
\item $Y = \Spec\,A$. 
\end{enumerate}
Since $\beta$ is an $\mathbb A^{n-1}$-localising morphism (Theorem \ref{t-Bertini1}(2)), 
the problem is reduced to the case when 
there exists a sequence of morphisms of affine schemes 
\[
Y^{\gen}_{\psi} =: Y_{\ell+1} \xrightarrow{\beta_{\ell+1}} Y_{\ell} \xrightarrow{\beta_{\ell}} \cdots 
\xrightarrow{\beta_2} Y_1 \xrightarrow{\beta_1} Y_0 := Y, \qquad A_{j} := \Gamma(Y_j, \MO_{Y_j})
\]
such that the following hold. 
\begin{enumerate}
\item[(a)] For any $j \in \{1, ..., \ell\}$, 
$A_j = A_{j-1}[t]$ and $\beta_j : Y_j \to Y_{j-1}$ coincides with the  projection $Y_j = Y_{j-1} \times_k \mathbb A^1_k \to Y_{j-1}$. 
\item[(b)] 
$A_{\ell+1} = S^{-1}A_{\ell}$ for some multiplicative subset $S$ of 
$A_{\ell}$ and $\beta_{\ell+1} : Y_{\ell+1} \to Y_{\ell}$ coincides with 
the morphism induced from the natural ring homomorphism 
$A_{\ell} \to S^{-1}A_{\ell}, a \mapsto a/1$. 
\end{enumerate}
Set $X_j := X \times_Y Y_j$ for each $j \in \{1, ..., \ell,  \ell+1\}$ and we inductively define an effective $\R$-divisor $\Delta_j$ on $X_j$ 
as the pullback of $\Delta_{j-1}$ via the induced morphism $\alpha_j : X_j \to X_{j-1}$. 
Fix $j \in \{1, ..., \ell, \ell+1 \}$ and assume that $(X_{j-1}, \Delta_{j-1})$ is (P). 
It suffices to show that $(X_j, \Delta_j)$ is (P). 

Assume $1 \leq j \leq \ell$. 
Take compactifications $Y'_j := Y_{j-1} \times_k \mathbb P^1_k$ and $X'_j :=X_{j-1} \times_k \mathbb P^1_k$ 
of $Y_j$ and $X_j$, respectively. 
Set $\Delta'_j$ to be the pullback of $\Delta_{j-1}$ to $X'_j$. 
If (P) is (i) or (ii), then $(X'_j, \Delta'_j)$ is (P) by the same argument as in \cite[Proposition 3.1]{GT16}, 
which in turn implies that $(X_j, \Delta_j)$ is (P). 
Let us treat the case when (P) is (iii). 
By the case (ii), which  we have already settled, we may assume that $(X_j, (1-\epsilon)\Delta_j)$ is globally $F$-split for any rational number $0 < \epsilon <1$. 
It follows from \cite[Theorem 3.2]{CTW18} that $(X_j, \Delta_j)$ is globally $F$-split. 
This completes the proof for the case when $1 \leq j \leq \ell$.

Assume $j =\ell+1$. 
In this case, for any $e>0$, the following diagram is cartesian: 
\begin{equation}\label{e1-Bertini-GFR}
\begin{CD}
X_{\ell+1} @>\alpha_{\ell+1} >> X_{\ell}\\
@VV F^e V @VV F^e V\\
X_{\ell+1} @>\alpha_{\ell+1} >> X_{\ell}, 
\end{CD}
\end{equation}
where each $F^e$ denotes the $e$-th iterated absolute Frobenius morphism. 
We first handle the case when (P) is (ii). 
We have a splitting 
\begin{equation}\label{e2-Bertini-GFR}
{\rm id} : \MO_{X_{\ell}} \to F^e_*\MO_{X_{\ell}}( \ulcorner (p^e-1)\Delta_{\ell} \urcorner) \xrightarrow{\sigma} \MO_{X_{\ell}}
\end{equation}
for some $e>0$ and $\sigma$. 
It holds that 
\begin{eqnarray*}
\alpha_{\ell+1}^* F^e_*\MO_{X_{\ell}}( \ulcorner (p^e-1)\Delta_{\ell} \urcorner)
&\simeq& 
F^e_* \alpha_{\ell+1}^* \MO_{X_{\ell}}( \ulcorner (p^e-1)\Delta_{\ell} \urcorner)\\
&\simeq &
F^e_* \MO_{X_{\ell}}( \alpha_{\ell+1}^* \ulcorner (p^e-1)\Delta_{\ell} \urcorner)\\
&=&
F^e_* \MO_{X_{\ell}}(\ulcorner (p^e-1)  \Delta_{\ell+1} \urcorner), 
\end{eqnarray*}
where the first isomorphism follows from the cartesian diagram (\ref{e1-Bertini-GFR}). 
Applying $\alpha_{\ell+1}^*(-)$ to (\ref{e2-Bertini-GFR}), we obtain 
\[
{\rm id} : \MO_{X_{\ell+1}} \to F^e_*\MO_{X_{\ell+1}}( \ulcorner (p^e-1)\Delta_{\ell+1} \urcorner) 
\xrightarrow{\alpha_{\ell+1}^*\sigma} \MO_{X_{\ell+1}}. 
\]
Hence, $(X_{\ell+1}, \Delta_{\ell+1})$ is globally sharply $F$-split. 
This completes the case when (P) is (ii). 
If (P) is (iii), then we can apply the same argument as in the case (ii). 
Assume that (P) is (i). 
In this case, by applying \cite[Corollary 3.10]{SS10}, the problem is reduced to the case when (P) is (ii). 
\end{proof}

\begin{thm}\label{t-Bertini3.2}
Let $R$ be an $F$-finite noetherian $\F_p$-algebra. 
Let $X$ be an integral normal scheme of finite type over $R$, 
$\Delta$ an effective $\R$-divisor, 
and 
$\varphi : X \to \mathbb P^n_R$ an $R$-morphism. 
Let $(P)$ be one of the following properties
\begin{enumerate}
\setcounter{enumi}{22}
\item Strongly $F$-regular. 
\item Sharply $F$-pure. 
\item $F$-pure. 
\end{enumerate}
If $(X, \Delta)$ is (P), then 
$(X_{\varphi}^{\gen}, \beta^*\Delta)$ is (P). 
\end{thm}

\begin{proof}
The assertion follows from Theorem \ref{t-Bertini-GFR}. 
\end{proof}

\begin{thm}\label{t-Bertini4.1}
Let $R$ be an excellent ring admitting a dualising complex. 
Let $X$ be an integral normal scheme of finite type over $R$, 
$\Delta$ an $\R$-divisor such that $K_X+\Delta$ is $\R$-Cartier, 
and 
$\varphi : X \to \mathbb P^n_R$ an $R$-morphism. 
Set $X_{R(n)} := X \times_{R} R(n)$. 
Let $\Delta_{R(n)}$ be the pullback of $\Delta$ on $X_{R(n)}$. 
Then the following hold. 
\begin{enumerate}
\item 
If $(X, \Delta)$ is lc, then $(X_{R(n)}, \Delta_{R(n)}+X^{\gen}_{\varphi})$ is lc. 
\item If $(X, \Delta)$ is klt, then 
$(X_{R(n)}, \Delta_{R(n)}+X^{\gen}_{\varphi})$ is plt. 
\end{enumerate}
\end{thm}

\begin{proof}
Recall that we have the following commutative diagram in which each square is cartesian: 
\[
\begin{CD}
X_{R(n)} @>>> X \times_R (\mathbb P^n_R)^* @>>> X\\
@VVV @VVV @VVV\\
\Spec\,R(n)@>>> (\mathbb P^n_R)^* @>>> \Spec\,R.
\end{CD}
\]
Set $\Delta_{(\mathbb P^n_R)^*}$ to be the pullback of $\Delta$ on $X \times_R (\mathbb P^n_R)^*$. 
Note that $X^{\univ}_{\varphi} \to X$ is a smooth morphism (Theorem \ref{t-Bertini1}(1)).

Let us show (1). 
It follows from \cite[(4) in page 47]{Kol13} that $(X \times_R (\mathbb P^n_R)^*, \Delta_{(\mathbb P^n_R)^*} +X^{\univ}_{\varphi})$ is lc. 
Therefore, also its localisation $(X_{R(n)}, \Delta_{R(n)}+X^{\gen}_{\varphi})$ is lc. 
Thus (1) holds. 
The assertion (2) follows from the same argument as in (1). 
\end{proof}

\begin{thm}\label{t-Bertini4.2}
Let $R$ be an $F$-finite noetherian $\F_p$-algebra. 
Let $X$ be an integral normal scheme of finite type over $R$, 
$\Delta$ an effective $\R$-divisor,  
and 
$\varphi : X \to \mathbb P^n_R$ an $R$-morphism. 
Set $X_{R(n)} := X \times_{R} R(n)$. 
Let $\Delta_{R(n)}$ be the pullback of $\Delta$ on $X_{R(n)}$. 
Then the following hold. 
\begin{enumerate}
\item 
If $(X, \Delta)$ is sharply $F$-pure, then $(X_{R(n)}, \Delta_{R(n)}+X^{\gen}_{\varphi})$ is sharply $F$-pure. 
\item If $R$ is a field and $(X, \Delta)$ is strongly $F$-regular, then 
$(X_{R(n)}, \Delta_{R(n)}+X^{\gen}_{\varphi})$ is purely $F$-regular. 
\end{enumerate}
\end{thm}

\begin{proof}
We use the same notation as in the proof of Theorem \ref{t-Bertini4.1}.

Let us show (1). 
Taking an affine open cover of $X$, we may assume that $X$ is affine. 
Furthermore, the problem is reduced to the case when $\Delta$ is a $\Q$-divisor. 
In this case, $(X, \Delta)$ is globally sharply $F$-split. 
By \cite[Theorem 4.3(ii)]{SS10}, we may assume that $(p^e-1)(K_X+\Delta)$ is Cartier for some $e \in \Z_{>0}$. 
Therefore, 
\begin{enumerate}
\item[(a)] 
$(X \times_R (\mathbb P^n_R)^*, \Delta_{(\mathbb P^n_R)^*})$ is sharply $F$-pure and  
\item[(b)] 
$(X^{\univ}_{\varphi}, \Delta_{(\mathbb P^n_R)^*}|_{X^{\univ}_{\varphi}})$ is sharply $F$-pure. 
\end{enumerate}
By (a), $(X \times_R (\mathbb P^n_R)^*, \Delta_{(\mathbb P^n_R)^*} +X^{\univ}_{\varphi})$ is sharply $F$-pure outside $X^{\univ}_{\varphi}$. 
By (b) and inversion of adjunction \cite[Main Theorem(iii)]{Sch09}, 
$(X \times_R (\mathbb P^n_R)^*, \Delta_{(\mathbb P^n_R)^*} +X^{\univ}_{\varphi})$ is sharply $F$-pure 
around $X^{\univ}_{\varphi}$. 
To summarise, $(X \times_R (\mathbb P^n_R)^*, \Delta_{(\mathbb P^n_R)^*} +X^{\univ}_{\varphi})$ is sharply $F$-pure. 
Hence, its localisation $(X_{R(n)}, \Delta_{R(n)} +X^{\gen}_{\varphi})$ is sharply $F$-pure. 
Thus (1) holds. 

Let us show (2). 
By the same argument as in (1), we may assume that $X$ is affine and $(p^e-1)(K_X+\Delta)$ is Cartier for some $e \in \Z_{>0}$ (use \cite[Theorem 4.3(i)]{SS10} instead of \cite[Theorem 4.3(ii)]{SS10}). 
Then the following hold. 
\begin{enumerate}
\item[(c)] 
$(X \times_R (\mathbb P^n_R)^*, \Delta_{(\mathbb P^n_R)^*})$ is strongly $F$-regular. 
\item[(d)] 
$(X^{\univ}_{\varphi}, \Delta_{(\mathbb P^n_R)^*}|_{X^{\univ}_{\varphi}})$ is strongly $F$-regular. 
\end{enumerate}
By (c), $(X \times_R (\mathbb P^n_R)^*, \Delta_{(\mathbb P^n_R)^*} +X^{\univ}_{\varphi})$ is strongly $F$-regular outside $X^{\univ}_{\varphi}$. 
By (d) and inversion of adjunction \cite[Theorem A]{Das15}, 
$(X \times_R (\mathbb P^n_R)^*, \Delta_{(\mathbb P^n_R)^*} +X^{\univ}_{\varphi})$ is purely $F$-regular 
around $X^{\univ}_{\varphi}$. 
To summarise, $(X \times_R (\mathbb P^n_R)^*, \Delta_{(\mathbb P^n_R)^*} +X^{\univ}_{\varphi})$ is purely $F$-regular. 
Hence, its localisation $(X_K, \Delta_K +X^{\gen}_{\varphi})$ is purely $F$-regular. 
Thus, (2) holds. 
\end{proof}

\begin{cor}\label{c-sa-perturb1}
Fix $\mathbb K \in \{\Q, \R\}$. 
Let $R$ be an excellent ring admitting a dualising complex. 
Let $X$ be an integral normal scheme of finite type over $R$, 
$\Delta$ an $\R$-divisor such that $K_X+\Delta$ is $\R$-Cartier, and 
$D$ a semi-ample $\mathbb K$-Cartier $\mathbb K$-divisor. 
Let (P) be one of the following properties. 
\begin{enumerate}
\item Lc. 
\item Klt. 
\end{enumerate}
If $(X, \Delta)$ is (P), 
then there exist $n \in \Z_{>0}$ and 
an effective $\mathbb K$-Cartier $\mathbb K$-divisor $D'$ 
such that $D' \sim_{\mathbb K} D_{R(n)}$ and 
$(X_{R(n)}, \Delta_{R(n)} + D')$ is (P), 
where $X_{R(n)} := X \times_R R(n)$ and 
we denote by $\Delta_{R(n)}$ and $D_{R(n)}$ 
the pullbacks of $\Delta$ and $D$ to  $X_{R(n)}$, respectively. 
\end{cor}

\begin{proof}
We first treat the case when $\mathbb K=\Q$. 
In this case, $|nD|$ is a base point free Cartier divisor for some $n \in \Z_{>0}$. 
Hence we may assume that $D$ itself is a base point free Cartier divisor. 
Then the case (1) directly follows from Theorem \ref{t-Bertini4.1}(1).  
Since $|2D|$ is base point free, 
there exists an $R$-morphism $\varphi : X \to \mathbb P^n_R$ such that 
$\MO_X(2D) \simeq \varphi^*\MO_{\mathbb P^n_R}(1)$. 
For the generic member $X^{\gen}_{\varphi} \sim 2D_{R(n)}$, we set $D' := \frac{1}{2} X^{\gen}_{\varphi} \sim_{\Q} D_{R(n)}$. 
Let us treat the case (2). 
Assume that $(X, \Delta)$ is klt. 
Since $(X_{R(n)}, \Delta_{R(n)})$ is klt and $(X_{R(n)},  \Delta_{R(n)} + X^{\gen}_{\varphi})$ is lc, 
their midpoint $(X_{R(n)}, \Delta_{R(n)} +D')$ is klt. 
Thus the case (2) holds true.

Let us consider the case when $\mathbb K = \R$. 
By definition, we have $D = \sum_{i=1}^{\ell} a_i D_i$ 
for some $a_1, ..., a_{\ell} \in \R_{>0}$ and semi-ample $\Q$-Cartier $\Q$-divisors $D_1, ..., D_{\ell}$. 
Then the problem is reduced to the case when $\mathbb K = \Q$. 
\end{proof}

\begin{cor}\label{c-sa-perturb2}
Fix $\mathbb K \in \{\Q, \R\}$. 
Let $R$ be an $F$-finite noetherian $\F_p$-algebra. 
Let $X$ be an integral normal scheme of finite type over $R$, 
$\Delta$ an $\R$-divisor such that $K_X+\Delta$ is $\R$-Cartier, and 
$D$ a semi-ample $\mathbb K$-Cartier $\mathbb K$-divisor. 
Let (P) be one of the following properties. 
\begin{enumerate}
\item Sharply $F$-pure.
\item Strongly $F$-regular.
\end{enumerate}
If $(X, \Delta)$ is (P), 
then there exist $n \in \Z_{>0}$ and 
an effective $\mathbb K$-Cartier $\mathbb K$-divisor $D'$ 
such that $D' \sim_{\mathbb K} D_{R(n)}$ and 
$(X_{R(n)}, \Delta_{R(n)} + D')$ is (P), 
where $X_{R(n)} := X \times_R R(n)$ and we denote by $\Delta_{R(n)}$ and $D_{R(n)}$ 
the pullbacks of $\Delta$ and $D$ to  $X_{R(n)}$, respectively. 
\end{cor}

\begin{proof}
We can apply the same argument as in Corollary \ref{c-sa-perturb1} 
after replacing Theorem \ref{t-Bertini4.1} by Theorem \ref{t-Bertini4.2}. 
\end{proof}

\begin{rem}
Assume that (P) is (1) or (2) in Corollary \ref{c-sa-perturb1}. 
If $R$ is a field $k$ of characteristic zero, 
Corollary \ref{c-sa-perturb1} is known to hold for (P) even without taking 
the base change $R \subset R(n)$ \cite[Lemma 5.17]{KM98}. 
If $R$ is a field $k$ of positive characteristic, the same assertion holds under the following assumptions \cite[Theorem 1]{Tan17}. 
\begin{enumerate}
\renewcommand{\labelenumi}{(\alph{enumi})}
\item $k$ is an $F$-finite field containing an infinite perfect field of characteristic $p>0$. 
\item $X$ is projective over $k$. 
\item There exists a log resolution of $(X, \Delta)$. 
\end{enumerate}
\end{rem}

\subsection{Examples}

In this subsection, we study the following three examples. 
\begin{enumerate}
\item The Frobenius morphism of $\mathbb P^n_k$ in positive characteristic (Example \ref{e-Frob}). 
\item Fibrations over $\mathbb P^1$ in positive characteristic (Example \ref{e-fib-P1}). 
\item Counterexamples to the Bertini theorem for base point free linear systems in mixed characteristic (Proposition \ref{p-mixed-cex}, Example \ref{e-mixed-cex}). 
\end{enumerate}
The author was not able to find a similar example to (3) in the literature. 

\subsubsection{Positive characteristic}

\begin{ex}\label{e-Frob}
Let $k$ be an algebraically closed field of characteristic $p>0$ 
and set  $X := X':= \mathbb P^n_k$. 
Let $F_{X/k}:X = \mathbb P^n_k \to X'= \mathbb P^n_k$ be the relative Frobenius morphism
\[
F_{X/k} : \mathbb P^n_k \to \mathbb P^n_k, \qquad [a_0: \cdots : a_n] \mapsto [a_0^p : \cdots : a_n^p], 
\]
which is a $k$-morphism. 
In what follows, we use the identifications via the following isomorphisms 
\begin{itemize}
\item $(\mathbb P^n_k)^* \simeq \mathbb P^n_k = \Proj\,k[s_0, ..., s_n]$, and 
\item $\Frac (\mathbb P^n_k)^* \simeq \Frac (\Proj\,k[s_0, ..., s_n]) = k(t_1, ..., t_n)$ with $t_i := s_i /s_0$. 
\end{itemize}
We have the following diagram in which each square is cartesian (Remark \ref{r-functorial}): 
\[
\begin{CD}
X^{\gen}_{F_{X/k}} @>>> X'^{\gen}_{{\rm id}}  @>>> \Spec\,k(n) = \Spec\,k(t_1, ..., t_n)\\
@VVV @VVV @VVV\\
X^{\univ}_{F_{X/k}} @>>> X'^{\univ}_{{\rm id}} @>>> (\mathbb P^n_k)^* \\
@VVV @VVV\\
X \times_k (\mathbb P^n_k)^* @>F_{X/k} \times {\rm id} >> X' \times_k (\mathbb P^n_k)^*\\
@VVV @VVV\\
X @>F_{X/k}>> X'. 
\end{CD}
\]
By 
\[
X'^{\univ}_{{\rm id}} = \{ s_0 x_0 + \cdots + s_nx_n =0 \} \subset 
\Proj\,k[x_0, ..., x_n] \times_k \Proj\,k[s_0, ..., s_n] = \mathbb P^n_k \times_k (\mathbb P^n_k)^*, 
\]
we obtain 
\[
X^{\univ}_{F_{X/k}} = \{s_0x_0^p+ \cdots + s_nx_n^p =0\} \subset 
\Proj\,k[x_0, ..., x_n] \times_k \Proj\,k[s_0, ..., s_n] =\mathbb P^n_k \times_k (\mathbb P^n_k)^*
\]
and 
\[
X^{\gen}_{F_{X/k}} = \{x_0^p+ t_1 x_1^p + \cdots + t_nx_n^p =0\} \subset 
\Proj\,k[x_0, ..., x_n] \times_k k(t_1, ..., t_n) = \mathbb P^n_k \times_k k(t_1, ..., t_n). 
\]
Note that $X^{\univ}_{F_{X/k}}$ is a smooth projective variety over $k$ 
(Theorem \ref{t-Bertini1}(1)) and 
$X^{\gen}_{F_{X/k}}$ is a regular projective variety over $k(n) = k(t_1, ..., t_n)$  (Theorem \ref{t-Bertini2}(4)). 
Although the generic member $X^{\gen}_{F_{X/k}}$ is regular, 
any member of (the linear system corresponding to) the morphism 
$F_{X/k}: \mathbb P^n_k \to \mathbb P^n_k$ is not reduced. 
Indeed, given $\lambda: \Spec\,k \to (\mathbb P^n_k)^* \xrightarrow{\simeq} \Proj\,k[s_0, ..., s_n]$ and $\lambda(\Spec\,k) = [a_0 : \cdots : a_n]$, the corresponding member 
\[
X^{\lambda}_{F_{X/k}} := X^{\univ}_{F_{X/k}} \times_{(\mathbb P^n_k)^*, \lambda} \Spec\,k
\]
is a non-reduced hypersurface of degree $p$: 
\[
X^{\lambda}_{F_{X/k}}=\{ a_0x^p_0+ \cdots + a_n x^p_n = 0 \} \subset \mathbb P^n_k. 
\]
\end{ex}

\begin{ex}\label{e-fib-P1}
We work over an algebraically closed field of characteristic $p>0$. 
Let $X$ be a smooth projective variety over $k$ and let $\pi:X \to \mathbb P^1_k$ be a surjective $k$-morphism. 
In this case, the induced morphism $\alpha : X^{\univ}_{\pi} \to X$ is an isomorphism (Theorem \ref{t-Bertini1}(1)). 
Furthermore, $\pi : X \to \mathbb P^1_k$ and the projection $\rho : X^{\univ}_{\pi} \to  (\P^1_k)^*$ are the same up to isomorphisms, because 
two morphisms $\pi \circ \alpha : X^{\univ}_{\pi} \to \P^1_k$ and $\rho : X^{\univ}_{\pi} \to (\P^1_k)^*$ 
are exchanged by using $X^{\univ}_{\pi} \to \P^1_k \times_k (\P^1_k)^*$ and 
the involution of $\P^1_k \times_k (\P^1_k)^*$ switching $\P^1_k$ and $(\P^1_k)^* (\simeq \P^1_k)$. 
Therefore, the generic member $X^{\gen}_{\pi}$ is isomorphic to the generic fibre 
$X \times_{\mathbb P^1_k} \Spec\,K(\mathbb P^1_k)$ 
of $\pi$. 
Since $X$ is regular, so is $X^{\gen}_{\pi}$. 
For example, if $\pi:X \to \mathbb P^1$ is a quasi-elliptic fibration (see \cite[Definition 7.6]{Bad01} for the definition), then 
its generic fibre is a regular projective curve of arithmetic genus one whose base change to the algebraic closure 
$\overline{K(\mathbb P^1)}$ is a cuspidal cubic curve. 
\end{ex}

\subsubsection{Mixed characteristic}

\begin{prop}\label{p-mixed-cex}
Let $p$ be a prime number and let $R$ be a discrete valuation ring 
of mixed characteristic $(0, p)$ whose residue field is a perfect field. 
Set 
\begin{itemize}
\item 
$\mathbb P^2_{x_0, x_1, x_2} := \Proj\,R[x_0, x_1, x_2] (=\mathbb P^2_{R}),$ 
\item 
$\mathbb P^2_{y_0, y_1, y_2} := \Proj\,R[y_0, y_1, y_2] (=\mathbb P^2_{R}),$ and  
\item 
$X := \{ x_0y_0^p + x_1y_1^p + x_2y_2^p= 0\} \subset \mathbb P^2_{x_0, x_1, x_2} \times_{R} \mathbb P^2_{y_0, y_1, y_2}$. 
\end{itemize} 
Then the following hold. 
\begin{enumerate}
\item $X$ is regular. 
\item 
Let $H$ be a hyperplane of $\mathbb P^2_{x_0, x_1, x_2}$ over $R$. 
For the induced morphism 
\[
\varphi : X \to \mathbb P^2_{x_0, x_1, x_2} =\mathbb P^2_{R}, 
\]
the scheme-theoretic inverse image $\varphi^{-1}(H)$ is not regular. 
\end{enumerate}
\end{prop}


\begin{proof}
Let $\pi$ be a uniformiser of $R$ and set $\m := \pi R$ and $\kappa := R/\m$. 
Let $v_{\pi} : R \to \Z \cup \{ \infty \}$ be the additive $\pi$-adic valuation with $v_{\pi}(\pi)=1$. 

Let us show (1). 
By the Jacobian criterion for smoothness, 
we can check that $X \times_{R} \kappa$ is smooth over $\kappa$. 
Since an effective Cartier divisor $X \times_{R} \kappa$ on $X$ is regular, 
$X$ is regular around $X \times_{R} \kappa$. 
The non-regular locus $X_{{\rm non-reg}}$ of $X$ is a closed subset of $X$ because $R$ is an excellent ring \cite[Scholie 7.8.3(iii)]{EGAIV2}. 
Therefore, if $X_{{\rm non-reg}} \neq \emptyset$, then 
$X_{{\rm non-reg}} 
\cap (X \times_{R} \kappa) \neq \emptyset$, which is a contradiction. 
Hence $X$ is regular. 
Thus (1) holds. 

Let us show (2). 
We can write $H =\{ a_0x_0+a_1x_1+a_2x_2=0\} \subset \mathbb P^2_{x_0, x_1, x_2}$ for some $a_0, a_1, a_2 \in R$. 
Since $H$ is flat over $R$ (cf. Subsection \ref{ss-notation}(\ref{ss-n-hyperplane})), 
we have that $v_{\pi}(a_0) =0, v_{\pi}(a_1)=0$, or $v_{\pi}(a_2)=0$. 
Hence we may assume that $a_0=-1$, i.e. $H = \{ a_1x_1+ a_2x_2 =x_0\}$. 
Since $\kappa$ is a perfect field, we can find $a'_1, a'_2 \in R$ such that 
$v_{\pi}((-a'_1)^p + a_1) \geq 1$ and $v_{\pi}( (-a'_2)^p + a_2) \geq 1$. 
We may assume that $v_{\pi}( (-a'_2)^p + a_2) \geq  v_{\pi}((-a'_1)^p + a_1) \geq 1$. 

We then obtain 
\[
\varphi^{-1}(H) = \{ (a_1x_1+ a_2x_2)y_0^p + x_1y_1^p + x_2y_2^p=0 \} \subset 
\mathbb P^1_{x_1, x_2} \times_R \mathbb P^2_{y_0, y_1, y_2}. 
\]
Consider the following open subset of $\varphi^{-1}(H)$ 
\[
\varphi^{-1}(H) \cap D_+(x_2) \cap D_+(y_0) = \{ f := x_1y_1^p + y_2^p + a_1x_1+ a_2 =0\} 
\subset \mathbb A^3_{x_1, y_1, y_2},  
\]
where the resulting polynomial $x_1y_1^p + y_2^p + a_1x_1+ a_2$ is obtained by the substitution 
$ x_2=y_0=1$. 
For $y'_1 := y_1+a'_1$ and $y'_2 := y_2+a'_2$, the following holds: 
\begin{eqnarray*}
f&=& x_1y_1^p + y_2^p + a_1x_1+ a_2 \\
&=& x_1(y'_1 -a'_1)^p + (y'_2 - a'_2)^p + a_1x_1+ a_2 \\
&=& x_1(y'_1 -a'_1)^p -x_1( -a'_1)^p + (y'_2 - a'_2)^p- ( -a'_2)^p\\
&+& (( -a'_1)^p+a_1)x_1+ ( ( -a'_2)^p+ a_2) \\
& \equiv & (( -a'_1)^p+a_1)x_1+ ( ( -a'_2)^p+ a_2)\mod I^2, 
\end{eqnarray*}
where $I := (\pi, y'_1, y'_2)$. 
We treat the following two cases separately: 
\begin{enumerate}
\item[(i)] $ v_{\pi}((-a'_1)^p + a_1) \geq 2$. 
\item[(ii)] $ v_{\pi}((-a'_1)^p + a_1) =1$. 
\end{enumerate}

(i) Assume $v_{\pi}((-a'_1)^p + a_1) \geq 2$. 
We then obtain $v_{\pi}( (-a'_2)^p + a_2) \geq  v_{\pi}((-a'_1)^p + a_1) \geq 2$. 
In particular, we get $(-a'_1)^p + a_1 \in I^2$ and $(-a'_2)^p + a_2 \in I^2$. 
Hence it holds that $f \in \m^2$ for the maximal ideal $\m := (\pi, x_1, y'_1, y'_2)$ of $R[x_1, y_1, y_2]$. 
Then $\varphi^{-1}(H)$ is not regular at $\m$ in this case. 

(ii) 
Assume $v_{\pi}( (-a'_1)^p + a_1) =1$. 
We can write $(-a'_1)^p + a_1 =\pi a''_1$ and $ (-a'_2)^p + a_2 =\pi a''_2$ 
for some $a''_1 \in R^{\times}$ and $a''_2 \in R$. 
We then have 
\begin{eqnarray*}
f & \equiv & (( -a'_1)^p+a_1)x_1+ ( ( -a'_2)^p+ a_2)\\
&=& \pi ( a''_1x_1+a''_2) \\
&\equiv & 0 \mod \mathfrak n^2
\end{eqnarray*}
for the maximal ideal $\mathfrak n := (\pi, a''_1x_1+a''_2, y'_1, y'_2)$ of $R[x_1, y_1, y_2]$. 
Hence, $\varphi^{-1}(H)$ is not regular at $\mathfrak n$. 
Thus (2) holds. 
\end{proof}

\begin{ex}\label{e-mixed-cex}
We use the same notation as in Proposition \ref{p-mixed-cex}. 
Let us compute the generic member $X^{\gen}_{\varphi}$ of 
\[
\varphi: X =  \{ x_0y_0^p + x_1y_1^p + x_2y_2^p= 0\} \to \mathbb P^2_{x_0, x_1, x_2}. 
\]
Under the identification $(\mathbb P^2_{x_0, x_1, x_2})^* \simeq \mathbb P^2_{s_0, s_1, s_2}$, 
we have that  
\[
X^{\univ}_{\varphi} = \{ x_0y_0^p + x_1y_1^p + x_2y_2^p =0\} 
\cap \{ s_0x_0+s_1x_1+s_2x_2 =0\} 
\]
\[
\subset \mathbb P^2_{x_0, x_1, x_2} \times_R \mathbb P^2_{y_0, y_1, y_2} \times_R \mathbb P^2_{s_0, s_1, s_2}. 
\]
For  the purely transcendental extension  $R(t_1, t_2)$ of $R$ 
with $t_1 := s_1/s_0, t_2 := s_2/s_0$, we have 
\[
X^{\gen}_{\varphi} \simeq \{ -(t_1x_1+t_2x_2)y_0^p + x_1y_1^p + x_2y_2^p =0\} \subset (\mathbb P^1_{x_1, x_2} \times_R \mathbb P^2_{y_0, y_1, y_2}) \times_R R(t_1, t_2). 
\]
Note that $R(t_1, t_2)$ is a discrete valuation ring of mixed characteristic $(0, p)$ 
(Theorem \ref{t-tr-ext-maximal}(3), Theorem \ref{t-tr-bc-conne}(1), Theorem \ref{t-tr-ext-dim}). 
Although $\pi R(t_1, t_2)$ is a maximal ideal of $R(t_1, t_2)$ (Theorem \ref{t-tr-ext-maximal}(3)), 
its residue field $R(t_1, t_2)/\pi R(t_1, t_2) \simeq (R/\pi R)(s'', t'')$ is no longer a perfect field 
(Proposition \ref{p-tr-ext-residue}). 
\end{ex}

\begin{rem}\label{r-mixed-cex}
We use the same notation as in Proposition \ref{p-mixed-cex}. 
Assume that $R := \Z_p$. 
Proposition \ref{p-mixed-cex} claims that 
there exists no hyperplane $H \subset \mathbb P^2_{\Z_p}$ over $\Z_p$ such that $\varphi^{-1}(H)$ is regular for the morphism $\varphi : X \to \mathbb P^2_{\Z_p}$ as in Proposition \ref{p-mixed-cex}. 
It is tempting to hope that we can find  such a hyperplane  $H \subset \mathbb P^2_{R'}$  over $R'$ such that $(\varphi \times_{\Z_p} R')^{-1}(H)$ is regular if we take a suitably  larger local field $R'$,  
which is the integral closure of $\Z_p$ in $K$ for some field extension $\Q_p \subset K$ of finite degree. 
However, Proposition \ref{p-mixed-cex}(2) claims that this is impossible. 
Example \ref{e-mixed-cex} shows that we can find such a member if it is allowed to take 
purely transcendental extensions. 
\end{rem}

\begin{rem}
If $R$ is a noetherian local domain and $X$ is a regular closed subscheme of $\mathbb P^n_R$, 
then it is known that there exists a hypersurface $H$ over $R$ such that $X \cap H$ is regular 
\cite[Theorem 2.15]{BMP}. 
\end{rem}

\section{Bertini theorems over fields}\label{s-density}

\subsection{Density of good members}

If a linear system $V \subset H^0(X, L)$, which is a finite-dimensional $k$-vector subspace, is base point free, 
then we have observed in Section \ref{s-main} that its generic member  has good singularities. 
The purpose of this subsection is to prove that many members of 
$V \otimes_k k' \subset H^0(X_{k'}, L_{k'})$ has a similar property 
for a purely transcendental field extension $k \subset k'$ with ${\rm tr.deg}_k k' \gg 0$, 
where $X_{k'} :=X \times_k k'$ and $L_{k'}$ denotes the pullback of $L$ to $X_{k'}$ 
(Theorem \ref{t-dense-good}, Theorem \ref{t-d-Bertini234}). 
The key result is the following lemma.

\begin{lem}\label{l-dense-affine}
Set $Z := \mathbb A^n_k = \Spec\,k[z_1, ..., z_n]$. 
Let $k \subset k'$ be a purely transcendental extension (whose transcendental degree is not necessarily finite). 
Assume that one of the following holds. 
\begin{enumerate}
\item[(a)] $k$ is an infinite field and ${\rm tr.deg}_k k' \geq n$. 
\item[(b)] ${\rm tr.deg}_k k' \geq n+1$.
\end{enumerate}
For $Z_{k'} := Z \times_k k'$, we set 
\[
\Lambda := \{ \lambda \in Z_{k'}(k')\,|\, {\rm (I)} \text{ and }{\rm (II)}\text{ hold for }\lambda\}. 
\]
\begin{enumerate}
\item[(I)] The image of $\lambda$ to $Z$ is the generic point of $Z$. 
\item[(II)] The field extension $K(Z) \subset \kappa(\lambda)$ induced by (I) is a purely transcendental extension.  
\end{enumerate}
Then $\Lambda$ is a dense subset of $Z_{k'}(k')$ with respect to the Zariski topology. 
\end{lem}

\begin{proof}
Assume (a). 
Fix a transcendental basis $\{ t_i \}_{i \in I}$ for the purely transcendental extension $k \subset k'$, i.e. $k' = k( \{ t_i \}_{i \in I})$ and $\{ t_i \}_{i \in I}$ is algebraically independent over $k$. 
By $|I| = {\rm tr.deg}_k k' \geq n$, we can take distinct elements $t_{i_1}, ..., t_{i_n}$. Set 
\[
\Lambda':= \{ (z_1 -c_1 t_1, ..., z_n - c_nt_n) \in \Spm k'[z_1, ..., z_n]\,|\, (c_1, ..., c_n) \in (k^{\times})^n \}. 
\]
Then $\Lambda' \subset Z_{k'}(k')$. 
It suffices to show that 
\begin{enumerate}
\item[(i)] If $\lambda \in \Lambda'$, then (I) holds for $\lambda$. 
\item[(ii)] If $\lambda \in \Lambda'$, then (II) holds for $\lambda$.  
\item[(iii)] $\Lambda'$ is dense in $Z_{k'}(k')$. 
\end{enumerate}
Indeed, (i) and (ii) imply that $\Lambda' \subset \Lambda$. 
Then, combining with (iii), we see that $\Lambda$ is dense in  $Z_{k'}(k')$.

Fix $\m \in \Lambda'$. 
We can write $\m = (z_1 - c_1 t_{i_1}, ..., z_n - c_nt_{i_n})$ 
for some $(c_1, ..., c_n) \in (k^{\times})^n$. 
Let us show (i). 
To this end, it suffices to prove that $\m \cap k[z_1, ..., z_n] = \{0\}$. 
Pick an element $f(z_1, ..., z_n) \in \m \cap k[z_1, ..., z_n]$. Then 
\[
f(z_1, ..., z_n) = \sum_{j=1}^n (z_j - c_j t_{i_j}) g_j(z_1, ..., z_n) 
\]
for some $g_j(z_1, ..., z_n) \in k'[z_1, ..., z_n]$. 
Substituting $z_1 = c_1t_{i_1}, ..., z_n = c_nt_{i_n}$, we obtain the equation 
$f(c_1t_{i_1}, ..., c_nt_{i_n}) =0$ in $k'$. 
By $f(z_1, ..., z_n) \in k[z_1, ..., z_n]$ and the fact that $c_1t_1, ..., c_nt_n$ are algebraically independent over $k$, 
we obtain $f(z_1, ..., z_n) =0$ in $k[z_1, ..., z_n]$, 
which implies $\m \cap k[z_1, ..., z_n] = \{0\}$. 
Thus (i) holds. 

Let us show (ii). 
By (i), we have the induced field extension: 
\[
K(Z)=k(z_1, ..., z_n) \hookrightarrow \kappa(\m) = k'[z_1, ..., z_n]/\m \xrightarrow{\simeq} k', \qquad z_j \mapsto c_j t_{i_j}.
\]
Since $c_1t_{i_1}, ..., c_nt_{i_n}$ is a subset of a transcendental basis for the purely transcendental extension $k \subset k'$, 
the induced field extension $K(Z) \hookrightarrow \kappa(\m)$ is a purely transcendental extension. 
Thus (ii) holds. 

Let us show (iii). 
Fix a proper closed subset $W$ of $Z_{k'}(k')$. 
Since $k^{\times}$ is an infinite set, there exists $c_1 \in k^{\times}$ such that 
\[
\{c_1t_{i_1} \} \times_{k'} \mathbb A_{k'}^{n-1} \not\subset W. 
\]
By induction on $n$, we see that $(c_1t_{i_1}, ..., c_nt_{i_n}) \not\in W$ for some $(c_1, ..., c_n) \in (k^{\times})^n$. 
Hence (iii) holds. 
This completes the proof for the case when (a) holds.

Assume (b). 
Fix a transcendental basis $\{ t_i \}_{i \in I}$ for the purely transcendental extension $k \subset k'$. 
By $|I| = {\rm tr.deg}_k k' \geq n+1$, we can take distinct elements $t_{i_0}, t_{i_1}, ..., t_{i_n}$. 
Set 
\[
\Lambda':= \{ (z_1 - t_{i_1} -t_{i_0}^{d_1}, ..., z_n - t_{i_n} -t_{i_0}^{d_n}) \in \Spm k'[z_1, ..., z_n]\,|\, (d_1, ..., d_n) \in (\Z_{>0})^n \}. 
\]
Then the same argument as in (a) works after replacing $\Lambda'$ by $\Lambda''$. 
\end{proof}

\begin{nota}\label{n-dense-good}Let $k \subset k'$ be a purely transcendental field extension. 
Assume that one of the following holds. 
\begin{enumerate}
\item[(a)] $k$ is an infinite field and ${\rm tr.deg}_k k' \geq n$. 
\item[(b)] ${\rm tr.deg}_k k' \geq n+1$.
\end{enumerate}
Take a $k'$-rational point $\lambda \in (\mathbb P^n_{k'})^*(k')$ and 
consider the morphism 
\[
\widetilde{\lambda} :  \Spec\,k' \xrightarrow{\lambda} (\mathbb P^n_{k'})^* \xrightarrow{\pi} (\mathbb P^n_{k})^*, 
\]
where $\lambda : \Spec\,k' \to (\mathbb P^n_{k'})^*$ denotes the closed immersion 
whose image is $\lambda$ and the latter morphism $\pi$ is the natural projection. 
We set 
\[
\Lambda := \{ \lambda \in (\mathbb P^n_{k'})^*(k')\,|\, (I), (II)\}.
\]
\begin{enumerate}
\item[(I)] The image of $\lambda$ to $(\mathbb P^n_{k})^*$ is equal to the generic point of  $(\mathbb P^n_{k})^*$.
\item[(II)] The field extension $K( (\mathbb P^n_{k})^*) \subset \kappa(\lambda)$ induced by (I) is a purely transcendental extension. 
\end{enumerate}

Let $X$ be a scheme of finite type over $k$ and 
let $\varphi : X \to \mathbb P^n_k$ be a $k$-morphism. 
We define $X^{\lambda}_{\varphi}$ 
as the closed subshceme of $X_{k'} := X \times_k k'$ corresponding to $\widetilde{\lambda}$, 
which completes the following cartesian diagram: 
\[
\begin{CD}
X^{\lambda}_{\varphi} @>>> X^{\univ}_{\varphi}\\
@VVV @VVV\\
\Spec\,k' @>\widetilde{\lambda} >> (\mathbb P^n_{k})^*. 
\end{CD}
\]
Set $X^{\widetilde{\lambda}}_{\varphi} := X^{\lambda}_{\varphi}$. 
\end{nota}

\begin{thm}\label{t-dense-good}
We use Notation \ref{n-dense-good}. 
Then the following hold. 
\begin{enumerate}
\item 
If $\lambda \in \Lambda$, then $X^{\lambda}_{\varphi} \simeq X^{\gen}_{\varphi} \times_{K( (\mathbb P^n_{k})^*)} \kappa(\lambda)$. 
\item 
$\Lambda$ is a dense subset of $(\mathbb P^n_{k'})^*(k')$ with respect to the Zariski topology. 
\end{enumerate}
\end{thm}

\begin{proof}
The assertion (1) follows from the fact that the right and big squares in the following diagram are cartesian  
\[
\begin{CD}
X^{\lambda}_{\varphi} @>>> X^{\gen}_{\varphi} @>>> X^{\univ}_{\varphi}\\
@VVV @VVV @VVV\\
\Spec\,\kappa(\lambda) @>>> \Spec\,K(\mathbb P^n_k) @>>> \mathbb P^n_k. 
\end{CD}
\]
The assertion (2) holds by Lemma \ref{l-dense-affine}. 
\end{proof}

\begin{thm}\label{t-d-Bertini234}
We use Notation \ref{n-dense-good}. 
For $\lambda \in \Lambda$, 
let $\beta_{\lambda} : X^{\lambda}_{\varphi} \to X$ be the induced morphism. 
Then the following hold. 
\begin{enumerate}
\item 
Let (P) be one of the properties listed in Theorem \ref{t-Bertini2}.
If $X$ is (P) and $\lambda \in \Lambda$, then $X_{\varphi}^{\lambda}$ is (P). 
\item 
Assume that $X$ is a normal variety over $k$. 
Let $\Delta$ be an $\R$-divisor on $X$ such that $K_X+\Delta$ is $\R$-Cartier 
and let $\Delta_{k'}$ be the pullback of $\Delta$ to $X_{k'}$. 
\begin{enumerate}
\item[(2.a)]  
Let (Q) be one of the properties listed in Theorem \ref{t-Bertini3.1}.
If $(X, \Delta)$ is (Q) and $\lambda \in \Lambda$, then $(X_{\varphi}^{\lambda}, \beta_{\lambda}^*\Delta)$ is (Q).  
\item[(2.b)] 
Let (Q)' be one of klt and lc. 
If $(X, \Delta)$ is (Q)' and $\lambda \in \Lambda$, then 
$(X_{k'}, \Delta_{k'} + X_{\varphi}^{\lambda})$ is (Q), 
\end{enumerate}
\item 
Assume that $k$ is an $F$-finite field of characteristic $p>0$ 
and $X$ is a normal variety over $k$. 
Let $\Delta$ be an effective $\R$-divisor on $X$. 
\begin{enumerate}
\item[(3.a)] 
Let (R) be one of the properties listed in Theorem \ref{t-Bertini3.2}.
If $(X, \Delta)$ is (R) and $\lambda \in \Lambda$, then $(X_{\varphi}^{\lambda}, \beta_{\lambda}^*\Delta)$ is (R).
\item[(3.b)] 
Let (R)' be one of sharply $F$-pure and strongly $F$-regular. 
If $(X, \Delta)$ is (R)' and $\lambda \in \Lambda$, then 
$(X_{k'}, \Delta_{k'} + X_{\varphi}^{\lambda})$ is (R)'.  
\end{enumerate}
\end{enumerate}
Here $\Delta_{k'}$ denotes the pullback of $\Delta$ to $X_{k'}$. 
\end{thm}

\begin{proof}
The assertion (1) follows from Theorem \ref{t-Bertini2}, 
Theorem \ref{t-dense-good}, and the fact that (P) is stable 
under taking base changes by purely transcendental field extensions. 
Indeed, if ${\rm tr.deg}_k k' <\infty$, then the composite morphism 
\[
X^{\lambda}_{\varphi} \to  X^{\gen}_{\varphi} \to X
\] 
is an $\mathbb A^{\ell}$-localising morphism for some $\ell \in \Z$ (Proposition \ref{p-A-local}(5)), and hence Theorem \ref{t-Bertini2} holds even after replacing 
$X^{\gen}_{\varphi}$ by $X^{\lambda}_{\varphi}$. 
The remianing case, i.e. ${\rm tr.deg}_k k' =\infty$, is reduced to this case 
because we can write $k' = \bigcup_{i \in I} k_i$ for  
some filtered set $\{k_i\}_{i \in I}$ consisting of purely transcendental extensions $k \subset k_i$ 
of finite transcendental degree.

The other assertions (2.a), (3.a), (2.b), and (3.b) 
hold by the same argument 
after replacing Theorem \ref{t-Bertini2} by 
Theorem \ref{t-Bertini3.1}, Theorem \ref{t-Bertini3.2}, Theorem \ref{t-Bertini4.1}, 
and Theorem \ref{t-Bertini4.2}, respectively. 
\end{proof}

\begin{rem}\label{r-dense-good}
We use Notation \ref{n-dense-good}. 
Assume that ${\rm tr.deg}_k k' <\infty$. 
For $\lambda \in \Lambda$,  the composite morphism 
\[
X^{\lambda}_{\varphi} \to  X^{\gen}_{\varphi} \to X
\] 
is an $\mathbb A^{\ell}$-localising morphism for some $\ell \in \Z$ (Proposition \ref{p-A-local}(5)). 
In particular, Theorem \ref{t-Bertini1}(3) holds even after replacing 
$X^{\gen}_{\varphi}$ by $X^{\lambda}_{\varphi}$. 
\end{rem}

\subsection{Generic members for linear systems}

The purpose of this subsection is to recall the definitions and some fundamental properties of the universal families $X^{\univ}_{L, V}$ and 
the generic members $X^{\gen}_{L, V}$ of linear systems. 
For flexibility, we consider a linear map $V \to H^0(X, L)$ instead of 
a usual linear system $V \subset H^0(X, L)$, 
which makes it easier to establish some functorial statements (cf. Proposition \ref{p-base-change}). 
Comparison to $X^{\univ}_{\varphi}$ and $X^{\gen}_{\varphi}$ introduced in Section \ref{s-main}, 
an advantage of $X^{\univ}_{L, V}$ and $X^{\gen}_{L, V}$ is to allow us to treat linear sytems with base points (cf. Remark \ref{r-2-gene-defs}). 

\begin{dfn}\label{d-ls-gen-member}
Let $X$ be a scheme of finite type over $k$, $L$ an invertible sheaf, and $\theta:V \to H^0(X, L)$ a $k$-linear map from a finite-dimensional $k$-vector space $V$ with $\overline V := \theta(V) \neq 0$.  
By $\mathbb P(\overline V) ={\rm Proj}\,(S(\overline{V}^*))$ (Subsection \ref{ss-notation}(\ref{sym-prod})), 
we have the universal closed subscheme: 
\[
X^{\univ}_{L, \overline V}
\subset X \times_k \mathbb P(\overline V). 
\]
For the induced dominant morphism $\mathbb P(\theta):\mathbb P(V) \setminus \mathbb P(\Ker(\theta)) \to \mathbb P(\overline{V})$ (cf. Remark \ref{r-proj-induced}), 
we define $X^{\univ}_{L, V}$ as the inverse image of $X^{\univ}_{L, \overline V}$ 
by $X \times (\mathbb P(V) \setminus \mathbb P(\Ker(\theta))) \to X \times_k \mathbb P(\overline V)$, 
i.e. we have the following cartesian squares: 
\[
\begin{CD}
X^{\univ}_{L, V} @>>> X^{\univ}_{L, \overline V}\\
@VVV @VVV\\
X \times_k (\mathbb P(V) \setminus \mathbb P(\Ker(\theta))) @>{\rm id} \times \mathbb P(\theta)>> X \times_k \mathbb P(\overline V)\\
@VVV @VVV\\
\mathbb P(V) \setminus \mathbb P(\Ker(\theta)) @>\mathbb P(\theta)>> \mathbb P(\overline V).
\end{CD}
\]
Let $X^{\gen}_{L, V}$ be the generic fibre of $X^{\univ}_{L, V} \to \mathbb P(V) \setminus \mathbb P(\Ker(\theta))$: 
\[
X^{\gen}_{L, V} := X^{\univ}_{L, V} \times_{\mathbb P(V) \setminus \mathbb P(\Ker(\theta))} \Spec\,K(\mathbb P(V)).
\] 
By definition, we obtain 
\[
X^{\gen}_{L, V} = X^{\gen}_{L, \overline{V}} \times_{K(\mathbb P(\overline{V}))} K(\mathbb P(V)).
\]
\end{dfn}

\begin{rem}\label{r-proj-induced}
Let $f:V \to W$ be a $k$-linear map of finite-dimensional $k$-vector spaces. 
\begin{enumerate}
\item If $f(V) \neq \{0\}$ i.e. $V \neq \Ker f$, then we have the corresponding rational map which is defined on $\mathbb P(V) \setminus \mathbb P(\Ker f)$: 
\[
\mathbb P(f):\mathbb P(V) \dashrightarrow \mathbb P(W), \qquad 
\mathbb P(f):\mathbb P(V) \setminus \mathbb P(\Ker f) \to \mathbb P(W).
\]
\end{enumerate}
Assume that $f(V) \neq \{0\}$. 
\begin{enumerate}
\setcounter{enumi}{1}
\item 
$f$ is injective if and only if $\mathbb P(f):\mathbb P(V) \dashrightarrow \mathbb P(W)$ is a morphism. 
Furthermore, if $f$ is injective, then $\mathbb P(f):\mathbb P(V) \to \mathbb P(W)$ is a closed immersion. 
\item $f$ is surjective if and only if $\mathbb P(f)$ is dominant. 
\end{enumerate}
\end{rem}

\begin{rem}\label{r-2-gene-defs}
Let $X$ be a scheme of finite type over $k$ and 
let $\varphi : X \to \mathbb P^n_k$ be a $k$-morphism. 
Recall that $\varphi$ is defined by $L := \varphi^*\MO_{\mathbb P^n_k}(1)$ and 
$s_0, ..., s_n \in H^0(X, L)$: 
\[
\varphi : X \to \mathbb P^n_k, \qquad x \mapsto [s_0(x) : \cdots : s_n(x)]. 
\]
For an $(n+1)$-dimensional $k$-vector space $V = \bigoplus_{i=0}^n k e_i$ 
with a $k$-linear basis $e_0, ..., e_n$, 
consider the following $k$-linear map 
\[
\theta : V = \bigoplus_{i=0}^n k e_i \to H^0(X, L), \qquad e_i \mapsto s_i. 
\]
We then have $X^{\gen}_{\varphi} = X^{\gen}_{L, V}$. 
\end{rem}

\begin{rem}\label{r-gen-member1}
We use the same notation as in Definition \ref{d-ls-gen-member}. 
Fix a $k$-linear basis $v_0, ..., v_n$ of $V$: 
\[
V = kv_0 \oplus \cdots \oplus k v_n. 
\]
Let $\varphi_0, ..., \varphi_n$ be their images by $\theta:V \to H^0(X, L)$. 
By $\overline V = \theta(V)$, we have $\varphi_0, ..., \varphi_n \in \overline V$. 
For the induced isomorphism $\zeta: V \xrightarrow{\simeq} k^{n+1}, \varphi_i \mapsto e_i$ where $e_0, ...,e_n$ denotes the standard basis, we have the induced isomorphism $\mathbb P(\zeta): \mathbb P(V) \xrightarrow{\simeq} \mathbb P^n = \Proj\,k[t_0, ..., t_n]$. 
Under the identification $\mathbb P(V) = \mathbb P^n_k =\Proj\,k[t_0, ..., t_n]$ via $\mathbb P(\zeta)$, we have 
\begin{itemize}
\item $X_{L, V}^{\univ} = \{ t_0\varphi_0 + t_1 \varphi_1 + \cdots + t_n\varphi_n = 0\} \subset X \times_k (\mathbb P(V) \setminus \mathbb P(\Ker(\theta))$, and 
\item $X_{L, V}^{\gen}  = \{ \varphi_0 + (t_1/t_0) \varphi_1 +\cdots + (t_n/t_0) \varphi_n = 0\} \subset X \times_k K(\mathbb P(V))$. 
\end{itemize}
where 
$K(\mathbb P(V)) = k(t_1/t_0, ..., t_n/t_0)$. 
\end{rem}

\begin{prop}\label{p-base-change}
Let $f:X \to Y$ be a $k$-morphism of schemes which are of finite type over $k$. 
Let $L_Y$ be an invertible sheaf on $Y$ and set $L_X :=f^*L_Y$. 
Let $\theta_Y: V \to H^0(Y, L_Y)$ be a $k$-linear map from a finite-dimensional $k$-vector space $V$, 
so that we have $k$-linear maps: 
\[
\theta_X: V \xrightarrow{\theta_Y} H^0(Y, L_Y) \xrightarrow{f^{\sharp}} H^0(X, L_X). 
\]
Assume that $\theta_X(V) \neq 0$. 
Then the following hold. 
\begin{enumerate}
\item 
We have the following cartesian diagram: 
\[
\begin{CD}
X^{\univ}_{L_X, V} @>>> Y^{\univ}_{L_Y, V}\\
@VVa_X V @VVa_Y V\\
X \times_k (\mathbb P(V) \setminus \mathbb P(\Ker(\theta_X)))  @>f \times_k i >> 
Y \times_k (\mathbb P(V) \setminus \mathbb P(\Ker(\theta_Y))), 
\end{CD}
\]
where $i: \mathbb P(V) \setminus \mathbb P(\Ker(\theta_X)) \hookrightarrow \mathbb P(V) \setminus \mathbb P(\Ker(\theta_Y))$ 
denotes the induced open immersion and both $a_X$ and $a_Y$ are the induced closed immersions. 
\item 
We have the following cartesian diagram: 
\[
\begin{CD}
X^{\gen}_{L_X, V} @>>> Y^{\gen}_{L_Y, V}\\
@VVb_X V @VVb_Y V\\
X \times_k K(\mathbb P(V))  @>f \times_k K(\mathbb P(V)) >> Y \times_k K(\mathbb P(V)), 
\end{CD}
\]
where $b_X$ and $b_Y$ denote the induced closed immersions. 
\end{enumerate}
\end{prop}

\begin{proof}
Taking the fibres over $K(\mathbb P(V))$,  (1) immediately implies (2). 
Therefore, it suffices to show (1). 
Set $n := \dim V -1$. 
We have $n = \dim V-1 \geq \dim \theta_X(V) -1 \geq 0$. 
Let $\overline{V}_X$ and $\overline{V}_Y$ be the images of $V$ on $H^0(X, L_X)$ and $H^0(Y, L_Y)$, respectively. 
Fix a $k$-linear basis $v_0, ..., v_n$ of $V$.  
Let $\varphi_0, ..., \varphi_n$ and $\psi_0, ..., \psi_n$ be 
the images on $\overline{V}_X$ and $\overline{V}_Y$, respectively. 
Under the identification $\mathbb P(V) = \Proj\,k[t_0, ..., t_n]$ as in Remark \ref{r-gen-member1}, 
we obtain 
\[
X^{\univ}_{L_X, V_X} = \{ t_0\varphi_0 + \cdots + t_n\varphi_n = 0 \} \subset X \times_k (\mathbb P(V) \setminus \mathbb P(\Ker(\theta_X))) 
\]
and 
\[
Y^{\univ}_{L_Y, V_Y} = \{ t_0\psi_0 + \cdots + t_n\psi_n = 0 \} \subset Y \times_k (\mathbb P(V) \setminus \mathbb P(\Ker(\theta_Y))).
\]
By $f^{\sharp}(\psi_0)=\varphi_0, ..., f^{\sharp}(\psi_n)=\varphi_n$, 
the scheme-theoretic inverse image  of the closed subscheme 
$Y^{\univ}_{L_Y, V_Y} =\{ t_0\psi_0 + \cdots + t_n\psi_n = 0 \}$ is equal to 
$X^{\univ}_{L_X, V_X} = \{ t_0\varphi_0 + \cdots + t_n\varphi_n = 0 \}$. 
\end{proof}

\subsection{Generic members vs General members}\label{ss-general-generic}

In this subsection, we summarise which property is stable under taking general members. 

\subsubsection{Characteristic zero} 

We first discuss the case of characteristic zero. 
In this case, almost all the properties are stable under taking general members. 

\begin{thm}\label{t-char-0-summary}
Let $k$ be an algebraically closed field of characteristic zero, $X$ a scheme of finite type over $k$, 
$L$ an invertible sheaf, and $V \subset H^0(X, L)$ a finite-dimensional $k$-vector subspace which is base point free. 
Let (P) be one of (1)--(25) in Theorem \ref{t-Bertini2}, 
Theorem \ref{t-Bertini3.1}, and  Theorem \ref{t-Bertini3.2}  
except for (11)--(12), (15)--(16), and (23)--(25). 
If $X$ is (P), then general members of $V$ are (P) (note that the latter five notions (15)--(16) and (23)--(25) are defined only in positive characteristic). 
\end{thm}

\begin{proof}
It is well known that the cases (1)--(8) hold true. 
As for (9)--(10), see \cite[Theorem 1, Theorem 2, and the paragraph immediately after Theorem 2]{CGM86} 
(note that seminormality coincides with weak normality in characteristic zero). 
It follows from \cite[Lemma 5.17]{KM98} that (17)--(22) hold. 

Let us treat the case (13). 
Let $f: X' \to X$ be a resolution of singularities. 
Assume that $f_*\MO_{X'} = \MO_X$ and $R^if_*\MO_{X'} =0$ for any $i>0$. 
Fix a general member $D$ of $V$, so that $D' := f^*D$ is smooth. 
By an exact sequence 
\[
0 \to \MO_{X'}(-D') \simeq \MO_{X'} \otimes f^*\MO_X(-D) \to \MO_{X'} \to \MO_{D'} \to 0, 
\]
it holds that $f_*\MO_{D'} = \MO_D$ and $R^if_*\MO_{D'}=0$ for $i>0$, i.e. $D$ has rational singularities. 
A similar argument works also for the case (14) by using Schwede's characterisation 
of Du Bois singularities 
(Subsubsection \ref{sss-sing-scheme}(\ref{Schwede-duBois})). 
\end{proof}

\begin{rem}\label{r-char-0-summary}
\begin{enumerate}
\item[(a)]  
Let $X$ be an irreducible scheme of finite type over an algebraically closed field $k$, $L$ an invertible sheaf, and $\theta:V \to H^0(X, L)$ a $k$-linear map from a finite-dimensional $k$-vector space $V$ with $\overline V := \theta(V) \neq 0$.  
If $\overline V$ is base point free and $\dim \overline{\varphi(X)} \geq 2$ 
for the induced morphism $\varphi : X \to \mathbb P(V)$, 
then any general member is irreducible \cite[4) of Theoreme 6.3]{Jou83}. 
It follows from \cite[Theoreme 4.10]{Jou83} that also the generic fibre $X^{\gen}_{L, V}$ is geometrically irreducible. 
\item[(b)]  
Assume that (P) is one of (11) and (12). 
Then (P) is not stable under taking general hyperplane sections. 
Indeed, if $X$ is a smooth projective curve and $L$ is a very ample invertible sheaf with $\deg L \geq 2$, 
then general members of $|L|$ are not connected. 
\end{enumerate}
\end{rem}

\subsubsection{Positive characteristic}

Let $k$ be an algebraically closed field of characteristic $p>0$, 
$X$ a scheme of finite type over $k$, 
$L$ an invertible sheaf, and $V \subset H^0(X, L)$ is a finite-dimensional $k$-vector subspace which is base point free. 
In Table \ref{table-char-p}, we use the following notation. 

\begin{table}
\caption{Bertini in positive characteristic}\label{table-char-p}
\begin{tabular}{|c|c|c|c|c|} \hline
\# &  & \hspace{3mm} {\rm Generic}\hspace{3mm} & \text{Very ample} & 
\text{Base point free} \\ \hline \hline
(1) & $R_m$ & Yes &  Yes & No \\ \cline{2-4} 
\hline 
(2) &$S_m$ & Yes &  Yes & Yes \\ \cline{2-4} 
\hline 
(3) & $G_m$ & Yes &  Yes & Yes \\ \cline{2-4} 
\hline 
(4) &Regular& Yes &  Yes & No \\ \cline{2-4} 
\hline 
(5) &Cohen--Macaulay  & Yes &  Yes & Yes \\ \cline{2-4} 
\hline 
(6) &Gorenstein  & Yes &  Yes & Yes \\ \cline{2-4} 
\hline 
(7) &Reduced  & Yes &  Yes & No \\ \cline{2-4} 
\hline 
(8) &Normal & Yes &  Yes & No \\ \cline{2-4} 
\hline 
(9) &Seminormal& Yes &  No & No \\ \cline{2-4} 
\hline 
(10) &Weakly normal  & Yes &  No & No \\ \cline{2-4} 
\hline 
(11) &Irreducible or empty  & Yes &  No & No \\ \cline{2-4} 
\hline 
(12) &Integral or empty  & Yes &  No & No \\ \cline{2-4} 
\hline 
(13) &{Resolution-rational} & Yes & ?  & No \\ \cline{2-4} 
\hline 
(15) &$F$-rational  & Yes & ? &  No  \\ \cline{2-4} 
\hline 
(16) &$F$-injective  & Yes & No &  No  \\ \cline{2-4} 
\hline 
(17) &Terminal  & Yes &  ? & No \\ \cline{2-4} 
\hline 
(18) &Canonical  & Yes &  ? & No \\ \cline{2-4} 
\hline 
(19) &Klt  & Yes &   ? & No \\ \cline{2-4} 
\hline 
(20) &Plt & Yes &  ? & No \\ \cline{2-4} 
\hline 
(21) &Dlt & Yes &  ? & No \\ \cline{2-4} 
\hline 
(22) &Lc & Yes &  ? & No \\ \cline{2-4} 
\hline 
(23) & Strongly $F$-regular  & Yes & Yes & No   \\ \cline{2-4} 
\hline 
(24) & Sharply $F$-pure  & Yes & Yes &  No  \\ \cline{2-4} 
\hline 
(25) & $F$-pure  & Yes & ? & No  \\ \cline{2-4} 
\hline 
  \end{tabular}
\end{table}
\begin{enumerate}
\renewcommand{\labelenumi}{(\roman{enumi})}
\item 
The column \lq\lq Generic\rq\rq \, gives an answer to whether the property is stable under taking the generic member 
$X^{\gen}_{L, V}$. 
\item 
The column \lq\lq Very ample\rq\rq \, gives an answer to whether the property is stable under taking general members of $V$ 
if $L$ is very ample and $V=H^0(X, L)$.
\item 
The column \lq\lq Base point free\rq\rq \, gives an answer to whether the property is stable under taking general members of 
$V$. 
\item The question mark \lq\lq  ?\rq\rq \, means that the question is open as far as the author knows. 
\item For (17)--(25), we consider the corresponding questions for pairs $(X, \Delta)$, 
where $X$ is a normal variety over $k$ and $\Delta$ is an $\R$-divisor on $X$ ($\Delta$ is assumed to be effective for (23)--(25)). 
\end{enumerate}

\begin{proof}(of Table \ref{table-char-p}) 
Concerning the column \lq\lq Generic\rq\rq, we may apply  Theorem \ref{t-Bertini2}, Theorem \ref{t-Bertini3.1}, and Theorem \ref{t-Bertini3.2}. 
Hence let us consider the remaining columns \lq\lq Very ample\rq\rq and \lq\lq Base point free\rq\rq. 

We first treat the column \lq\lq Base point free\rq\rq. 
All of \lq\lq No\rq\rq are guaranteed by the existence of a base point free linear system 
on a smooth variety such that an arbitrary member is not geometrically reduced (cf. Example \ref{e-Frob}). 
It is well known that (5) and (6) hold true. 
After removing suitable closed subsets, (2) and (3) follow from (5) and (6), respectively. 

We now handle the column \lq\lq Very ample\rq\rq. 
We have already proven the cases (2) and (3). 
(4) is nothing but the classical Bertini theorem. 
Then (4) implies (1). 
Then each of (5)--(7) is a combination of (1)--(3). 
As for (9) and (10), see  \cite[Remark 2.6]{CGM83} and \cite[Corollary 1 and Corollary 4]{CGM89}, respectively. 
For (11)--(12), see Remark \ref{r-char-0-summary}. 
It follows from \cite[Theorem 7.5]{SZ13} that (16) holds true. 
Concerning (23)--(24), 
we may assume that $\Delta$ is a $\Q$-divisor by enlarging the coefficients of $\Delta$. 
Then we may apply \cite[Theorem 6.1]{SZ13}. 
\end{proof}

\begin{rem}
Contrary to the situation in characteristic zero, 
it is not known whether singularities {in the} 
minimal model program (i.e. (17)--(22)) 
are stable under taking general hyperplane sections as written in Table \ref{table-char-p}. 
In the three-dimensional case, \cite{ST20} has settled some of the cases. 
\end{rem}

\begin{rem}
As for (25), 
$F$-purity is stable under taking very general hyperplane sections 
because $(X, (1-\frac{1}{n})\Delta)$ is sharply $F$-pure for any $n \in \Z_{>0}$. 
\end{rem}

\subsection{Base points}

\begin{prop}\label{p-base-pt}
Let $X$ be a scheme of finite type over $k$, $L$ an invertible sheaf, and $\theta:V \to H^0(X, L)$ a $k$-linear map from a finite-dimensional $k$-vector space $V$ with $\overline V := \theta(V) \neq 0$.  
Set $K:=K(\mathbb P(V))$ and $X_K := X \times_k K$. 
Then the following hold. 
\begin{enumerate}
\item Fix a closed point $x \in X$ and let $x'$ be its inverse image to $X \times_k K$. 
Note that $x'$ is a closed point of $X_K$ (cf. Lemma \ref{l-tr-ext-normaln}). 
Then $x \in \Bs(\overline V)$ if and only if $x' \in X^{\gen}_{L, V}$. 
\item 
If $\overline V$ is base point free, then $Y \times_k K \not\subset X^{\gen}_{L, V}$ 
for any non-empty closed subscheme $Y$ of $X$ (in particular, $X^{\gen}_{L, V}$ does not contain the inverse image 
of an arbitrary closed point on $X$). 
\item 
If $k$ is algebraically closed and $\overline V$ is not base point free, then 
$X^{\gen}_{L, V}$ contains a $K(\mathbb P(V))$-rational point. 
%
\end{enumerate}
\end{prop}

\begin{proof}
Let us show (1). 
It is clear that $x \in \Bs(\overline V)$ implies $x' \in X^{\gen}_{L, V}$. 
Assuming that $x \not\in \Bs(\overline V)$, let us show $x' \not\in X^{\gen}_{L, V}$. 
Replacing $X$ by $X \setminus \Bs(\overline V)$, 
the problem is reduced to the case when $\overline V$ is base point free, and hence 
$X^{\gen}_{L, V} = X^{\gen}_{\varphi}$ for the morphism $\varphi : X \to \mathbb P^n_k$ 
induced by the linear system $\overline V$ (Remark \ref{r-2-gene-defs}). 
We then obtain $x' \not\in X^{\gen}_{L, V}$ by Theorem \ref{t-gene-avoid}. 
This completes the proof of (1). 

The assertions (2) and (3) immediately follow from (1). 
\end{proof}

\begin{rem}\label{r-size}
We use the same notation as in Proposition \ref{p-base-pt}. 
Assume that $X$ is integral. 
By $\overline{V} \neq 0$,  we have $\Bs(\overline V) \subsetneq X$. 
Therefore, $X^{\gen}_{L, V}$ is an effective Cartier divisor on $X_K$ (Remark \ref{r-gene-hyperplane}, Proposition \ref{p-base-pt}(1)). 
In particular, it holds that $X^{\gen}_{L, V}$ is either empty or of pure codimension one. 
\end{rem}

\subsection{Picard groups for generic members}

{

\begin{prop}\label{p-LHST-surje}
Let $k$ be a field and let $X$ be a regular projective variety over $k$. 
Let $\varphi : X \to \P^n_k$ be a $k$-morphism, where $n$ is a non-negative integer. 
Then the map 
\[
\beta^* : \Pic\,X \to \Pic\,X^{\gen}_{\varphi}, \qquad L \mapsto \beta^*L
\]
is surjective, where $\beta : X^{\gen}_{\varphi} \to X$ denotes the induced morphism. 
\end{prop}

\begin{proof}
If $n=0$, then 
we get 
$X^{\univ}_{\varphi} = X^{\gen}_{\varphi} = \emptyset$ and $\Pic\,X^{\gen}_{\varphi} = \{0\}$. 
In what follows, we assume $n >0$. 
We have the following commutative diagram in which all the squares are cartesian (Definition \ref{d-gen-member}): 
\[
\begin{tikzcd}
	X^{\gen}_{\varphi}  & X^{\univ}_{\varphi} \\
	X \times_k k(n) & X \times_k (\mathbb P^n_k)^* & X\\
	\Spec\,k(n) & (\mathbb P^n_k)^* & \Spec\,k.
	\arrow[from=3-2, to=3-3]
	\arrow[from=3-1, to=3-2]
	\arrow[from=2-3, to=3-3]
	\arrow[from=2-2, to=3-2, "{\rm pr}_2"']
	\arrow[from=2-1, to=3-1]
	\arrow[from=2-1, to=2-2]
	\arrow[hook, from=1-2, to=2-2, "j"']
	\arrow[hook, from=1-1, to=2-1]
	\arrow[from=1-1, to=1-2, "\gamma"]
	\arrow[from=1-2, to=2-3, "\alpha"]
	\arrow[from=2-2, to=2-3, "{\rm pr}_1"']
 \arrow[from=1-1, to=2-3, "\beta", bend left=15mm]
\end{tikzcd}
\]
Then we get the following commutative diagram:  
\begin{equation}\label{e1-LHST-surje}
\begin{tikzcd}
0 \arrow[r] & \Pic\,X \arrow[r, "\alpha^*"] & \Pic\,X^{\univ}_{\varphi} \arrow[r, "\cdot \ell"] & {[\kappa(x):k]}\Z \arrow[r] &0\\
0 \arrow[r] & \Pic\,X \arrow[r, "\pr_1^*"]\arrow[u, equal] & \Pic\,(X \times_k (\P^n_k)^*) \arrow[r, "\cdot \ell"]\arrow[u, "j^*"] & {[\kappa(x):k]}\Z \arrow[r]\arrow[u, equal] &0, 
\end{tikzcd}
\end{equation}
where $\ell$ is a line contained in a fibre $\P^{n-1}_{\kappa(x)}$ of the $\P^{n-1}$-bundle $X^{\univ}_{\varphi} \to X$ 
over a closed point $x \in X$  
(Remark \ref{r-tangent-bdl}). 
This fibre $\P^{n-1}_{\kappa(x)}$ is a hyperplane on 
the projective space $\{x \} \times_k (\P^n_k)^* (\simeq \P^n_{\kappa(x)})$, 
and hence the lower sequence in (\ref{e1-LHST-surje}) is exact. 
Since we have $X^{\univ}_{\varphi} \simeq \P_X(E)$ for some locally free sheaf $E$ of rank $n$ 
(Remark \ref{r-tangent-bdl}), 
also the upper sequence  in (\ref{e1-LHST-surje}) is exact. 
By the snake lemma, $j^*$ is an isomorphism. 
We obtain the following group homomorphisms:  
\[
\Pic\,X \times \Pic\,(\P^n_k)^* \xrightarrow{\pr_1^* \times \pr_2^*, \simeq}  \Pic\,(X \times_k (\P^n_k)^*) \xrightarrow{j^*, \simeq} \Pic\,X^{\univ}_{\varphi} \xrightarrow{\gamma^*} \Pic\,X^{\gen}_{\varphi}.  
\]
Note that $\gamma^*$ is surjective, because $X^{\univ}_{\varphi}$ is a regular variety and 
$X_{\varphi}^{\gen}$ is the generic fibre of $\pr_2 \circ j : X^{\univ}_{\varphi} \to (\P^n_k)^*$. 
Hence the induced composite group homomorphism 
\[
\Pic\,X \times  \Pic\,(\P^n_k)^* \to \Pic\,X^{\gen}_{\varphi}, \qquad (L, M) \mapsto  
\gamma^*j^*\pr_1^*L  \otimes \gamma^*j^*\pr_2^*M
\]
is surjective. 
Since $\pr_2 \circ j \circ \gamma : X_{\varphi}^{\gen} \to (\P^n_k)$ factors through $\Spec\,k(n)$ 
for a field $k(n)$, 
we get $\gamma^*j^*\pr_2^*M  \simeq (\pr_2 \circ j \circ \gamma)^*M \simeq \MO_{X^{\gen}_{\varphi}}$ 
for every $M \in \Pic\,(\P^n_k)^*$. 
Then $\beta^*:\Pic\,X  \to \Pic\,X^{\gen}_{\varphi}$ is surjective by  $\gamma^*j^*\pr_1^*L = \beta^*L$. 
\end{proof}

}


\begin{bibdiv}
\begin{biblist*}

\bib{AM69}{book}{
   author={Atiyah, M. F.},
   author={Macdonald, I. G.},
   title={Introduction to commutative algebra},
   publisher={Addison-Wesley Publishing Co., Reading, Mass.-London-Don
   Mills, Ont.},
   date={1969},
}


\bib{Bad01}{book}{
   author={B\u{a}descu, Lucian},
   title={Algebraic surfaces},
   series={Universitext},
   note={Translated from the 1981 Romanian original by Vladimir Ma\c{s}ek and
   revised by the author},
   publisher={Springer-Verlag, New York},
   date={2001},
   pages={xii+258},
}


\bib{BMP}{article}{
   author={Bhatt, Bhargav},
   author={Ma, Linquan},
   author={Patakfalvi, Zsolt},
   author={Schwede, Karl},
   author={Tucker, Kevin},
   author={Waldron, Joe},
   author={Witaszek, Jakub},
   title={Globally +-regular varieties and the minimal model program for threefolds in mixed characteristic},
   journal={preprint},
   eprint={arXiv:2012.15801v2},
}




\bib{CTW18}{article}{
   author={Cascini, Paolo},
   author={Tanaka, Hiromu},
   author={Witaszek, Jakub},
   title={Klt del Pezzo surfaces which are not globally $F$-split},
   journal={Int. Math. Res. Not. IMRN},
   date={2018},
   number={7},
   pages={2135--2155},
   issn={1073-7928},
}


\bib{CR11}{article}{
   author={Chatzistamatiou, Andre},
   author={R\"{u}lling, Kay},
   title={Higher direct images of the structure sheaf in positive
   characteristic},
   journal={Algebra Number Theory},
   volume={5},
   date={2011},
   number={6},
   pages={693--775},
   issn={1937-0652},
   review={\MR{2923726}},
   doi={10.2140/ant.2011.5.693},
}
\bib{CGM83}{article}{
   author={Cumino, Caterina},
   author={Greco, Silvio},
   author={Manaresi, Mirella},
   title={Bertini theorems for weak normality},
   journal={Compositio Math.},
   volume={48},
   date={1983},
   number={3},
   pages={351--362},
}

\bib{CGM86}{article}{
   author={Cumino, Caterina},
   author={Greco, Silvio},
   author={Manaresi, Mirella},
   title={An axiomatic approach to the second theorem of Bertini},
   journal={J. Algebra},
   volume={98},
   date={1986},
   number={1},
   pages={171--182},
   issn={0021-8693},
}

\bib{CGM89}{article}{
   author={Cumino, Caterina},
   author={Greco, Silvio},
   author={Manaresi, Mirella},
   title={Hyperplane sections of weakly normal varieties in positive
   characteristic},
   journal={Proc. Amer. Math. Soc.},
   volume={106},
   date={1989},
   number={1},
   pages={37--42},
   issn={0002-9939},
}

\bib{Das15}{article}{
   author={Das, Omprokash},
   title={On strongly $F$-regular inversion of adjunction},
   journal={J. Algebra},
   volume={434},
   date={2015},
   pages={207--226},
   issn={0021-8693},
   review={\MR{3342393}},
   doi={10.1016/j.jalgebra.2015.03.025},
}

\bib{DM}{article}{
   author={Datta, Rankeya},
   author={Murayama, Takumi},
   title={Permanence properties of $F$-injectivity},
   journal={preprint},
   eprint={arXiv:1906.11399v2},
}


\bib{GT16}{article}{
   author={Gongyo, Yoshinori},
   author={Takagi, Shunsuke},
   title={Surfaces of globally $F$-regular and $F$-split type},
   journal={Math. Ann.},
   volume={364},
   date={2016},
   number={3-4},
   pages={841--855},
}


\bib{GT80}{article}{
   author={Greco, S.},
   author={Traverso, C.},
   title={On seminormal schemes},
   journal={Compositio Math.},
   volume={40},
   date={1980},
   number={3},
   pages={325--365},
}
\bib{EGAIV2}{article}{
   author={Grothendieck, A.},
   title={\'{E}l\'{e}ments de g\'{e}om\'{e}trie alg\'{e}brique. IV. \'{E}tude locale des sch\'{e}mas et
   des morphismes de sch\'{e}mas. II},
   language={French},
   journal={Inst. Hautes \'{E}tudes Sci. Publ. Math.},
   number={24},
   date={1965},
   pages={231},
}

\bib{Har77}{book}{
   author={Hartshorne, Robin},
   title={Algebraic geometry},
   note={Graduate Texts in Mathematics, No. 52},
   publisher={Springer-Verlag, New York-Heidelberg},
   date={1977},
}


\bib{HH94}{article}{
   author={Hochster, Melvin},
   author={Huneke, Craig},
   title={$F$-regularity, test elements, and smooth base change},
   journal={Trans. Amer. Math. Soc.},
   volume={346},
   date={1994},
   number={1},
   pages={1--62},
   issn={0002-9947},
}
   

\bib{Jou83}{book}{
   author={Jouanolou, Jean-Pierre},
   title={Th\'{e}or\`emes de Bertini et applications},
   language={French},
   series={Progress in Mathematics},
   volume={42},
   publisher={Birkh\"{a}user Boston, Inc., Boston, MA},
   date={1983},
   pages={ii+127},
}



\bib{Kol13}{book}{
   author={Koll{\'a}r, J{\'a}nos},
   title={Singularities of the minimal model program},
   series={Cambridge Tracts in Mathematics},
   volume={200},
   note={With a collaboration of S\'andor Kov\'acs},
   publisher={Cambridge University Press, Cambridge},
   date={2013},
}





\bib{KM98}{book}{
   author={Koll{\'a}r, J{\'a}nos},
   author={Mori, Shigefumi},
   title={Birational geometry of algebraic varieties},
   series={Cambridge Tracts in Mathematics},
   volume={134},
   publisher={Cambridge University Press, Cambridge},
   date={1998},
}





\bib{Man80}{article}{
   author={Manaresi, Mirella},
   title={Some properties of weakly normal varieties},
   journal={Nagoya Math. J.},
   volume={77},
   date={1980},
   pages={61--74},
}

\bib{Mat89}{book}{
   author={Matsumura, Hideyuki},
   title={Commutative ring theory},
   series={Cambridge Studies in Advanced Mathematics},
   volume={8},
   edition={2},
   note={Translated from the Japanese by M. Reid},
   publisher={Cambridge University Press, Cambridge},
   date={1989},
   pages={xiv+320},
}

\bib{Mur21}{article}{
   author={Murayama, Takumi},
   title={The gamma construction and asymptotic invariants of line bundles over arbitrary fields},
   journal={Nagoya Math. J.},
   volume={242},
   date={2021},
   pages={165--207},
}







\bib{ST20}{article}{
   author={Sato, Kenta},
   author={Takagi, Shunsuke},
   title={General hyperplane sections of threefolds in positive
   characteristic},
   journal={J. Inst. Math. Jussieu},
   volume={19},
   date={2020},
   number={2},
   pages={647--661},
}

\bib{Sch07}{article}{
   author={Schwede, Karl},
   title={A simple characterization of Du Bois singularities},
   journal={Compos. Math.},
   volume={143},
   date={2007},
   number={4},
   pages={813--828},
   issn={0010-437X},
}

\bib{Sch09}{article}{
   author={Schwede, Karl},
   title={$F$-adjunction},
   journal={Algebra Number Theory},
   volume={3},
   date={2009},
   number={8},
   pages={907--950},
}

\bib{SS10}{article}{
   author={Schwede, Karl},
   author={Smith, Karen E.},
   title={Globally $F$-regular and log Fano varieties},
   journal={Adv. Math.},
   volume={224},
   date={2010},
   number={3},
   pages={863--894},
   issn={0001-8708},
}

\bib{SZ13}{article}{
   author={Schwede, Karl},
   author={Zhang, Wenliang},
   title={Bertini theorems for $F$-singularities},
   journal={Proc. Lond. Math. Soc. (3)},
   volume={107},
   date={2013},
   number={4},
   pages={851--874},
}


\bib{Sei50}{article}{
   author={Seidenberg, A.},
   title={The hyperplane sections of normal varieties},
   journal={Trans. Amer. Math. Soc.},
   volume={69},
   date={1950},
   pages={357--386},
   issn={0002-9947},
}










\bib{Tan17}{article}{
   author={Tanaka, Hiromu},
   title={Semiample perturbations for log canonical varieties over an
   $F$-finite field containing an infinite perfect field},
   journal={Internat. J. Math.},
   volume={28},
   date={2017},
   number={5},
   pages={1750030, 13},
}








\bib{Vel95}{article}{
   author={V\'{e}lez, Juan D.},
   title={Openness of the F-rational locus and smooth base change},
   journal={J. Algebra},
   volume={172},
   date={1995},
   number={2},
   pages={425--453},
   issn={0021-8693},
}


\end{biblist*}
\end{bibdiv}
\end{document}